\documentclass[11pt]{article}
\usepackage{amsfonts}
\usepackage{amsmath}
\usepackage{amsthm}
\usepackage{fullpage}
\usepackage{hyperref}
\usepackage{graphicx}
\usepackage{amssymb,verbatim}
\usepackage{pstricks}

\newtheorem{theorem}{Theorem}[section]
\newtheorem{proposition}[theorem]{Proposition}
\newtheorem{lemma}[theorem]{Lemma}
\newtheorem{remark}[theorem]{Remark}
\newtheorem{example}[theorem]{Example}
\newtheorem{problem}{Problem}

\newcommand{\RR}{\mathbb R}

\newtheorem{defin}{Definition}
\def\s0{{ s_0}}
\def\tnabla{{\widetilde\nabla}}
\def\ts0{{\tilde s_0}}

\def\eq#1{(\ref{#1})}
\def\nn{\nonumber}

\def\({\left(\begin{array}{cccccc}}
\def\){\end{array}\right)}

\def\bes{\begin{eqnarray}}
\def\ees{\end{eqnarray}}

\newcommand{\del}{\partial}
\newcommand{\pd}[2]{\frac{\partial{#1}}{\partial{#2} }}
\newcommand{\beq}{\begin{equation}}
\newcommand{\eeq}{\end{equation}}
\newcommand{\bea}{\begin{eqnarray}}
\newcommand{\eea}{\end{eqnarray}}
\newcommand{\baln}{\begin{align}}
\newcommand{\ealn}{\end{align}}
\newcommand{\beann}{\begin{eqnarray*}}
\newcommand{\eeann}{\end{eqnarray*}}

\newcommand{\lam}{\ensuremath{\lambda}}

\newcommand{\aaa}{\ensuremath{a}}

\newcommand{\R}{\ensuremath{\mathfrak{R}}}

\newcommand{\X}{\ensuremath{\mathcal{X}}}
\newcommand{\fs}{\ensuremath{\mathcal{F}}}
\newcommand{\fstriv}{\fs^{\text{triv}}}
\newcommand{\fsid}{\fs^{\text{id}}}
\newcommand{\oo}{\ensuremath{\mathcal{O}}}
\newcommand{\uu}{\ensuremath{\mathcal{U}}}
\newcommand{\ww}{\ensuremath{\mathcal{W}}}
\newcommand{\vr}{\ensuremath{\mathbf{r}}}
\newcommand{\vs}{\ensuremath{\mathbf{s}}}
\newcommand{\vt}{\ensuremath{\mathbf{t}}}
\newcommand{\vf}{\ensuremath{\mathbf{f}}}
\newcommand{\tJ}{\ensuremath{\tilde{J}}}
\newcommand{\cl}{\ensuremath{\mathcal{L}}}

\newcommand{\pf}{\begin{proof}}
\newcommand{\foorp}{\end{proof}}

\DeclareMathOperator{\rank}{rank}

\DeclareMathOperator{\Span}{span}

\newcommand{\Addresses}{{
  \bigskip
  \footnotesize

  M.~Benfield, \textsc{San Diego, CA 
} \texttt{(mike.benfield@gmail.com).}

  \medskip

 H.~K.~Jenssen, \textsc{Department of
Mathematics, Penn State University}
 \texttt{(jenssen@math.psu.edu).}

  \medskip

  I.~A.~Kogan, \textsc{Department of Mathematics, North Carolina State University}\texttt{(iakogan@ncsu.edu).}

}}
\makeatletter
\let\@fnsymbol\@arabic
\makeatother

\title{Jacobians with prescribed eigenvectors}

\author{Michael Benfield\thanks{M.~Benfield was partially supported by NSF grant
DMS-1311743.} \and Helge Kristian Jenssen\thanks{H.~K.~Jenssen was partially supported by NSF grant
DMS-1311353.} \and Irina A. Kogan\thanks{I.~A.~Kogan was partially supported by NSF grant
DMS-1311743.}}

\date{}

\begin{document}
\maketitle

\abstract{Let $\Omega\subset \RR^n$ be open and let $\R$ be a partial frame on $\Omega$, that is a set of $m$ linearly 
independent vector fields  prescribed on $\Omega$ ($m\leq n$). We consider 
the issue of describing the set of  all maps $F:\Omega\to\RR^n$ with the property that each of the given
 vector fields is an eigenvector of the Jacobian matrix of $F$. 
By introducing a coordinate independent definition of the Jacobian, we 
obtain an intrinsic formulation of the problem, which leads to an overdetermined PDE 
system, whose compatibility conditions can be expressed  in an intrinsic, coordinate independent manner.
To analyze this system we formulate and prove a generalization of the classical 
Frobenius integrability theorems. 
The size and structure of the solution set of this system depends on the properties of the partial frame, in particular, whether or not it is in  involution.
A particularly nice subclass of  involutive partial frames, called rich, can be completely analyzed.  Involutive, but non-rich case is somewhat harder to handle.  We provide a  complete answer in the case of $m=3$ and arbitrary $n$, as well as some general results  for arbitrary $m$.  
 The 
non-involutive case is far more challenging, and we only obtain a comprehensive analysis in the case $n=3$, $m=2$.
Finally, we provide explicit examples illustrating the various possibilities.  Our initial motivation for considering this problem comes from  the geometric study of hyperbolic conservative systems  in one spatial dimension.} 



\tableofcontents
\section{Introduction}\label{sect-intro}
The present work deals with the construction of maps $F:\Omega\subset \RR^n\to\RR^n$ 
whose Jacobian matrix has a partially prescribed set of eigenvector fields on $\Omega$. 
We consider this problem locally, i.e.\ in a sufficiently small neighborhood of a given point 
in $\Omega$. The case when the full frame of $n$ independent eigenvectors is prescribed has been considered 
in \cite{jk1}.  The  generalization to a partially prescribed set of eigenvector fields allows a  greater degree of flexibility in constructing such maps $F$  and, in particular, 
 permits maps $F$ whose Jacobian matrix is not diagonalizable. Another difference from the pervious work is that
all the overdetermined system of PDEs arising in the current paper are analyzed using smooth\footnote{We employ $C^1$ integrability theorems, but to avoid some technicalities $C^\infty$-smoothness is assumed throughout the paper.}   integrability theorems and, in particular, a recently proved
 generalization of the Frobenius  theorem (see Section~\ref{sect-integrthm}). This theorem allows us  to remain in
the smooth category,  while in \cite{jk1} we 
 appealed in some cases to the Cartan-K\"ahler theorem, which 
requires analyticity assumptions.

Our motivation stems from the study of initial value problems for one dimensional 
conservative systems of the form 
\beq\label{claw}
	u_t+F(u)_x=0,\qquad u(0,x)=u_0(x),
\eeq
where $t\in\RR$ and $x\in\RR$ are the independent variables, $u=u(t,x)\in \RR^n$
is a vector of unknowns, and the flux function $F$ is defined on some open set in 
$\RR^n$ and takes values in $\RR^n$. It is an outstanding open problem to provide 
an existence theory for the Cauchy problem for \eq{claw} which is general enough 
to cover nonlinear systems of physical interest and  
of initial data $u_0(x)$ of ``large'' total variation. Such a theory is in place for near-equilibrium 
solutions (Glimm's theorem \cite{glimm}):  global-in-time existence of a weak solution
is guaranteed, provided the initial data $u_0(x)$ have sufficiently small total variation. 
(For detailed accounts of this theory see \cite{smoller,bressan,daf}.) 
A key ingredient in the proof is the use of Riemann problems, i.e.\ initial value 
problems for \eq{claw} where the data $u_0(x)$ consists of two constant 
states $u^\pm$, separated by a jump discontinuity,
\beq\label{claw-data}
u_0(x)=
\left\{\begin{array}{ll}
	u^- & x<0\\
	u^+ & x>0.
\end{array}\right.
\eeq
By knowing how to solve Riemann problems one can, via an approximation
scheme, solve general (small variation) Cauchy problems. 

The solution of \eq{claw}-\eq{claw-data} is a self-similar (function of $x/t$) 
fan of $n$ waves emanating from the origin. These waves are determined
from the eigenvalues and eigenvectors of the Jacobian matrix $[D_uF]$.
It is, therefore, of interest to gain an understanding of how the eigenstructure 
of $F$ induce properties of the solutions $u(t,x)$. It is then a basic question 
to what extent one can prescribe some or all of the eigenvectors 
of the flux $F$.

The present work is concerned with this last, purely geometric problem. 
 A precise 
formulation of the problem is provided in Section~\ref{sect-problem}. This first formulation, ``Problem~\ref{pr},'' 
makes use of a chosen coordinate system. Section~\ref{sect-geom} provides the geometric framework
required to obtain a  coordinate-free formulation. We also state and prove of a generalization of the Frobenius integrability theorem,
 which we use in this paper.
In Section~\ref{sect-intrinsic}, we give an intrinsic (coordinate independent) 
definition of the Jacobian, and use it to reformulate Problem~\ref{pr} in an intrinsic manner (see Problem~\ref{pr2}).
Exploiting the coordinate independent formulation we treat, in Section~\ref{sect-involutive},  the case when 
the prescribed, partial frame of eigenvectors-to-be, is in involution. In this  case,  the integrability condition  of the $\fs(\R)$ system lead to a closed algebro-differential system on eigenvalues-to-be $\lam$'s.
Section~\ref{sect-noninvol} analyzes the simplest non-involutive case of two prescribed vector fields 
in $\RR^3$. Finally, Section~\ref{sect-ex} provides a list of examples that illustrate the 
results from the earlier sections.
\section{Problem formulation}\label{sect-problem}
 In this paper,  $[D_u\Psi]$ denotes the Jacobian matrix of a map  $\Psi$ from an 
open subset $\Omega\subset \RR^{n}$ to $\RR^{n}$, relative to coordinates 
$u$. That is,
\[[D_u\Psi]=\left[\pd{\Psi^i}{u^j}\right]_{i,j=1,\dots,n}.\]
We use the notation  $[D_u\Psi]|_{u=\bar u}$, or simply  $[D_u\Psi]|_{\bar u}$, 
when the matrix is evaluated at a point $\bar u$.
 We consider the following problem:

\begin{problem}\label{pr}
	Given an open subset $\Omega\subset \RR^n$ 
	on which we fix a coordinate system $u=(u^1,\dots,u^n)$ and a point 
	$\bar u\in \Omega$. Let $\R=\{R_1, \dots, R_m\}$ be a set of $m\leq n$ 
	smooth vector valued functions $R_i: \Omega\to \RR^n$ which are 
	linearly independent at $\bar u$. Then: describe the set $\fs(\R)$ of 
	all smooth vector-valued functions 
	\[F(u)=[F^1(u),\dots, F^n(u)]^T\]
	defined near $\bar u$ and with the property that 
	$R_1(u), \dots, R_m(u)$ are right eigenvectors of the 
	Jacobian matrix $[D_uF]|_u$ throughout 
	a neighborhood of $\bar u$. In other words, we ask that there 
	exist smooth, scalar functions $\lam^i$ such that 
	\begin{equation} \label{eq:eigenvectors}
  		\left[D_uF\right] R_i(u) = \lambda^i(u)  R_i(u),\qquad i=1,\dots,m,
  	\end{equation}
	holds on a neighborhood of $\bar u$.
\end{problem}
As outlined in the Introduction, we are motivated by the construction
of flux functions $F$ in systems of conservation laws of the form \eq{claw}.
The system \eq{claw} is called \emph{hyperbolic} on $\Omega$ provided 
the Jacobian matrix $[D_uF]$ has a basis of real eigenvectors at each 
$u\in\Omega$, and it is called \emph{strictly hyperbolic} if, in addition, 
all its eignvalues are distinct at each $u\in\Omega$.  
We adopt the term \emph{flux} for a vector-function satisfying 
\eqref{eq:eigenvectors}, with adjectives \emph{hyperbolic, strictly 
hyperbolic or non-hyperbolic}  depending on the structure of 
eigenvectors and eigenvalues of  $[D_uF]$, as described above.

In the list below, we clarify what we mean by ``describe'' in Problem~\ref{pr} 
and  make some preliminary observations about Problem \ref{pr}:
\begin{enumerate}
\item (\emph {PDE system}) Equations \eqref{eq:eigenvectors} comprise a system of $m\,n$ first order PDEs on $n+m$ unknown functions $\lam^i$ and $F^j$:
\begin{equation} \label{eq-fr-coord-e}
  		\sum R_i^k \,\pd{F^j}{u^k} = \lambda^i \,R_i^j, \quad \text{ for } i=1,\dots,m,\, j=1,\dots, n,
  		\end{equation}
where $R_i(u)=[R_i^1(u),\dots, R_i^n(u)]^T$, $i=1,\dots,m$.  This system is overdetermined for all $n\geq m$, such that $n>2$ and $m\geq 2$. Although derivatives of $\lam$'s do not appear in the equations, these functions are not arbitrary parameters, but must, in turn, satisfy certain differential equations arising as differential consequence of \eqref{eq-fr-coord-e}.

\item (\emph {Vector space structure})   Let $F_1, F_2\in \fs(\R)$, have  the domains of definitions $\Omega_1$ and $\Omega_2$, respectively.  Since $\bar u$ belongs to both  $\Omega_1$ and $\Omega_2$, the sum $F_1+F_2$ is defined on the non-empty open neighborhood $\Omega_1\cap\Omega_2$  of $\bar u$. It is easy to check that  $F_1+F_2$ belongs  $\fs(\R)$ and that $\bar a\,F_1\in \fs$, where $\bar a$ is any real number, belongs to $\fs(\R)$. Thus $\fs(\R)$ is a \emph{vector space} over $\RR$.  We will see below, that in some instances  this is a finite dimensional vector space, while in others it is an infinite dimensional space. In the latter case, we describe the ``size''  of $\fs(\R)$ in terms of the number of arbitrary functions of a certain number of variables appearing in the general solution of \eqref{eq-fr-coord-e}. These arbitrary functions prescribe the values of $F$  and $\lam$'s along certain submanifolds of $\Omega$.
To obtain these results we use the integrability theorem stated in Section~\ref{sect-integrthm}.

\item(\emph {Scaling invariance}) Since eigenvectors are defined up to \emph{scaling}, it is clear that 
\beq\label{scaling}
\fs(R_1, \ldots, R_m)=\fs(\alpha^1\,R_1, \ldots, \alpha^m R_m)
\eeq
 for any nowhere zero smooth functions $\alpha^i  \text{ on } \Omega$.

\item (\emph{Trivial solutions}) For  any $n+1$ constants  $\bar\lam, \bar a^1, \dots, \bar a^n  \in\RR$ the ``trivial''  flux 
  \begin{equation} \label{tflux}
  F(u)=\bar \lam
\left[
\begin{array}{c}
 u^1   \\
  \vdots   \\
  u^n  
\end{array}
\right]
+\left[\begin{array}{c}
\bar a^1   \\
  \vdots   \\
 \bar a^n  
\end{array}\right]
\end{equation}
 satisfies  \eqref{eq:eigenvectors}.
The set of such trivial solutions, denoted by $\fstriv$, is an (n+1)-dimensional vector subspace of $\fs(\R)$.

 \item (\emph{Triviality is generic}) It is worthwhile  emphasizing that   when $n>2$ and $m\geq 2$, the compatibility conditions for $\fs(\R)$-system are closed, and thus almost all frames admit only trivial fluxes.
 One of the goals of the paper is to determine the properties of the frames that allow them to possess  non-trivial fluxes, and in particular strictly hyperbolic fluxes.
 
 \end{enumerate}

The vector space $\fs(\R)$ will be called the \emph{flux space}.
We are, of course, only interested in non-trivial fluxes, and particularly 
in strictly hyperbolic fluxes due to their central role in the 
theory of conservation laws.


The next remark addresses the coordinate dependence of our formulation of Problem~\ref{pr}. 
In Section~\ref{sect-fr}, we formulate  a coordinate independent version  (Problem~\ref{pr2}), 
which, when expressed in \emph{an affine  system of coordinates} 
(see Definition~\ref{def-affine}) coincides with Problem~\ref{pr}.  
This intrinsic definition allows us to apply a geometric approach to analyze 
the solution set of PDE system \eqref{eq-fr-coord-e}.

\begin{remark}[Coordinate dependence of the problem formulation]
Assume $F(u)\in \fs(\R)$  for  {$\R=\{R_1,\dots, R_m\}$}, i.e.\ there exist $\lam^1(u), \dots, \lam^m(u)$, such that system  \eqref{eq:eigenvectors} is satisfied.
 Let a change of variables be described by  a local diffeomorphism $$u=\Phi(w).$$ 
 It is then not true, in general, that  $\tilde F(w)=F(\Phi(w))$  belongs to $\fs(\tilde \R)$, where $\tilde \R=\{\tilde R_1(w),\dots, \tilde R_m(w)\}$, with 
$\tilde R_i(w)= R_i(\Phi(w))$. Indeed:
\[[D_w\left (F\circ \Phi\right)]\, \tilde R_i= [D_uF]|_{u=\Phi(w)}
[D_w\Phi] \, R_i(\Phi(w)).\]
In general,  $R_i(\Phi(w))$ is not an eigenvector of $ [D_uF]|_{u=\Phi(w)}
[D_w\Phi]$.

Even if we transform  $R_i(u)$'s, by treating them, more appropriately, as vector-fields:
\[ R_i^*(w)= [D_w\Phi]^{-1}\, R_i(\Phi(w)),\]
 then
 \begin{align*}
	[D_w\left (F\circ \Phi\right)]\,  R^*_i&= [D_uF]|_{u=\Phi(w)}
	[D_w\Phi] [D_w\Phi]^{-1}\, R_i(\Phi(w))\\
 	&= [D_uF]|_{u=\Phi(w)} R_i(\Phi(w))\\
 	&
  	= \lam^i(\Phi(w))R_i(\Phi(w))=\lam^i(\Phi(w))\,[D_w\Phi]R^*_i(w),
\end{align*}
and we see that $R^*_i(w)$ is not an eigenvector of $[D_w\left (F\circ \Phi\right)]$, unless it is an eigenvector of $[D_w\Phi]$.
\end{remark}
We also recall that the property of  a matrix being  the Jacobian 
matrix of some map is also coordinate dependent:
\begin{remark}[Coordinate dependence of 
	the  property of a matrix being a Jacobian matrix]\label{rem-coord-jac}
	Assume $A(u)=[D_uF]$ for some smooth map 
	$F:\Omega\to\RR^n$, and let a change of coordinates be 
	given by a diffeomorphism $u=\Phi(w)$. Then, it is not necessarily 
	the case that matrix $A(\Phi(w))$  is a Jacobian matrix of any 
	map in $w$-coordinates.
\end{remark}
On the other hand, it is still possible to give a coordinate independent definition 
of the Jacobian linear map, as we do in Section~\ref{sect-intrinsic-jac}. This is 
used to obtain a coordinate-independent formulation of Problem~\ref{pr}.
We exploit this by working in frames that are adapted to the problem at hand, and we  
use the following geometric preliminaries.

\section{Geometric preliminaries}\label{sect-geom}
Most of the notions and results, reviewed in this section, can be found in a standard differential geometry text-book. We included them to set up notation, as well as to make a paper self-contained.   The notable exception is Section~\ref{sect-integrthm}, where we state  and prove a generalization of the Frobenius integrability theorem. 
\subsection{Vector fields, flows, partial frames, involutivity, richness}
\label{sect-vf}

It will be useful for us to give an intrinsic, coordinate free definition of a vector field as a linear first order differential operator on the set of functions
\begin{defin}
A smooth vector field $\vr$ on $\Omega$ is an $\RR$-linear  map from the set of smooth functions $C^\infty(\Omega)$  to itself that satisfying the product rule.

\end{defin}  

The set of all smooth vector  fields will be denoted as   $\mathcal X(\Omega)$, and it is an  infinite dimensional vector space over $\RR$ and a free $n$-dimensional module over $C^\infty(\Omega)$.
Relative to any coordinate system, a vector field $\vr$ evaluated at a point $\bar u\in \Omega$ becomes a vector in $\RR^n$. We say that  vector fields $\vr_1,\dots, \vr_m$ are\emph{ independent } at $\bar u\in \Omega$ if vectors $\vr_1(\bar u),\dots, \vr_m(\bar u)\in \RR^n$ are independent over $\RR$ relative to one and, therefore, to all coordinate systems. 
\begin{defin} A set of smooth vector fields $\vr_1,\dots,\vr_m$ on $\Omega$ is called a \emph{partial frame} on $\Omega$ if they are  independent for all $\bar u\in \Omega$. If $m=n$, then this set is called a \emph{frame}.
\end{defin}
It is easy to show that a frame comprises a basis of the module $\mathcal X(\Omega)$  over $C^\infty(\Omega)$, i.~e. for any smooth vector field $\vr\in \mathcal X(\Omega)$ there are smooth functions $R^1,\dots, R^n\in C^\infty(\Omega)$, called the components of $\vr$ relative to frame   $\vr_1,\dots,\vr_n$, such that 
$\vr=R^1\,\vr_1+\dots+R^n\,\vr_n$. For a fixed coordinate system $u^1,\dots, u^n$, the frame $\left\{\pd{}{u^1},\dots, \pd{}{u^n}\right\}$ of partial derivatives  is called \emph{a coordinate frame}, but as we see below using non-coordinate frames  can simplify a problem.
The Lie bracket of two vector-fields can be defined as the commutator operator on functions.
\begin{defin}
Given two smooth vector field, their Lie bracket is the map $C^\infty(\Omega)\to C^\infty(\Omega)$ defined by
$$[\vr_1,\vr_2]\phi=\vr_1(\vr_2(\phi))-\vr_2(\vr_1(\phi)).$$
\end{defin}
A standard calculation shows that $[\vr_1,\vr_2]$ is a vector field,~i.e. a \emph{first order} linear differential operator.  
 Skew symmetry of the  Lie bracket is obvious and Jacobi identity can be checked by an explicit calculation. Therefore, $\mathcal X(\Omega)$ has a structure of an infinite-dimensional real Lie algebra. 
Given a frame $\vr_1,\dots,\vr_n$, we can write the following \emph{structure equations}:
$$[\vr_i,\vr_j]=\sum_{k=1}^n c_{ij}^k \vr_k,$$
where $c_{ij}^k$, such that  $c_{ij}^k=-c_{ji}^k$, are smooth functions on $\Omega$, called \emph{structure coefficients}, or \emph{structure functions}. In the conservation laws literature, these functions are called \emph{interaction coefficients} because of their role in wave interaction formulas \cite{glimm}. The Jacobi identity imply the following relationship on the structure coefficients:
\begin{align}\label{c-jacobi}
	& r_l\big(c_{jk}^i\big)+r_k\big(c_{lj}^i\big)+r_j\big(c_{kl}^i\big)
	\qquad\qquad\qquad\qquad\text{(Jacobi)}\nn\\
	&\quad+\sum_{s=1}^n\big( c_{jk}^sc_{ls}^i+c_{lj}^sc_{ks}^i+c_{kl}^sc_{js}^i\big)=0
	\qquad 1\leq i, j, k,l \leq n.
\end{align}

We will define two classes of partial frames with especially nice properties:
\begin{defin}[Involutive frame]\label{def-inv}
We say that  a  partial frame $\R=\{\vr_1,\dots,\vr_m\}$ is \emph{in involution} if $[\vr_i,\vr_j]\in\Span_{C^{\infty}(\Omega)} \R$ for all $\vr_i,\vr_j\in\R$.\end{defin}
The  proof of the following proposition can be found  in the proof of Theorem~6.5 of Spivak \cite{spivak:v1}. 
\begin{proposition}\label{prop-inv-distr} Let  $\vr_1,\dots,\vr_m$ be a partial frame in involution on $\Omega$, then there is a commutative partial frame such that   $\tilde\vr_1,\dots,\tilde\vr_m$ on some open $\Omega'\subset\Omega$, such that 
$$\Span_{\RR}\{\vr_1|_u,\dots,\vr_m|_u\}=\Span_{\RR}\{\tilde\vr_1|_u,\dots,\tilde\vr_m|_u\} \text{ for all } u\in\Omega'.$$
\end{proposition}

\begin{proposition}\label{prop-comm-frame}\emph{(Theorem 5.14 in \cite{spivak:v1})} If $\vr_1,\dots,\vr_m$ is a commutative partial frame on $\Omega$, then in a neighborhood of each point $\bar u\in\Omega$ there  exist coordinate functions $v^1,\dots, v^n$, such that 
$$\vr_i=\pd{}{v^i},\quad i=1,\dots, m.$$
\end{proposition}
\begin{defin}[Rich frame]\label{def-rich}
We say that  a  partial frame $\R=\{\vr_1,\dots,\vr_m\}$ is \emph{rich} if every pair of its vector fields is in involution, i.e. $[\vr_i,\vr_j]\in\Span_{C^{\infty}(\Omega)}\{\vr_i,\vr_j\}$ for all $i,j=1,\dots, m$.
\end{defin}

In  Lemma~\ref{lem-rich-coord},  we show that every rich partial frame $\vr_1,\dots,\vr_m$ can be scaled to become a commutative frame and so around each point one can find  coordinates $w^1,\dots, w^n$, and non-zero functions $\alpha^1, \dots, \alpha^n$ such that 
$$\alpha^i  \vr_i=\,\pd{}{w^i},\quad i=1,\dots, m.$$
Classically, a conservative system is  called rich if there are coordinate functions, called {\em Riemann invariants}, in which the system  is diagonalizable.  For definitions, and the 
fact that richness of a conservative system is equivalent to the richness of its eigenframe in the sense of our definition, we refer to \cite{serre2}, and Section 7.3 in \cite{daf}. 
Riemann invariants are exactly the coordinates  appearing in Lemma~\ref{lem-rich-coord}, in the case of full frame: $n=m$.
The term \emph{rich}  refers to a large family of extensions (companion conservation laws) that  strictly hyperbolic diagonalizable systems possess \cite{daf,serre2}.


\subsection{Connection, symmetry, flatness, affine coordinates}
We defined vector fields as directional derivatives of smooth functions. More generally, one can define directional derivatives of vector fields themselves
by introducing the notion of a covariant derivative. We will use this notion to give a coordinate free definition of Jacobians and to express Problem~\ref{pr} in a non-coordinate  frame, which make it easier to find its solution.

\begin{defin}A \emph{connection} $\nabla$ on $\Omega$ is an 
$\RR$-bilinear map 
\[\nabla\colon \mathcal{X}(\Omega)\times \mathcal{X}(\Omega)\to\mathcal{X}(\Omega)\qquad\qquad
	(\vr,\vs)\mapsto \nabla_\vr \vs\]
such that for any smooth function $\phi$ on  $\Omega$
\beq\label{connection}
	\nabla_{\phi\, \vr}\vs=\phi\nabla_\vr\vs, \qquad \nabla_\vr (\phi\, \vs)=\vr(\phi)\,\vs+\phi\nabla_\vr\vs\,.
\eeq
The vector field $ \nabla_\vr \vs$ is called the covariant derivative of $\vs$ in the direction of $\vr$.
\end{defin}
Given a  connection $\nabla$ and a  frame, $\{\vr_1,\dots,\vr_n\}$, for all $i,j\in\{1,\dots,n\}$ we can write 
\beq\label{eq-gam}\nabla_{\vr_i}\vr_j=\sum_{k=1}^n\Gamma^k_{ij}\vr_k,\eeq
for some  smooth functions $\Gamma^k_{ij}$,  called \emph{connection  components, or Christoffel symbols}. 
Conversely,  due to $\RR$-bilinearity and \eq{connection}, for any choice of a frame and $n^3$ functions 
$\Gamma^k_{ij},\, i,j,k=1,\dots,n$, formula \eqref{eq-gam}  uniquely defines a connection on $\Omega$.
\begin{defin}[affine coordinates]\label{def-affine} Given a connection $\nabla$, coordinate systems, such that relative to the corresponding coordinate frame  all Christoffel symbols for $\nabla$ are zero, are called   \emph{affine}.
\end{defin}
\begin{defin}[symmetry and flatness] A connection  $\nabla$  is \emph{symmetric}
if  for all $\vr,\vs\in \mathcal{X}(\Omega)$:
\beq\label{T0} \nabla _\vr \vs-\nabla_\vs \vr=[\vr,\vs]. \eeq

 A connection  $\nabla$  is \emph{flat}
if  for all $\vr,\vs, \vt\in \mathcal{X}(\Omega)$:
\beq
                           \label{R0}\nabla _\vr \nabla_\vs\, \vt-\nabla_\vs\nabla_\vr \vt=\nabla_{[\vr,\vs]}.
\eeq
\end{defin}
The above conditions are equivalent to the following relationships among the structure functions and   Christoffel symbols relative to  an \emph{arbitrary} frame:  for all  $i,j,k,s=1,\dots, n$,
\begin{eqnarray}\label{T0-gamma} \Gamma_{ij}^k-\Gamma_{ji}^k=c_{ij}^k &\qquad&\mbox{Symmetry}\\
                           \label{R0-gamma} r_s\big(\Gamma^j_{ki}\big)-r_k\big(\Gamma^j_{si}\big)
	=\sum_{l=1}^n \big(\Gamma^j_{kl}\Gamma^l_{si}-\Gamma^j_{sl}
	\Gamma^l_{ki}-c^l_{ks}\Gamma^j_{li}\big)&\qquad& \mbox{Flatness}.
\end{eqnarray}
A  well known  result, stated, for instance, in Proposition~1.1 in \cite{hirohiko}, implies that any flat and symmetric connection on $\Omega$ admits an affine system of coordinates, and that any two affine coordinate systems are related by an affine transformation:
\begin{proposition}\label{prop-flat-symm} A connection $\nabla$ on an $n$-dimensional manifold $M$ is symmetric and flat  (has properties \eq{T0} and \eq{R0}) if and only if, it can be covered  with an atlas of affine coordinate systems.
 
Two coordinate system $u=(u^1,\dots, u^n)$ and $w=(w^1,\dots, w^n)$, on an open subset $\Omega\subset M$, are affine if and only if  $[w^1,\dots, w^n]^T=C\,[u^1,\dots, u^n]^T+\bar b$, where $u$ and $w$ are treated as column vectors, $C\in\RR^{n\times n}$ is an $n\times n$ invertible   matrix and $\bar b\in \RR^n$ is a constant vector.  
\end{proposition}

Throughout the paper, we  will use a particular connection, denoted $ \tnabla$, defined by  setting all Christoffel symbols to be zero, relative to the coordinate frame corresponding  to  the coordinate system $u^1,\dots,u^n$ fixed in Problem~\ref{pr}:
\beq\label{eq-affine} \tnabla_{\pd{}{u^i}}{\pd{}{u^j}}=0,\, \text{ for all } i,j=1,\dots, n.\eeq

If the column vectors $R$ and $S$ are components of vector fields $\vr$ and $\vs$, respectively,   in  an \emph{affine} coordinate system, then 
\beq\label{eq-tr-conn}\text{ the components of }  \tilde\nabla_\vr \vs \text{ are given by the column vector $\vr(S)$},\eeq
where $\vr$ is applied to each component of $S$.

\subsection{Integrability  theorems} 
\label{sect-integrthm}

To analyze the ``size'' of the flux-space $\fs(\R)$ in Problem~\ref{pr}, we use  two integrability theorems: generalized Frobenius Theorem and Darboux Theorem. 
  
The classical Frobenius theorem has three equivalent formulations: PDE formulations, vector-field formulation, and differential form formulation (see Spivak \cite{spivak:v1} Theorems 6.1, 6.5 and 7.14  Warner \cite{warner} Theorem 1.60, Remark 1.61, Theorem 2.32). 
For our generalization, we start with a vector-field formulation, Theorem~\ref{th-vf-frob}, and then, as a consequence, prove its PDE formulation, Theorem~\ref{th-pde-frob}, which we  use further in the paper. When $n=m$, both theorems are equivalent to the corresponding local versions of   classical Frobenius theorem. A formulation of an appropriate  global foliation version of the  Theorem~\ref{th-vf-frob}, as well as its differential form formulation are of interest, but fall outside of the scope of this paper. 

In his thesis \cite{mike:thesis}, the first author proved  Theorem~\ref{th-pde-frob}  directly,  using contractive maps and Picard type argument. In the current proof, Picard type argument is hidden in the existence and uniqueness result for a flow of vectors field.  A weaker version of Theorem~\ref{th-pde-frob}  (with right hand-sides of \eqref{eq-frob} independent of $\phi$'s) appears in Lee \cite{lee:m}, Theorem~19.27.

The vector field generalization of the Frobenius Theorem, which we rigorously formulate below,  states that, given a local partial frame in involution  $\vs_1,\dots,\vs_m$ on an open subset  $\mathcal O$ of $\RR^{n+p}$, where  integers $1\leq m\leq n$ and $p\geq 1$,  and an  $(n-m)$-dimensional  embedded submanifold $\Lambda\subset\mathcal{ O}$, not tangent to any of the given vector fields $\vs$'s,  one can locally extend $\Lambda$ to an $n$-dimensional submanifold $\Gamma$,  tangent to each of  the $\vs$'s at every point.  Moreover, such an extension is locally unique.  If $n=m$, then $\Lambda$ is a single point and we get a statement which is equivalent to a local vector-field version of the classical Frobenius theorem.
\begin{theorem}\label{th-vf-frob} Let $\vs_1,\dots,\vs_m$ be a partial frame in involution  defined on an open subset $\mathcal O\subset\RR^{n+p}$, where $1\leq m\leq n$ and $p\geq 1$.  Let  $\Lambda\subset \mathcal O$ be an $(n-m)$-dimensional embedded submanifold, such that
\beq\label{eq-span}\Span_{\RR}\{\vs_1|_z, \dots, \vs_m|_z\}\oplus T_{z}\Lambda\cong\RR^{n}\eeq  for every point $z\in \Lambda$. Then for every point $\bar z\in \Lambda$, there exists an open neighborhood $\mathcal O_{\bar z}\subset \mathcal O$  and an  $n$-dimensional submanifold $\Gamma_{\bar z}$ of $\RR^{n+p}$, such that
\begin{enumerate}
\item[ 1)] $\Lambda\cap\oo_{\bar{z}}= \Lambda\cap\Gamma_{\bar{z}} $;
\item[2)] $\vs_i|_z\in T_z\Gamma_{\bar z}$,  for all $i=1,\dots, m$ and  for every point $z\in \Gamma_{\bar z}$. 
\end{enumerate}
 Manifold $\Gamma_{\bar z}$ is locally unique, i.e. if there is another $n$-dimensional manifold $\Gamma_{\bar z}'$, satisfying the two  conditions stated above, then $\Gamma_{\bar z}\cap \Gamma_{\bar z}' $ is also an $n$-dimensional manifold satisfying these conditions. 
\end{theorem}
\begin{proof}  Since $\vs_1,\dots,\vs_m$ are in involution then,  on an open neighborhood $\mathcal O'_{\bar z}\subset \RR^{n+p}$  of $\bar z \in \Lambda$, by Proposition~\ref{prop-inv-distr}, there exists a partial \emph{commutative} frame $\tilde \vs_1,\dots,\tilde\vs_m$, such that $\Span _{\RR}\{\vs_1|_{z},\dots,\vs_m|_{z}\} =\Span _{\RR}\{\tilde \vs_1|_z,\dots,\tilde \vs_m|_{ z}\}$ for all $z\in \oo'_{\bar z}$.
Since $\Lambda$ is an embedded submanifold, by shrinking $\mathcal O'_{\bar z}$ we may assume that $\Lambda'_{\bar z}=\Lambda\cap\oo'_{\bar z}$ is a coordinate neighborhood of $\bar z$ in $\Lambda$, i.e.~there exists an open set  $\ww'\subset \RR^{n-m}$  and a diffeomorphism $$\psi\colon \ww'\to \Lambda'_{\bar z},$$ such that $\psi(0)=\bar z$. Let $ {\mathcal B}_{\varepsilon}=(-\varepsilon,\varepsilon)^m$ denote an open box in $\RR^m$ with sides $2\varepsilon$ centered at the origin.
 Let $\exp^{\epsilon \vs}(z)$ denote the flow of a vector field  $\vs$ on $\mathcal O$, i.e.
\beq\label{eq-flow} \frac{d}{d\epsilon}\exp^{\epsilon \vs} (z)=\vr\big|_{\exp^{\epsilon \vs(z)}},\qquad \qquad \exp^{0\cdot \vs}(z)=z.\eeq
Then there exists an $\varepsilon>0$, such that the map $\Psi\colon \ww' \times {\mathcal B}_{\varepsilon} \to \RR^{n+p}$:
\beq
\label{eq-extend}\Psi(w, \epsilon)=\left(\exp^{\epsilon_1\,\tilde\vs_1} \circ
{\cdots}\circ\exp^{\epsilon_m\,\tilde\vs_m}\right)
\, \left(\psi(w)\right).\eeq
is defined for all $w\in \ww'$ and $\epsilon \in {\mathcal B}_{\varepsilon}$
The map  $\Psi$ is smooth (see \cite{lang93}, pp.~371-- 379).
Let $D\Psi|_0\colon T|_{0}\RR^n \to T|_{\bar z}\RR^{n+p}$  denote  the differential of $\Psi$ at the origin  in $\RR^n$  and let $D\psi|_0\colon T|_{0}\RR^{n-m} \to T|_{\bar z}\RR^{n+p}$  denote  the differential of $\psi$ at the origin in $\RR^{n-m}$.   Then vectors
$$D\Psi|_{0}\left(\pd{}{w^i}\right)=\left.\pd{}{w^i}\right|_{(w,\epsilon)=0}\,\Psi=\left.\pd{}{w^i}\right|_{w=0}\psi=D\psi|_0\left(\pd{}{w^i}\right), \text{ where  } i=1,\dots, {n-m},$$
span the tangent space $T|_{\bar z}\Lambda$. On the other hand,
$$D\Psi|_0\left(\pd{}{\epsilon^j}\right)=\left.\pd{}{\epsilon^j}\right|_{\epsilon^j=0}\,\exp^{\epsilon_j\,\tilde\vs_j} (\bar z)=\tilde \vs_j|_{\bar z},  \text{ for  } j=1,\dots {m}.$$
Therefore, due to \eqref{eq-span}, $D\Psi|_0$ has maximal rank $n$ at $0\in \RR^n$. Then, there exists  an open subset of $\uu\subset \ww'\times{\mathcal B}_{\varepsilon}\subset \RR^n$ containing the origin, such that the restriction $\Psi|_{\uu}\colon \uu\to \oo_{\bar z}$ is  an injective immersion. Define $\Gamma_{\bar z}=\Psi(\uu)$. 
By construction, $\Gamma_{\bar z}$ is an $n$-dimensional submanifold of $\RR^{n+p}$.

We next show that $\Gamma_{\bar z}$ satisfies the tangency property 2)~of the theorem, i.e.~ $\vs_i|_z\in T_z\Gamma_{\bar z}$,  for all $i=1,\dots m$ and  for every point $z\in \Gamma_{\bar z}$.    Since $\Psi$ is an injective immersion, for any $z\in \Gamma_{\bar z}$, there exists unique 
 $(w,\epsilon)\in \uu$, such that $z=\Psi(w,\epsilon)$. Since commutativity of the vector fields implies commutativity of the flows, we can pull $\exp^{\epsilon_j\,\tilde\vs_j}$ to the most left in \eqref{eq-extend}. Then, by \eq{eq-flow}  of an integral curve:
\begin{align*}
D\Psi|_{(w,\epsilon)}\left(\pd{}{\epsilon^j}\right)&=\left.\pd{}{\epsilon^j}\right|_{(w,\epsilon)}
\exp^{\epsilon_j\,\tilde\vs_j}\circ  \cdots \circ \exp^{\epsilon_{j-1}\,\tilde\vs_{j-1}}\circ \exp^{\epsilon_{j+1}\,\tilde\vs_{j+1}}\circ \cdots \circ \psi(w)\\
&=\vs_{j}|_z.
\end{align*}
It remains to construct  $\oo_{\bar z}$, such that the intersection property 1)~of the theorem is satisfied. Let $\ww$ be the projection of $\bar \uu=\uu\cap (\ww'\times \{0\})$ on $\RR^{n-m}$. Then  $\ww\subset\ww'\subset \RR^{n-m}$ is an open subset  containing the origin. Define
   $\Lambda_{\bar z}:=\psi(\ww)$, then $\Lambda_{\bar z}$ is  coordinate a chart on $\Lambda$ centered at  $\bar z$ . By construction,  $\Lambda_{\bar z}=\Psi(\bar\uu)=\Lambda\cap\Gamma_z$.   Since $\Lambda$ is embedded,  there exists 
open subset $\oo_z\subset\oo'_{\bar z}\subset\RR^{n+p}$ such that $\Lambda_{\bar z}=\oo_z\cap\Lambda$.   Then, by construction, $\Lambda_{\bar z}= \Lambda\cap \Gamma_{\bar z}=\Lambda\cap \oo_{\bar z}$.

Local uniqueness of $\Gamma_{\bar z}$ follows from the uniqueness of the  integral curve of a given vector field originating at a given point, combined with the fact that any  submanifold tangent to $\vs_j$ must  contain an open interval of the integral curve of $\vs_j$  originating at each  point of the submanifold.
\end{proof}
We now formulate and prove a PDE version of the generalized Frobenius Theorem.  A PDE system on $p$ functions of $n$ variables, considered in this theorem, prescribes derivative of each unknown function in the directions of $m\leq n$ vector fields comprising an involutive partial frame. We will call such systems  to be of    \emph{generalized Frobenius type}. The theorem claims that under some natural integrability conditions, there is a unique solution of this system with an initial data prescribed along an $m$-dimensional manifold transversal to the given partial frame. For $n=m$, this theorem is equivalent to the classical PDE version of the Frobenius Theorem (Theorem 6.1 in \cite{spivak:v1}).
\begin{theorem}[Generalized Frobenius] \label{th-pde-frob} Let $\R=\{\vr_1,\dots,\vr_m\}$ be a partial frame  in involution on an open subset $\Omega\subset \RR^n$ with coordinates $(u^1,\dots, u^n)$. 
Let $\Theta\subset \RR^p$ be an open subset with coordinates $(\phi^1,\dots,\phi^p)$. Let $h^i_j$, $ i=1,\dots, p,\ j=1,\dots, m$, be given smooth functions on $\Omega\times\Theta$. Consider a system of differential equations:
\beq\label{eq-frob}
\vr_j(\phi^i(u))=h^i_j(u, \phi(u)), \quad i=1,\dots, p;\ j=1,\dots, m.
\eeq
Assume the following integrability conditions 
\beq\label{eq-ic}
\vr_j\left(\vr_k(\phi^i)\right)-\vr_k\left(\vr_j(\phi^i)\right)=\sum_{l=1}^m c_{jk}^l\vr_l(\phi^i) \quad i=1,\dots, p;\ j,k=1,\dots, m,
\eeq
 where the functions $c$'s are defined by
 \beq\label{eq-r-str} [\vr_j,\vr_k]=\sum_{l=1}^m c_{jk}^l\vr_l, \eeq
 are identically satisfied on $\Omega\times\Theta$ after substitution of $h^i_j(u, \phi)$ for  $\vr_j(\phi^i(u))$ as prescribed by the system \eqref{eq-frob}\footnote{The resulting equations, explicitly written down as \eqref{eq-ic2}, involve no derivatives of $\phi$.}.\\
Then for any point $\bar u\in \Omega$ and  for any smooth initial data prescribed  along any   embedded submanifold $ \Xi\subset \Omega$ of codimension $m$ containing $\bar u$ and  transversal\footnote{Here transversality means that $\Span_{\RR}\{\vr_1|_{\bar u}, \dots, \vr_m|_{\bar u}\}\oplus T_{\bar u}\Xi=\RR^n$ at very point $\bar u\in \Xi$, where  $T_{\bar u}\Xi$ denotes the tangent space to $\Xi$ at $\bar u$.}to $\R$, there is a unique smooth local solution of  \eqref{eq-frob}.  
In other words, given  arbitrary functions $(g^1,\dots,g^p)\colon \Xi\to\Theta$, there is an open subset $\Omega'\subset\Omega$, containing $\bar u$ and smooth functions $(\alpha^1,\dots, \alpha^p )\colon\Omega'\to \Theta$ satisfying \eqref{eq-frob},  such that $\alpha^i|_{\Xi\cap\Omega'}=g^i, i=1,\dots, p$.   \end{theorem} 
\begin{proof}
Before staring a proof, we expand conditions \eqref{eq-ic}. After the first substitution  of  the derivatives of $\phi$'s as prescribed by \eqref{eq-frob} into  \eqref{eq-ic}, we get  for $ i=1,\dots, p;\ j,k=1,\dots, m$:
\beq\label{eq-ic1}
\vr_j\left(h^i_k(u,\phi(u)\right)-\vr_k\left(h^i_j(u,\phi(u)\right)=\sum_{l=1}^m c_{jk}^l\,h_l^i \left(u,\phi(u)\right)
\eeq
Using the chain rule and again making a substitution  prescribed by  \eqref{eq-frob} for the derivatives of $\phi$'s we get:
\begin{eqnarray}\nonumber 
&& \sum_{l=1}^n \left(\pd{ h^i_k\left(u,\phi\right)}{u^l}\vr_j(u^l)-\pd{ h^i_j\left(u,\phi\right)}{u^l}\vr_k(u^l) \right)+ \sum_{s=1}^p \left(\pd{h^i_k\left(u,\phi\right)}{\phi^s} \,h^s_j\left(u,\phi\right)-\pd{h^i_j\left(u,\phi\right)}{\phi^s} \,h^s_k\left(u,\phi\right)\right)\\
\label{eq-ic2} &=&\sum_{l=1}^m c_{jk}^l(u)\,h_l^i \left(u,\phi\right).
\end{eqnarray}
In order to use Theorem~\ref{th-vf-frob}, for $j=1,\dots, m$, we define vector-fields $$\vs_j=\vr_j+  \sum _{i=1}^p\,h_j^i\pd{}{\phi^i}$$
on the open subset $\Omega\times\Theta\subset \RR^{n+p}$.  
Independence of $\vs_1,\dots,\vs_m$ follows from independence of $\vr_1,\dots, \vr_m, \pd{}{\phi^1},\dots, \pd{}{\phi^p}$. Using \eqref{eq-r-str} and \eqref{eq-ic1}, we can show that  $\vs_1,\dots,\vs_m$ are in involution, and in fact satisfy the same structure equations as  $\vr_1,\dots,\vr_m$. Indeed,  for all $j,k=1,\dots, m$:
$$[\vs_j,\vs_k]=[\vr_j,\vr_k]+ \sum_{i=1}^p \left(\vr_j(h_k^i) -\vr_k(h_j^i)\right)\, \pd{}{\phi^i}=\sum_{l=1}^m c_{jk}^l\, \left(\vr_l+  \sum _{i=1}^p\,h_l^i\pd{}{\phi^i}\right)=\sum_{l=1}^m c_{jk}^l\,\vs_l.$$
Define $$\Lambda=\{(u, g(u)\,|\,u\in\Xi\},$$
to be the graph of the map $g=(g^1,\dots,g^p)\colon \Xi\to \RR^{p}.$ Then $\Lambda$ is an $(n-m)$-dimensional embedded submanifold of $\RR^{n+p}$, such that  \eqref{eq-span} is satisfied.
Thus by Theorem~\ref{th-vf-frob}, there  exist an open neighborhood $\oo\subset \Omega\times\Theta$ of  $\bar z=g(\bar u)$ an $n$-dimensional manifold $\Gamma$, through the point $\bar z=g(\bar u)$, such that
$\Gamma\cap\Lambda =\oo\cap\Lambda$ and vector-fields $\vs_1,\dots,\vs_m$ are tangent to $\Gamma$ at every point of $\Gamma$. By possibly shrinking  $\oo$  and $\Gamma$ around $\bar z$, we may assume that 
$$\Gamma=\{(u, \alpha(u)\,|\,u\in\Omega'\},$$ 
is a graph of a map $\alpha=(\alpha^1,\dots, \alpha^p)\colon \Omega' \to \Theta$, where $\Omega'$ is an open  subset of $\Omega$, equal to the projection of $\oo$ to $\RR^n$. 
Since $\Gamma\cap\Lambda =\oo\cap\Lambda$,  we have $\alpha|_{\Xi\cap\Omega'}=g$.
Then since $\vs_1,\dots,\vs_m$ are tangent to $\Gamma$, we have
$$\vs_j\left(\phi^i-\alpha^i(u)\right)=0 \text{ for all }  u\in \Omega' \text { and }   i=1,\dots, p;\ j=1,\dots, m,$$
which is equivalent to 
$$\vr_j(\alpha^i)= h^i_j(\alpha^i(u),u), \text{ for all }  u\in \Omega' \text { and }   i=1,\dots, p;\ j=1,\dots, m. $$
and therefore $\alpha:\Omega' \to \Theta $ is  a solution of  the PDE system \eqref{eq-frob}. Its uniqueness follows from the local uniqueness of $\Gamma$.   
\end{proof}
We conclude this section by  stating another integrability theorem, appeared as Theorem III in Book III, Chapter I of \cite{darboux}.  The PDE system on $p$ functions of $n$-variables, considered in this theorem, prescribes  some subset of partial derivatives of each unknown functions. A subset of derivatives prescribed for one of the unknown functions, may differ from a subset prescribed for the other. We will call such systems to be  of  the \emph{Darboux type}. The Darboux theorem claims that provided the natural integrability conditions are satisfied, there is a unique solution for an appropriately prescribed initial data.  

\begin{theorem}\label{dar3} (Darboux \cite{darboux}) Let $\Omega\subset \RR^n$ and $\Theta\subset \RR^p$ be open subsets, let $\bar u=(\bar u^1,\dots,\bar u^n)\in \Omega$ be a fixed point, and let $h^i_j$, $ i=1,\dots, p,\ j\in {S_i}$, where $S_i\subset\{1,\dots, n\}$ is a fixed subset, be some given smooth functions on $\Omega\times\Theta$. Consider a system of differential equations on unknown functions $(\phi^1,\dots \phi^p)\colon \Omega\to\Theta$ of independent variables $u^1,\dots,u^n$: 
 \beq\label{sys_ii}
	\frac{\del \phi^i}{\del u^j} = h_{j}^i(u,\phi), \quad \ j\in {S_i},\,  i=1,\dots, p.
\eeq
 Assume that the system prescribes compatible second order mixed 
derivatives in the following sense:
\begin{itemize}
	\item[(C)] Whenever two distinct derivatives $\frac{\del \phi^i}{\del u^j}$ and 
	$\frac{\del \phi^i}{\del u^{k}}$ of the same unknown $\phi^i$ are present 
	on the left hand side of \eq{sys_ii}, then the equation
	\[\frac{\del}{\del u^{k}}\big[h_{j}^i(u,\phi(u))\big] = \frac{\del}{\del u^{j}}\big[h_{k}^i(u,\phi(u))\big]\]
	contains (after expanding each side using the chain rule) 
	only first order derivatives which appear in \eq{sys_ii}, and substitution from 
	\eq{sys_ii} for these first derivatives results in an identity in $u$ and $\phi$.
\end{itemize} 
Next, to describe the data, suppose a dependent variable $\phi^i$ appears 
differentiated in \eq{sys_ii} with respect to $u^{j_1},\dots,u^{j_{s}}$. 
Then, letting $\tilde u$ denote the remaining independent variables, 
we prescribe a smooth function $g^i(\tilde u)$ and require that
\beq\label{DDdata}
	\phi^i(u^1,\dots,u^n)\Big|_{u^{i_1}=\bar u^{i_1},\dots,\, u^{i_s}=\bar u^{i_s}}
	=g^i(\tilde u)\,.
\eeq
We make such an assignment of data for each $\phi^i$ that appears 
differentiated in \eq{sys_ii}. Then, under the compatibility condition (C), the problem 
\eq{sys_ii} - \eq{DDdata} has a unique, local smooth solution for $u$ near $\bar u$.
\end{theorem}

\begin{remark}
%
	If partial derivatives of all unknown functions are prescribed for the same set coordinates directions (i.e $S_1=\dots =S_p$), the Darboux type system  is  of the generalized  Frobenius type. Conversely, using Propositions~\ref{prop-inv-distr} and~\ref{prop-comm-frame}, one can show that for any system of  the generalized Frobenius type there is an equivalent  Darboux type system,  with all partial derivatives of all unknown functions prescribed for the same set coordinates directions.   In this case, integrabilty conditions (C) of Theorem~\ref{dar3} are equivalent to integrability conditions in Theorem~\ref{th-pde-frob}.  However, the manifold $\Xi$ along which the initial data is allowed to be prescribed in Theorem~\ref{th-pde-frob} is more general than the coordinate subspace for which the data is prescribed in Theorem~\ref{dar3}.
	\end{remark}

\section{Coordinate-free formulation of the problem}\label{sect-intrinsic}
In this section, we give an intrinsic (coordinate independent) formulation of Problem~\ref{pr}, which leads to a system of  differential equations written in terms of the frame adapted to the problem. We derive some differential consequences of this system, which, in particular, lead to a set of necessary conditions for the existence of strictly hyperbolic fluxes.  
\subsection{Intrinsic definition of the Jacobian and the $\fs(\R)$-system}\label{sect-intrinsic-jac}

We start with an intrinsic  definition of the Jacobian map $\X(\Omega)\to\X (\Omega)$ adapted  from Remark~2.15 in \cite{jk2}.

\begin{defin} \label{def-jac-free} Given  a  connection $\nabla$ on a smooth manifold $M$, the \emph{$\nabla$-Jacobian} of a vector field  $\vf\in\mathcal X(M)$ is the $C^{\infty}(M)$-linear map 
$ J\vf\colon \mathcal X(\Omega)\to {\mathcal X(\Omega)}$ 
defined by:
\beq\label{jac-free} 
	 \quad J\vf(\vs)=\nabla_{\vs} \vf, \quad \forall \vs\in \mathcal X(\Omega). 
\eeq
\end{defin}
 If  $\{\vr_1,\dots,\vr_n\}$ is a frame with Christoffel symbols 
$\Gamma_{ij}^k$ and $\vf=\sum_{i=1}^n F^i\,\vr_i$ then \eq {jac-free} implies:
\beq\label{Jr}
	 J\vf(\vr_j)=\sum_{i=1}^n\left(\vr_j (F^i)
	+\sum_{k=1}^n\Gamma^i_{jk}\,F^k\right)\,\vr_i.
\eeq

Let $(u^1,\dots, u^n)$ is an affine system of coordinates (see \eq{eq-affine}) relative to a flat symmetric connection  $\tnabla$,  let 
$\vf=\sum_{i=1}^n \tilde F^i(u)\pd{}{u^i}$, and  let  $\tJ \vf$ denote the $\tnabla$-Jacobian of $\vf$. Then a direct computation shows that 
$$\tJ\vf\left(\pd{}{u^j}\right)=\sum_{i=1}^n\pd{\tilde F^i}{u^j}\,\pd{}{u^i},$$ 
which corresponds  exactly to the $j$-th column 
vector of the usual Jacobian matrix $[D_u\,\tilde F]$ of  the vector valued function $\tilde F(u)=[\tilde F^1(u),\dots,\tilde F^n(u)]^T$.

Using the intrinsic definition of the Jacobian, we give an equivalent intrinsic formulation of  Problem~\ref{pr}, which  allows us to analyze it relative to a frame that is adapted to the problem.

\begin{problem}\label{pr2}
Given a partial frame $\R=\{\vr_1,\dots, \vr_m\}$ on an open subset $\Omega\subset \RR^n$ ($n\geq m$), with a flat symmetric connection $\tilde \nabla$, and a fixed point $\bar u\in \Omega$;
describe the set $\fs(\R)$ of smooth vector fields $\vf$  for which there exist an open neighborhood $\Omega'\subset \Omega$ of $\bar u$ and  smooth functions $\lam^i\colon \Omega' \to\RR$, such that 
\begin{equation} \label{eq-fr}
  		\tnabla_{\vr_i} \,\vf = \lambda^i \, \vr_i, \text{ on  } \Omega' \text{ for } i=1,\dots,m .
  		\end{equation}
\end{problem}
\begin{remark} 
 Problem~\ref{pr2} makes sense if we replace $\RR^n$ with an arbitrary manifold $M$, and replace $\tnabla$ with an arbitrary connection on  the tangent bundle of $M$. In particular, it would be of interest to consider this problem on a Rimannian manifold with the Riemannian connection. These generalizations, however, fall  outside of the scope of the current paper.

\end{remark}
From Proposition~\ref{prop-flat-symm}, we know that any flat and symmetric connection admits an affine system of coordinates.
If $F^1,\dots, F^n$ are  the components of $\vf$, and $R^1_i,\dots, R^n_i$ are the components of $\vr_i $  in an affine system of coordinates, then 
\eqref{eq-fr} turns in a system of $mn$ first order PDE's on $n+m$ unknown functions $F$'s and $\lam$'s:

\begin{equation} \label{eq-fr-coord}
  		r_i(F^j)  = \lambda^i \,R_i^j, \quad \text{ for } i=1,\dots,m,\, j=1,\dots, n,
  		\end{equation}
which is equivalent to \eqref{eq:eigenvectors},  because $r_i(F^j)=\sum R_i^k \,\pd{F^j}{u^k}$. Therefore, Problem~\ref{pr2} is, indeed, a coordinate-free formulation of  Problem~\ref{pr}, and 
we can call system \eqref{eq-fr} the  $\fs(\R)$-system. The set of vector fields  satisfying   \eqref{eq-fr} will be denoted $\fs(\R)$, and  the elements of this set will be  called fluxes for $\R$.  The set  fluxes always contains  the set of \emph{identity} fluxes $\fsid$, which we define by the property:
\beq\label{def-idf} \tnabla_{\vr}\vf=\vr \text { for all vector fields } \vr\in\mathcal X(\Omega).\eeq  
One can easily show that $\hat \vf\in \fsid$ if and only if relatively to any affine coordinates system $(u^1,\dots, u^n)$:
$$\hat \vf=[u^1,\dots, u^n]^T+\bar b, \text{ for some } \bar b\in\RR^n.$$
The previously defined vector space of trivial fluxes \eq{tflux}, in this more abstract setting, corresponds to the vector space
 \beq\label{def-triv} \fstriv=\{\vf\in \mathcal X(\Omega)\,|\, \forall \vr \in \mathcal X(\Omega)\, \exists \bar\lam\in \RR, \text {  such that }  \tnabla_{\vr}\vf=\bar\lam \vr\}.
 \eeq
Equivalently, one can say that $ \fstriv=\{\bar\lam\, \hat \vf\,|\, \bar\lam \in \RR, \, \hat\vf\in\fsid \}$ and, 
 clearly, $\fsid\subset\fstriv\subset \fs(\R)$ for any partial frame $\R$.

\subsection{Differential consequences of the $\fs(\R)$-system}\label{sect-fr}


We now derive the differential consequences of \eqref{eq-fr} implied by the flatness of the connection.
\begin{proposition}\label{prop-lam-a} Given a partial frame $\R=\{\vr_1,\dots,\vr_m\}$, assume that $\vf\in \fs(\R)$ is a flux, and $\vs_1,\dots, \vs_{n-m}$ is any completion of $\R$ to the full frame. Let the functions $\aaa^l_k$ be defined by
\beq\label{eq-s}\tnabla_{\vs_l}\vf=\sum_{k=1}^m \aaa_l^k \vr_k+ \sum_{k=m+1}^n \aaa_l^k\,\vs_t,\quad l=1,\dots, n-m.\eeq
Then the functions $\lam^i$, $i=1,\dots,m$, prescribed by \eq{eq-fr}, and  the functions $\aaa^l_k$, $l=1,\dots,n-m$, $k=1,\dots,n$,   satisfy the system of differential and algebraic equations:

\begin{align}
\label{gc1}{\vr_i} (\lambda^j)&= \Gamma^j_{ji}\,(\lambda^i-\lambda^j)+ \sum_{l=m+1}^n   a_l^j \, c_{ij}^l \text { for all  } 1\leq i\neq j\leq m\\
\label{gc2} \lambda^j\, \Gamma_{ij}^k-\lambda^i\, \Gamma_{ji}^k -c_{ij}^k \lambda^k&= \sum_{l=m+1}^n   a_l^k \, c_{ij}^l \text{ for all distinct triples } i, j ,k\in \{1,\dots, m\}\\
\label{gc3} \lambda^j\,\Gamma_{ij}^l-\lambda^i\,\Gamma_{ji}^l&=\sum_{t=m+1}^n   a_t^l \, c_{ij}^t \text { for all  } 1\leq i\neq j\leq m \text { and } l=m+1,\dots, n.
\end{align}
In the above equations, $c$'s and $\Gamma$'s are the structure functions and the Christoffel symbols for the frame:
\begin{eqnarray}
\label{eq-c} [\vr_i,\vr_j]&=&\sum_{k=1}^m c_{ij}^k \vr_k+ \sum_{l=m+1}^n c_{ij}^l\,\vs_l,\\
\label{eq-gamma}\tnabla_{\vr_i}{\vr_j}&=&\sum_{k=1}^m \Gamma_{ij}^k \vr_k+ \sum_{l=m+1}^n \Gamma_{ij}^l\,\vs_l.
\end{eqnarray}
\end{proposition}
\begin{proof}

 Flatness condition \eqref{R0} implies that
 \beq\label{R0f}\tnabla_{\vr_i}\tnabla_{\vr_j}\vf-\tnabla_{\vr_j}\tnabla_{\vr_i}\vf=\tnabla_{[\vr_i,\vr_j]}  \vf\quad \text { for all } i, j=1,\dots, m\eeq 
must hold on the solutions of \eqref{eq-fr}, and, therefore,
 \beq\label{R0f1}{\vr_i} (\lambda^j)\,\vr_j +\lambda^j\,\tnabla_{\vr_i}{\vr_j}-{\vr_j} (\lambda^i)\,\vr_i -\lambda^i\,\tnabla_{\vr_j}{\vr_i}=\tnabla_{[\vr_i,\vr_j]} \vf\eeq 

Using  \eq{eq-c} and \eq{eq-gamma}, we obtain that  \eqref{R0f1} is equivalent to:
\begin{eqnarray}\nonumber &&{\vr_i} (\lambda^j)\,\vr_j +\sum_{k=1}^m\lambda^j\, \Gamma_{ij}^k \vr_k+ \sum_{l=m+1}^n \lambda^j\Gamma_{ij}^l\,\vs_l-{\vr_j} (\lambda^i)\,\vr_i -\sum_{k=1}^m\lambda^i\, \Gamma_{ji}^k \vr_k- \sum_{l=m+1}^n \lambda^i\Gamma_{ji}^l\,\vs_l\\
\label{eq-almost} &=&\sum_{k=1}^m c_{ij}^k \tnabla_{\vr_k}\,\vf+ \sum_{l=m+1}^n c_{ij}^l\,\tnabla_{\vs_l} \,\vf
\end{eqnarray}
It remains to rewrite the right-hand side of \eqref{eq-almost} in terms of the frame using  \eqref{eq-fr} for the first sum and  \eq{eq-s} for the second sum: 
\begin{eqnarray*}&&{\vr_i} (\lambda^j)\,\vr_j +\sum_{k=1}^m\lambda^j\, \Gamma_{ij}^k \vr_k+ \sum_{l=m+1}^n \lambda^j\Gamma_{ij}^l\,\vs_l-{\vr_j} (\lambda^i)\,\vr_i -\sum_{k=1}^m\lambda^i\, \Gamma_{ji}^k \vr_k- \sum_{l=m+1}^n \lambda^i\Gamma_{ji}^l\,\vs_l\\
&=&\sum_{k=1}^m c_{ij}^k \lambda^k\,\vr^k+\sum_{k=1}^m \sum_{l=m+1}^n   a_l^k \, c_{ij}^l \vr_k + \sum_{l,t=m+1}^n   a_l^t \, c_{ij}^l \vs_t
\end{eqnarray*} 
Collecting coefficients of the frame vector fields, we get a system of  differential and algebraic  equations \eq{gc1}--\eq{gc3}.
\end{proof}
It is worthwhile emphasizing that, in general, the structure functions $c$'s and the Christoffel symbols $\Gamma$'s appearing in  \eq{gc1}--\eq{gc3}, depend on the completion of $\R$ to a full frame. 
\begin{remark}\label{rem-lambda-a}
We  note that  equations \eq{gc1}--\eq{gc3}  do not provide a complete set of integrability conditions for  the Frobenius-type system \eqref{eq-fr}, \eqref{eq-s}, because it does not include  conditions derived from   $\tnabla_{\vr_i}\tnabla_{\vs_j}\vf-\tnabla_{\vs_j}\tnabla_{\vr_i}\vf =\tnabla_{[\vr_i,\vs_j]}$ and $\tnabla_{\vs_i}\tnabla_{\vs_j}\vf-\tnabla_{\vs_j}\tnabla_{\vs_i}\vf =\tnabla_{[\vs_i,\vs_j]}$.  We will derived these additional conditions, in Section~\ref{sect-noninvol}, for $m=2$, $n=3$ case only, and we will observe how technical they become even in low dimensions.

 However, we will see in Section~\ref{sect-involutive}, that if $\R$ is an \emph{involutive} partial frame, then  \eq{gc1} -- \eq{gc3} simplify to a system which involves unknown functions $\lambda$s only, and this system does  provide a complete set of integrability conditions for  \eqref{eq-fr}. In the case of the full frame $(m=n)$, equations \eq{gc1} -- \eq{gc3} reduce to the $\lam$-system introduced in \cite{jk1}.
\end{remark}
 We can use equations \eq{gc1} -- \eq{gc3} to obtain  necessary conditions for $\fs(\R)$ to contain a strictly hyperbolic flux. As we will see below, these conditions are not sufficient except for the case of rich partial frames.
\begin{proposition}[necessary condition for strict hyperbolicity]\label{prop-nec-shg} Let $\R=\{\vr_1,\dots,\vr_m\}$ be a partial frame on $\Omega\subset \RR^n$ containing $\bar u$. If there is a strictly hyperbolic flux $\vf\in \fs(\R)$ on some open  neighborhood  $\Omega'$ of $\bar u$ then for each pair of indices $i\neq j\in\{1,\dots, m\}$ the following  equivalence condition holds:
\beq\label{eq-nec-shg}
\tnabla_{\vr_i}\vr_j\in \Span_{C^\infty(\Omega')}\{\vr_i, \vr_j\}\text{ if and only if } [\vr_i,\vr_j]\in \Span_{C^\infty(\Omega')}\{\vr_i, \vr_j\}\eeq
\end{proposition}
\begin{proof} If $\vf$ is strictly hyperbolic on $\Omega'$,  then $\R$ can be completed to a  frame of eignvectors $\vr_1, \dots,\vr_m$, $\vr_{m+1},\dots, \vr_n$, such that there exist functions $\lam^1,\dots, \lam^n\colon\Omega'\to\RR$,  with all distinct values  at each point of $\Omega'$,  and
$$\tnabla_{\vr_i} \vf=\lam^i\vr_i, \quad i=1,\dots, n.$$
In the statement of Proposition~\ref{prop-lam-a}, let $\vs_l=\vr_l$ for $l=m+1,\dots, n$. Then $a_l^i=\delta^i_l \lam^l$, where $\delta^i_l$ is  the Kronecker  delta function, and  the  algebraic  conditions \eqref{gc2}, \eqref{gc3}  become 
\beq\label{eq-necg} 
\Gamma_{ij}^k\,\lambda^j - \Gamma_{ji}^k\lambda^i\,  -c_{ij}^k \lambda^k =0  \text{ for all }  1\leq i\neq j\leq m \text { and } 1 \leq k \leq n, \text { such that } k\neq i \text { and }  k \neq  j.
\eeq
Let us first assume that  for some $i, j$, such that  $1\leq i\neq j\leq m$, we have $\tnabla_{\vr_i}\vr_j\in \Span_{C^\infty(\Omega')}\{\vr_i, \vr_j\}$ and  $[\vr_i,\vr_j]\notin \Span_{C^\infty(\Omega')}\{\vr_i, \vr_j\}$. Then,  from  the latter condition,   there exists  $k\in\{1,\dots,  n\} $, such that  $k\neq i \text { and }  k \neq  j$ and $c_{ij}^k\not\equiv 0$, while  the former condition implies that $\Gamma_{ij}^k\equiv 0$.  Symmetry of $\tnabla$ implies that $c_{ij}^k=-\Gamma^k_{ji}\not\equiv 0$, and then from \eqref{eq-necg} we have
$$   c_{ij}^k\,(\lam^i-\lam^k)\equiv0.$$
We then have $\lam^i=\lam^k$ at least somewhere in $\Omega'$, which contradicts our strict hyperbolicity assumption.

Let us now assume that for some $i, j$ such that  $1\leq i\neq j\leq m$, we have $\tnabla_{\vr_i}\vr_j\notin \Span_{C^\infty(\Omega')}\{\vr_i, \vr_j\}$ and  $[\vr_i,\vr_j]\in \Span_{C^\infty(\Omega')}\{\vr_i, \vr_j\}$. Then, from the former condition,  there  exists  $k\in\{1,\dots,  n\} $, such that  $k\neq i \text { and }  k \neq  j$ and $\Gamma_{ij}^k\not\equiv 0$, but from the latter condition we have $c_{ij}^k\equiv 0$.  Symmetry of $\tnabla$ implies that $\Gamma_{ij}^k=\Gamma^k_{ji}$, and then from \eqref{eq-necg} we have
$$   \Gamma_{ij}^k\,(\lam^j-\lam^i)\equiv0.$$
We then have $\lam^i=\lam^j$ at least somewhere  in $\Omega'$, which contradicts our strict hyperbolicity assumption.
\end{proof}

 
\section{Involutive partial frame} \label{sect-involutive}
As it is discussed in Remark~\ref{rem-lambda-a} above, the analysis of the $\fs(\R)$ system is much simpler when the partial frames $\R$ is  in involution. Partial frames of two  ``extreme'' sizes:  $m=1$, that is all partial frames consisting of a single vector field,  and $m=n$, that is all full frames,  fall into this category. 
In the former case $\R$  also trivially satisfies the definition of a rich partial frame - the case that we treat in Section~\ref{sect-rich} (we make some comments about $m=1$ case in Remark~\ref{rem-single}). The case of the full frame was considered in details in \cite{jk1} and some of the theorems of this paper are natural generalizations of these results.  We will first state the results that apply to arbitrary  involutive partial frames, then consider rich partial frames, which is a particularly nice case subclass of involutive  partial frames, and finally we will consider non-rich    involutive partial frames partial frames consisting of three vector-fields.
\subsection{Arbitrary involutive partial framas}\label{ssect-inv}
If a given partial frame $\R$ is in involution, then for any completion of $\R$ to a full frame $\{\vr_1,\dots, \vr_m, \vs_{m+1},\dots,\vs_{n}\}$, we have $ c_{ij}^l=0$   for all $i,j=1,\dots, m$, $l=m+1,\dots, n$ and, therefore,  $\Gamma_{ij}^l=\Gamma_{ji}^l$ due to the symmetry of the connection \eqref{T0}. Then \eqref{gc1} -- \eqref{gc3} simplify to
\begin{eqnarray}
\label{ic1}{\vr_i} (\lambda^j)&=& \Gamma^j_{ji}\,(\lambda^i-\lambda^j)\text { for all  } 1\leq i\neq j\leq m\\
\label{ic2} \lambda^j\, \Gamma_{ij}^k-\lambda^i\, \Gamma_{ji}^k -c_{ij}^k \lambda^k&=&0 \text{ for all distinct triples } i, j ,k\in \{1,\dots, m\}\\
\label{ic3} (\lambda^j-\lambda^i)\,\Gamma_{ji}^l&=&0\text { for all  } 1\leq i\neq j\leq m \text { and } l=m+1,\dots, n.
\end{eqnarray}
where $c$'s and $\Gamma's$ are defined by \eqref{eq-c} and \eqref{eq-gamma}, respectively.  Note that, due to involutivity of $\R$, functions $c_{ij}^k$,  $i, j ,k\in \{1,\dots, m\}$ do not depend on a choice of completion of $\R$ to a frame, while  $\Gamma$'s, in general, do depend on a choice of such completion. We will call   \eqref{ic1} -- \eqref{ic3} \emph{the $\lam$-system}, generalizing the terminology of \cite{jk1} to partial involutive frames.

The following proposition allows us, in the involutive case, to solve Problems~\ref{pr2} (and \ref{pr}) in two steps: first find (or describe the set of) all solutions $\lam$ of  system~\eqref{ic1} -- \eqref{ic3}, and then find (or describe the set of) all solutions    $\vf$ of \eqref{eq-fr} for a given set of functions $\lam$'s, satisfying \eqref{ic1} -- \eqref{ic3}.   This is possible because in the involutive case, equations  \eqref{ic1} -- \eqref{ic3} provide a complete set of the integrability conditions for the $\fs(\R)$-system  \eqref{eq-fr},as we show in the proof of the following proposition.
\begin{proposition}[$\lam(\R)$-system] \label{prop-inv} If a partial  frame $\R=\{\vr_1,\dots,\vr_m\}$ is in involution, then
\begin{enumerate}
\item[1)] For every $\vf\in\fs(\R)$, functions $\lam^1,\dots,\lam^m$ prescribed by  \eqref{eq-fr} satisfy  \eqref{ic1} -- \eqref{ic3}.
\item[2)] For every solution  $\lambda^1,\dots, \lambda^m$  of   \eqref{ic1} -- \eqref{ic3},  
and  any smooth initial data for $\vf$ prescribed  along any   embedded submanifold $\Xi\subset \Omega$ of codimension $m$ transverse to $\R$,
there is a unique smooth local solution of $\fs(\R)$-system \eqref{eq-fr}.
In other words, given  arbitrary  smooth vector field $\tilde \vf$ on $\Xi$, there is an open subset $\Omega'\subset\Omega$, containing $\Xi$ and unique smooth extension  $\vf$ of $\tilde \vf$ to $\Omega'$ satisfying \eqref{eq-fr}.
\end{enumerate}
\end{proposition}
\begin{proof} \hfill
\begin{enumerate}
\item[1)]  Equations  \eqref{ic1} -- \eqref{ic3} are differential consequences of  \eqref{eq-fr}, and, therefore, for every $\vf\in\fs(\R)$, functions $\lam^1,\dots,\lam^m$ prescribed by  \eqref{eq-fr} satisfy  \eqref{ic1} -- \eqref{ic3}.

\item[2)]Assume $\lambda^1,\dots, \lambda^m$  are solutions of   \eqref{ic1} -- \eqref{ic3}. In an affine system of coordinates $u=(u^1,\dots,u^n)$,  equations  \eqref{eq-fr} turn into \eqref{eq-fr-coord}.  To simplify the notation  we rewrite them  in a vector-form:
\begin{equation} \label{rF}
  		r_i(F)|_u  = \lambda^i(u) \,R_i(u), \quad \text{ for } i=1,\dots,m,
  		\end{equation}
where  we assume that
$F$  and $R_i$ are column vectors of the components of the vector-fields $\vf$ and $\vr_i$ respectively relative to the coordinate frame $\pd{}{u^1},\dots, \pd{}{u^n}$. 

The above system  is of the form \eq{eq-frob} described in  generalized Frobenius Theorem~\ref{th-pde-frob}.
The integrability conditions  become:
\beq\label{R0f2}{\vr_i} (\lambda^j)\,(R_j) +\lambda^j\,{\vr_i}{(R_j)}-{\vr_j} (\lambda^i)\,(R_i) -\lambda^i\,{\vr_j}{(R_i)}=\sum_{k=1}^m c_{ij}^k\lambda^k R_k \eeq
Recalling that in the affine coordinates we have  formula \eqref{eq-tr-conn} for covariant derivatives, we see that \eqref{R0f2} is equivalent to
\beq\label{R0f3}{\vr_i} (\lambda^j)\,\vr_j +\lambda^j\,\tnabla_{\vr_i}{\vr_j}-{\vr_j} (\lambda^i)\,\vr_i -\lambda^i\,\tnabla_{\vr_j}{\vr_i}=\sum_{k=1}^m c_{ij}^k\lambda^k \vr_k,\eeq 
which, when written out in components relative to a completion of $\R$ to a frame  $\vr_1,\dots,\vr_m,$ $\vs_{m+1},\dots, \vs_n$, is equivalent to  \eqref{ic1} -- \eqref{ic3}.
Components of the vector field $\tilde\vf$ provide the data   for $F$ of the type described  in Theorem~\ref{th-pde-frob}, and this theorem guarantees the existence of  a locally  unique solution of \eq{rF} with this data.
\end{enumerate}
\end{proof}


System  \eqref{ic1} -- \eqref{ic3} always has solution $\lam^1=\dots=\lam^m$, but existence of other solutions of  \eqref{ic1} --   \eqref{ic3} is a subtle question. Moreover, even for non-trivial solutions of \eqref{ic1} -- \eqref{ic3}, existence of hyperbolic and strictly hyperbolic fluxes is a subtle question.  

We note that  conditions \eqref{ic2} and \eqref{ic3} immediately provide us with  necessary conditions to existence of strictly hyperbolic solutions for Problem~\ref{pr}, in the case of involutive partial frames.  
\begin{proposition}[necessary condition for strict hyperbolicity in the involutive case]\label{prop-nec} If a partial  frame $\R=\{\vr_1,\dots,\vr_m\}$ is in involution, then the following  conditions must be satisfied  for all $1\leq i\neq j\leq m$ on some open neighborhood $\Omega' \subset \Omega$ of $\bar u$, in order for the $\fs(\R)$ set to contain a strictly hyperbolic flux:
\beq\label{eq-nec-sh2}
\tnabla_{\vr_i}\vr_j\in \Span_{C^\infty(\Omega')}\{\vr_i, \vr_j\}\Longleftrightarrow [\vr_i,\vr_j]\in \Span_{C^\infty(\Omega')}\{\vr_i, \vr_j\}\eeq
and
\beq\label{eq-nec-sh3}
\tnabla_{\vr_i}\vr_j\in \Span_{C^\infty(\Omega')}{\R},\eeq
\end{proposition}
\noindent As a side remark,  we observe that  involutivity implies that $[\vr_i,\vr_j]\in \Span_{C^\infty(\Omega')}{\R}$, and hence, due to the symmetry condition \eqref{T0}, we can replace condition $ 1\leq i\neq j\leq m$ in \eqref{eq-nec-sh3} with  $ 1\leq i< j\leq m$.
\begin{proof} Condition \eq{eq-nec-sh2} is the same  as  \eq{eq-nec-shg} proved earlier. If  for all open subsets $\Omega'\subset\Omega$, there are $1\leq i\neq j\leq m$, such that $\tnabla_{\vr_i}\vr_j\notin \Span_{C^\infty(\Omega')}{\R}$, then 
there exists $m+1\leq l \leq n$, such that $\Gamma_{ij}^l\not\equiv 0$ on $\Omega'$. From \eqref {ic3}, it then follows that $\lam^i=\lam^j$ at least somewhere on $\Omega'$ and therefore $\fs(\R)$ contains no stricitly hyperbolic fluxes.
\end{proof}
The above conditions are not sufficient as will be illustrated by  Example~5.3 in \cite{jk1}.
However, we can prove the following condition is sufficient.
\begin{proposition}[sufficient condition for strict hyperbolicity in the involutive case]\label{prop-suff} Assume, that functions  $\lam^1,\dots,\lam^m$ satisfying    \eqref{ic1} -- \eqref{ic3} are such that for some $\bar u\in \Omega$, all $m$ numbers $\lam^1(\bar u),\dots, \lam^m(\bar u)$ are distinct. Then on an open neighborhood $\bar u$ there exists a  strictly hyperbolic flux $\vf$, such that
$$ \tnabla_{\vr_i}\vf=\lam^i\vr_i,\quad i=1,\dots, m. $$
\end{proposition}
\begin{proof} Let $R_i$ be the column vector of components of $\vr_i$ in an affine system of coordinates $u^1,\dots, u^n$, and let $R=[R_1| \dots |R_m]$  be an $n\times m$ matrix comprised  of these column vectors. Since $\vr_i$, $i=1,\dots, m$ are independent at $\bar u$, there is a non-zero $m\times m$ minor of $R(\bar u)$. Due to continuity the same minor is non-zero on some open neighborhood of $\bar u$.   Let $\{i_1,\dots, i_m\}$ be the row indices of the submatrix corresponding to this minor. Up to permuting coordinate functions $u^1,\dots, u^n$ we may, in order to simplify the notation, assume that $i_j=j$. Then the set of  vector-fields $\vr_1,\dots, \vr_m, \pd{}{u^{m+1}},\dots,\pd{}{u_n}$ are independent and, therefore, a submanifold $\Xi$ defined by $u^i=\bar u^i$ for $i=1,\dots, m$ is  transversal to $\R$.   

 For $l\in\{m+1,\dots, n\}$, 
choose arbitrary constants $\bar \lambda^l$, such that all $n$ real numbers $\lam^1(\bar u),\dots,$ $\lam^m(\bar u)$, $\bar\lam^{m+1}, \dots,  \bar \lam^n$ are distinct. Define  
$$\tilde F(\bar u^1,\dots, \bar u^m, u^{m+1},\dots, u^n)=\left[0,\dots, 0, \,\bar\lam^{m+1}\, u^{m+1},\dots, \bar\lam^{n}\, u^{n}\right]^T$$
and let $F$ be an extension of $\tilde F$ such that  $[D_uF]\,R_i(u)=\lam^i\,R_i(u)$, where $i=1,\dots, m$ guaranteed by Proposition~\ref{prop-inv}. Then
$$[D_uF](\bar u)=
\left[
\begin{array}{cccccc}
  \pd{F^1}{u^1}(\bar u)&\dots &  \pd{F^1}{u^1}(\bar u)    &                                     & \dots                &  \\
  \vdots&\dots &  \vdots  &                                     & \dots                &  \\
   \pd{F^{m}}{u^1}(\bar u) &\dots & \pd{F^{m}}{u^m}(\bar u)  &                                    &   \dots              &  \\
  \pd{F^{m+1}}{u^1}(\bar u)&\dots & \pd{F^{m+1}}{u^{m}}(\bar u)                       &  \bar \lam^{m+1}  &                         &  \\
  \vdots                        & \vdots& \vdots                         &                        & \ddots               &  \\
  \pd{F^n}{u^1}(\bar u) &\dots &  \pd{F^n}{u^m}(\bar u) &                                     &                          &\bar\lam^n 
\end{array}
\right],
$$
where empty spaces a filled with zero. At the point $\bar u$, the matrix  $[D_uF]$ has $n$ distinct  real eigenvalues $\lam^1(\bar u),\dots, \lam^m(\bar u), \bar\lam^{m+1},\dots, \bar\lam^n$. Since the entries of   $[D_uF]$ are smooth real functions, a standard argument, involving the implicit function  theorem, implies that there is an open neighborhood $\Omega'\subset\Omega$ of $\bar u$,  such that at every point of $\Omega'$  the matrix $[D_uF]$ has $n$ distinct real eignevalues, and, therefore, $F$ is strictly hyperbolic on $\Omega'$.

\end{proof}

\begin{remark}[Single vector field case] \label{rem-single} \rm  When  $\R=\{\vr_1\}$, all three conditions \eqref{ic1} - \eqref{ic3} trivially hold. We, therefore, can assign $\lam^1$ to be any function on $\Omega$.  Then, by Proposition~\ref{prop-inv}, for every assignment  of the vector field $\tilde\vf$  on an $(n-1)$-dimensional 
 manifold  $\Xi$, transverse to $\vr_1$,  there exists  unique local vector field $\vf$  such $\tnabla_{\vr_1}\vf=\lambda^1 \,\vr_1$ and  $\vf|_{\Xi}=\tilde\vf|_{\Xi}$.   Thus the   general solution of the $\fs(\R)$-system \eqref{eq-fr} depends on one arbitrary function of $n$-variables (this is the functions $\lam^1$) and
 $n$ functions of $n-1$ variables, that locally describe the initial data  for  the vector field $\vf$.  Due to Proposition~\ref{prop-suff},  the $\fs(\R)$-set contains strictly hyperbolic fluxes.

\end{remark}

\begin{remark} [Full frame]\label{rem-full}  \rm If $\R$ is a full frame, it, of course, in involution.   In this case \eqref{ic3}, trivially holds  and the remaining equations,  \eqref{ic1} and  \eqref{ic2},
 comprise an algebro-differential system,  called the $\lam$-system and analyzed in details  in \cite{jk1}. According to Proposition~\ref{prop-inv},  for every solution of the $\lam$-system and for every assignment of the vector $\tilde\vf$ at a point $\bar u\in \Omega$, there exists a locally  unique solution $\vf$ of \eqref{eq-fr} such that $\vf|_{\bar u}$ is prescribed.   This can be also seen directly as follows. Since a full frame of eigenvector fields  are given, once eigenfunctions  are found, the Jacobian matrix $[D_uF]$ can be immediately obtained. The $i$-th row of $[D_uF]$ is the gradient of $F^i$, and $F^i$ itself can be recovered in the standard manner by solving a sequence of ODEs. If the value of $F^i(\bar u)$is prescribed, then functions $F^i$ is unique.
 
\end{remark}
 
\subsection{Rich partial frame}\label{sect-rich}
Rich frames comprise a  particularly nice subclass of involutive frames.   Recall that according to Definition~\ref{def-rich}, a partial frame $\R=\{\vr_1,\dots, \vr_m\}$, $1\leq m\leq n$,  is called rich, if it is pairwise in  involution: $[\vr_i,\vr_j]\in\Span\{\vr_i,\vr_j\}$ for all $i,j\in \{1,\dots,m\}$.   This case trivially includes all partial frames consisting of a single vector-field. Also this case includes all  involutive partial frames consisting of two vector fields. 

Let $\{\vr_1,\dots, \vr_m, \vs_{m+1},\dots,\vs_{n}\}$ to be any completion of $\R$ to a frame and let, as usual, use $c$ and $\Gamma$ to denote the corresponding structure functions and Christoffel symbols of the connection $\tnabla$, respectively. Since $\R$ is rich and  due to the symmetry of the connection, we have    
\beq \label{eq-c-gamma-rich} c_{ij}^l=0 \text{ and } \Gamma_{ij}^l=\Gamma_{ji}^l \text{ for all  distinct tripples $i,j,l$, such that $1\leq i,j\leq m$, $1\leq l\leq n$.}\eeq
Then the $\fs(\R)$-system   \eqref{ic1} -- \eqref{ic3} becomes
\begin{eqnarray}
\label{ic1r}\vr_i (\lambda^j)&=& \Gamma^j_{ji}\,(\lambda^i-\lambda^j)\text { for all  } 1\leq i\neq j\leq m\\
\label{ic2r} \Gamma_{ij}^l \,(\lambda^i-\lambda^j)&=&0\text { for all  } 1\leq i< j\leq m, \, 1\leq l\leq n \text{ such that } l\neq i\text{ and }l\neq j.
\end{eqnarray}

In the rich case, the necessary conditions for  the $\fs(\R)$-set to contain strictly hyperbolic fluxes, spelled out in the Proposition~\ref{prop-nec}, become
\beq\label{eq-rich-nec} \tnabla_{\vr_i}\vr_j\in \Span_{C^\infty(\Omega')}\{\vr_i,\vr_j\}\text{ for all }1\leq i\neq j\leq m.\eeq
Theorem~\ref{thm-rich-sh}  shows that, for a rich partial frame, this necessary conditions are also sufficient. Moreover, for the frames that satisfy \eqref{eq-rich-nec}, the proposition  describes the "size" of the set $\fs(\R)$. Theorem~\ref{thm-rich-nsh}  describes the "size" of the set $\fs(\R)$ for partial frames that do not satisfy \eqref{eq-rich-nec}, and therefore, do not admit strictly hyperbolic fluxes 

The following  lemma  allows  us to introduce a coordinate system adapted to a given rich partial frame and subsequently to invoke Darboux theorem to describe  the solution set of the $\fs(\R)$-system.
\begin{lemma} \label{lem-rich-coord}Assume a partial frame $\R=\{\vr_1,\dots,\vr_m\}$ on $\Omega$ is \emph{rich}, then in a neighborhood of every point $\bar u\in \Omega$ there exist
\begin{enumerate}
\item[1)] positive scalar functions $\alpha^1,\dots, \alpha^m$,  such that 
vector fields $\tilde \vr_i=\alpha^i\,\vr_i, \,\, i=1,\dots, n$ commute, i.e. $[\tilde \vr_i,\tilde\vr_j]=0$ for all $i,j\in\{1,\dots, m\}$;
\item[2)] local coordinate functions $(w^1,\dots, w^n)$,  such that   $\tilde \vr_i=\pd{}{w^i},\,\, i=1,\dots, m$.\end{enumerate}
\end{lemma}
\begin{proof} \begin{enumerate}
\item[1)] For a rich partial frame $\R$ the following structure equations hold:
$$[\vr_i,\vr_j]=c_{ij}^i\vr_i+c_{ij}^j\vr_j\quad i,j=1,\dots, m,$$
where structure functions $c_{ij}^k$ are independent of comletition of $\R$ to a frame.
We will show that the condition $[\tilde \vr_i,\tilde\vr_j]=0$ leads to a PDE system on $\alpha$'s of generalized Frobenius type. Indeed,
\begin{eqnarray}\nonumber [\tilde \vr_i,\tilde\vr_j]&=&[\alpha^i \vr_i,\alpha^j\vr_j]=\alpha^i\, \alpha^j\, [\vr_i,\vr_j]+\alpha^i \,\vr_i(\alpha^j)\vr_j -\alpha^j \,\vr_j(\alpha^i)\vr_i \\
\label{eq-alpha}&=& \alpha^j\left(\alpha^i c_{ij}^i-\vr_j(\alpha^i)\right)\vr_i - \alpha^i\left(\alpha^j c_{ji}^j-\vr_i(\alpha^j)\right)\vr_j.
\end{eqnarray}
Then $[\tilde \vr_i,\tilde\vr_j]=0$ if and only if $\beta^i=\ln(\alpha^i)$ satisfies the PDE system. 
\beq\label{eq-aij} \vr_j(\beta^i)=c_{ij}^i(u) \text{ for all }1\leq i\neq j\leq m  \eeq
To this system we add equations:
\beq\label{eq-aii} \vr_j(\beta^j)=0 \text{ for all }1\leq  j\leq m,  \eeq
making an additional requirement that, for each $i=1,\dots,m$, $\beta^i$ is constant along the integral curve of $r^i$. 
Since $c_{jj}^j=0$, we can combine \eqref{eq-aij} and \eqref{eq-aii}  in one system of $m^2$ equations on $m$ unknown functions $\beta$ of $n$ variables of generalized Frobenius type: 
\beq\label{eq-aall} \vr_j(\beta^i)=c_{ij}^i(u) \text{ for all }1\leq i, j\leq m  \eeq
We now write out the integrability conditions \eqref{eq-ic}, prescribed in Theorem~\ref{th-pde-frob}.  
\beq \vr_j(c^i_{ik})-\vr_k(c^i_{ij})=c_{jk}^j\, c^i_{ij}+c^k_{jk}\,c^i_{ik}  \text{ for all }1\leq i, j, k\leq m\eeq
and note that they are satisfied due to Jacobi  identities \eqref{c-jacobi}.

Due to Theorem~\ref{th-pde-frob}, we can prescribe any initial value for $\beta$'s along a submanifold transversal to $\R$  and get a unique solution of \eqref{eq-aall} on an open neighborhood of $\bar u$ with this initial data. Then positive functions $\alpha^i=e^{\beta^i}$ satisfy requirements of the theorem.
\item[2)] This is a direct consequence of Proposition~\ref{prop-comm-frame}.
\end{enumerate}
\end{proof}

Due to Lemma~\ref{lem-rich-coord}  and thanks to the  scaling invariance of Problems~\ref{pr} and~\ref{pr2}, we may assume that   the given rich partial frame is commutative.  
We then can use  a  local coordinate system $w^1,\dots,w^n$, such that    $\vr_i=\pd{}{w^i}$, for $i=1,\dots, m$. We complete $\R$ to a frame $\{\vr_1,\dots, \vr_m, \vs_{m+1},\dots,\vs_{n}\}$, where  $\vs_l=\pd{}{w^l}$, for $l=m+1,\dots, n$. 
The commutativity of the frame and  the symmetry of the connection $\tnabla$ imply the following conditions on the structure coefficients \eqref{eq-c} and Christoffel symbols \eqref{eq-gamma} for this frame:
\beq \label{eq-c-gamma-rich} c_{rs}^l=0 \text{ and } \Gamma_{rs}^l=\Gamma_{sr}^l \text{ for all  } l, s, r\in\{1,\dots,n\}.\eeq
 Then equations  \eqref{ic1} --
 \eqref{ic3}  become: 
\begin{eqnarray}
\label{ic1rw}\pd{}{w_i} (\lam^j)&=&  \Gamma^j_{ji}\,(\lambda^i-\lambda^j)\text { for all  } 1\leq i\neq j\leq m\\
\label{ic2rw} \Gamma_{ij}^l \,(\lambda^i-\lambda^j)&=&0\text { for all  } 1\leq i< j\leq m,  1\leq l\leq  n, \text{ such that } l\neq i\text{ and }l\neq j.
\end{eqnarray}
Assuming that the Christoffel symbols $\Gamma$'s and the unknown functions $\lam$'s are expressed in $w$-coordinates, we can treat \eqref{ic1rw} -- \eqref{ic2rw}, as a system of PDE's with simple linear constrains on the unknown functions $\lam$'s:

\begin{theorem}\label{thm-rich-sh} 
If a partial  frame $\R=\{\vr_1,\dots,\vr_m\}$ is rich  and it satisfies conditions \eqref{eq-rich-nec}, then the set $\fs(\R)$ of all local solutions  of  \eqref{eq-fr} near $\bar u$ depends on 
\begin{itemize}
\item $m$ arbitrary functions of $n - m + 1$ variables, prescribing, for $j=1,\dots, m$,   a function $\lam^j$ along an arbitrary $(n - m + 1)$-dimensional manifold $\Xi_j$ containing $\bar u$ and transverse to the set of vector-fields $\{\vr_1,\dots,\vr_{j-1},\vr_{j+1},\dots, \vr_m\}$;
\item $n$ functions of $n-m$ variables\footnote{Example~\ref{ex:basis} demonstrates that, when a general solution of an $\fs(\R)$-system is explicitly written out, some of the arbitrary functions of $n-m$ variables  may be absorbed into arbitrary functions of $n-m+1$ variables (a larger number of variables).  This is a standard phenomena  arising in applications of integrability theorems.}, prescribing components of a vector field $\vf$ along an arbitrary $(n - m)$-dimensional manifold $\Xi$ transverse to the partial frame $\R$.
\end{itemize}
The above data uniquely determines $\vf$ in an open  neighborhood  of $\bar u$. The $\fs(\R)$-set always contains strictly hyperbolic fluxes. 
\end{theorem}
\begin{proof}
\hfill
\begin{enumerate} 
\item[1)]As has been discussed above, after rescaling, we may assume that $\R$ is a commutative frame and we choose a coordinate system such that $\vr_i=\pd{}{w^i},\,\, i=1,\dots, m$.
Conditions \eqref{eq-rich-nec} are invariant under resaling of $\R$ and  imply that 
\beq \label{eq-gamma-distinct} \Gamma_{ij}^l\equiv 0 \text{ for all } 1\leq i\neq j\leq m,  1\leq l\leq  n, \text{ such that } l\neq i\text{ and }l\neq j,\eeq
and therefore equations \eqref{ic2rw}  trivially hold. Equations \eqref{ic1rw} is of the Darboux type and we proceed by verifying the  integrability conditions (C) stated in Theorem~\ref{dar3}. For this purpose we substitute partial derivatives prescribed by  \eqref{ic1rw}, into equality of mixed partials conditions condition: 
$$ \pd{}{w^k}\left(\pd{\lam^j}{w^i}\right)\equiv\pd{}{w^i}\left(\pd{\lam^j}{w^k}\right), \text {for all distinct triples $i,j, k\in \{1,\dots,m\}$}.$$
The first substitution leads to
$$ \pd{}{w^k}\left( \Gamma^j_{ji}\,(\lambda^i-\lambda^j)\right)\equiv \pd{}{w^i}\left( \Gamma^j_{jk}\,(\lambda^k-\lambda^j)\right), \text {for all distinct triples $i,j, k\in \{1,\dots,m\}$},$$
and the subsequent substitution (using abbreviated  notation  $\del_i=\pd{}{w^i}$) leads to the condition:
\begin{eqnarray}\label{int-k0}  
	\left(\del_{i}\Gamma^{j}_{jk} -\del_{k}\Gamma^{j}_{ji}\right)\lam^j
	&+&\left(\Gamma^{j}_{ji}\Gamma^{i}_{ik}+\Gamma^{j}_{jk}\Gamma_{ki}^{k}-\Gamma^{j}_{ji}\Gamma^{j}_{jk}-\del_{i}\Gamma^{j}_{jk}\right)\lam^{k}\\
	\nonumber&-&\left(\Gamma^{j}_{ji}\Gamma^{i}_{ik}+\Gamma^{j}_{jk}\Gamma_{ki}^{k} - \Gamma^{j}_{jk}\Gamma^{j}_{ji}
	-\del_{k}\Gamma^{j}_{ji}\right)\lam^i\equiv 0\,.
\end{eqnarray}
which must hold  for all  triples of pairwise distinct indices $i,j, k\in \{1,\dots,m\}$.
We will use   flatness condition \eqref{R0-gamma} to show that all $\lam$'s appear \eqref{int-k0}  with identically zero coefficients. 

We first substitute $s=j$  in  \eqref{R0-gamma} and  we assume that $i,j, k\in \{1,\dots,m\}$  are pairwise distinct indices. Then  using \eqref{eq-c-gamma-rich} and \eqref{eq-gamma-distinct}, we obtain that  for all  triples of pairwise distinct indices $i,j, k\in \{1,\dots,m\}$:
\beq
\label{eq-R01} -\del_k\Gamma^j_{ji}
	= \Gamma^j_{jk}\Gamma^j_{ji}-\Gamma^j_{ji}\Gamma^i_{ik}-\Gamma^j_{jk}\Gamma^k_{ki}.
\eeq
This immediately implies that the coefficient, $\Gamma^{j}_{ji}\Gamma^{i}_{ik}+\Gamma^{j}_{jk}\Gamma_{ki}^{k} - \Gamma^{j}_{jk}\Gamma^{j}_{ji}-\del_{k}\Gamma^{j}_{ji}$, of  $\lam^i$ in \eqref{int-k0} is identically zero.
Interchanging $k$ and $i$ in \eqref{eq-R01}, we obtain:
\beq
\label{eq-R02} -\del_i\Gamma^j_{jk}
	= \Gamma^j_{ji}\Gamma^j_{jk}-\Gamma^j_{jk}\Gamma^k_{ki}-\Gamma^j_{ji}\Gamma^i_{ik}.
\eeq
and, therefore, the coefficient of  $\lam^k$ in \eqref{int-k0} is identically zero.
We note that the right-hand sides of the identities  \eqref{eq-R01} and \eqref{eq-R02}  are equal and, therefore, the coefficient, $\del_{i}\Gamma^{j}_{jk} -\del_{k}\Gamma^{j}_{ji}$,  of $\lam^j$ in \eqref{int-k0} is identically zero.

Thus, we have verified the  integrability conditions (C) stated in Theorem~\ref{dar3} do hold for the PDE system  \eqref{ic1rw}. We conclude that,
for a fixed point $\bar u\in \Omega$, whose $u$-coordinates are $(\bar u^1,\dots,\bar u^n)$ and $w$-coordinates are $(\bar w^1,\dots,\bar w^n)$  and  any assignment of $m$ arbitrary functions of $n-m+1$ variables: 
 $$\tilde \lam ^i(\bar w^1,\dots, \bar w^{i-1}, w^i, \bar w^{i+1},\dots, \bar w^m, w^{m+1},\dots, w^n), \, i=1,\dots, m$$   on the subsets  $\Xi_i\subset \Omega$, where  $w^j=\bar w^j$, for $j=1,\,\dots,\, i-1, \, i+1,\,\dots,\, m$, there is a   unique local solution $\lambda^1,\dots, \lam^m$ of \eqref{ic1rw},  such that $\lam^i|_{\Xi_i\cap\Omega'}=\tilde \lam^i|_{\Xi_i\cap\Omega'}$ on some open subset $\Omega'\subset\Omega$  containing $\bar u$.
 Thus  the general solution $\lam$  of  \eqref{ic1rw} depends on $m$ arbitrary functions of $n - m + 1$ variables.

 \item[2)] \label{it2} Recalling that  for a  rich frame, satisfying \eq{eq-rich-nec}, the system  \eqref{ic1r} is equivalent to the $\lam$-system \eqref{ic1} -- \eqref{ic3}, we use   Proposition~\ref{prop-inv} to conclude  that  for any solution $\lam$  of  \eqref{ic1r} and  any smooth initial data for $\vf$ prescribed  along any   embedded submanifold $\Xi\subset \Omega$ of codimension $m$ transversal to $\R$,
there is a unique smooth local solution of $\fs(\R)$-system \eqref{eq-fr}. In local coordinates, the initial data can be  defined by $n$  functions (components of $\vf$) of $n-m$ variables (local coordinates on $\Xi$). Therefore, for a given solution $\lam$ of  \eqref{ic1r}, the general solution $\vf$ of      $\fs(\R)$-system \eqref{eq-fr} depends on $n$ arbitrary  functions of $n-m$ variables.

\item[3)]\label{it3}   We can always choose $\tilde\lam^1,\dots, \tilde\lam^m$ in  
Part~1)
 of the proof, such that all $m$ real numbers  
 $\tilde\lam^1(\bar u),\dots,\tilde\lam^m(\bar u)$ are distinct. Let    
 $\lam^1,\dots, \lam^m$ be the corresponding solutions of \eqref{ic1r}.  Then the existence of strictly hyperbolic fluxes in the $\fs(\R)$-set  follows from Proposition~\ref{prop-suff}. 
\end{enumerate}
 \end{proof}
We observe that  in single vector field case ($m=1$), the conclusion of Theorem~\ref{thm-rich-sh}  is consistent with the observation made in Remark~\ref{rem-single}. 
The first part of  the proof of  Theorem~\ref{thm-rich-sh} is  a rather straightforward generalization of the   proof of Theorem~4.3 in \cite{jk1}, where  the $\lam$-system \eqref{ic1r}  was considered in the case of the full frame ($m=n$).
 In a similar way, we can  generalize Theorem~4.4 in \cite{jk1} to treat the case when necessary conditions \eqref{eq-rich-nec} for strict hyperbolicity are not satisfied. In this case,  the algebraic relationship  \eqref{ic2r} implies that there exist $i,j\in\{1,\dots,m\}$, such that $i\neq j$ and $\lam^i\equiv \lam^j$, and therefore,   there are no strictly hyperbolic fluxes in the $\fs(\R)$-set. A rather involved description of the $\fs(\R)$-set is given by the following theorem, whose proof can be easily spelled out by combining the arguments in the proofs of Theorem~4.4 in \cite{jk1} and Theorem~\ref{thm-rich-sh} above. The argument is rather technical and is not reproduced here.
\begin{theorem}\label{thm-rich-nsh} 
Let $\R=\{\vr_1,\dots,\vr_m\}$ be a  rich  partial frame on an open subset $\Omega\subset\RR^n$, that \emph{does not} satisfy conditions \eqref{eq-rich-nec}.   Then
the system \eq{ic1r}  -- \eq{ic2r} imposes multiplicity conditions\footnote{ It is clear that for all $i\neq j$, such that  $\tnabla_{\vr_i}\vr_j\notin\Span\{\vr_i,\vr_j\}$, equations \eq{ic2r} imply a multiplicity condition $\lam^i=\lam^j$. Less obviously,  \eq{ic1r} may impose additional multiplicity conditions on $\lam$'s. See the proof of  Lemma~4.5 in \cite{jk1} for more details.} 
	on $\lam$'s in the following sense. There are disjoint subsets 
	$A_1, \dots, A_{s_0}\subset \{1,\dots,m\}$ ($s_0\geq1$) of cardinality two or more, and such 
	that \eq{ic1r}  -- \eq{ic2r} impose the equality $\lam^i=\lam^j$ if and only if 
	$i,\, j\in A_\alpha$ for some $\alpha\in\{1,\dots,s_0\}$. 
	Let $l=\sum_{\alpha=1}^{s_0} |A_\alpha|\leq m$ and $s_1=m-l$. By relabeling indices we 
	may assume that $\{1,\dots,m\}\setminus\bigcup_{\alpha=1}^{s_0}A_\alpha=\{1,\dots,s_1\}$.

 The set $\fs(\R)$ of all local solutions  of  \eqref{eq-fr} near $\bar u$ depends on 
\begin{itemize}
\item $s_1$ arbitrary functions $\tilde \lam^1,\dots,\tilde \lam^{s_1}$ of $n - m + 1$ variables, prescribing, for $j=1,\dots, s_1$,   the data  for  function $\lam^j$, so that  $\lam^j|_{\Xi_j}=\tilde\lambda^j$, where $\Xi_j$ is an arbitrary $(n - m + 1)$-dimensional manifold $\Xi_j$ containing $\bar u$ and  transverse to the set of vector-fields $\{\vr_1,\dots,\vr_{j-1},\vr_{j+1},\dots, \vr_m\}$;
\item $s_0$ arbitrary functions  $\kappa^1,\dots, \kappa^{s_0}$ of $m-n$ variables,  prescribing,   for  $j=s_1+1,\dots, m$, the data  for functions $\lam^j$, so that when   $j\in A_\alpha$ for some $\alpha =1,\dots, s_0$  when   $j\in A_\alpha$ for some $\alpha =1,\dots, s_0$, then  $\lam^j|_{\Xi_j}=\kappa^\alpha$, where $\Xi_j$ is an $(n-m)$-dimensional manifold passing through $\bar u$ and transverse to $\R$;
\item $n$ functions of $n-m$ variables prescribing components of a vector field $\vf$ along an arbitrary $(n - m)$-dimensional manifold $\Xi$ transverse to the partial frame $\R$.
\end{itemize}
The above data uniquely determines $\vf$ in an open  neighborhood  of $\bar u$. The $\fs(\R)$-set never contains strictly hyperbolic fluxes. \end{theorem}


\subsection{Non-rich involutive frames consisting of three vector-fields}

The lowest cardinality of a partial frame, for which involutive, non-rich scenario may  appear, is $m=3$ case. In \cite{jk1}, we treated the case  when $m=n=3$, i.e. the full frame case.
We now generalize these results to  $n\geq 3$. Generalization to  $m>3$  would require a consideration of a large number of cases and was not performed here.

 We first treat the case when $\R$ satisfies the necessary conditions of Proposition~\ref{prop-nec}  for the existence of strictly hyperbolic fluxes.  We choose an arbitrary completion of $\R$ to a frame and write out the $\lam$-system \eq{ic1} -- \eq{ic3}. The differential part \eq{ic1}
becomes:
\bea
	\nn                  \vr_2(\lam^1) &=& \Gamma_{12}^1(\lam^2-\lam^1)\\
	\nn                   \vr_3(\lam^1) &=& \Gamma_{13}^1(\lam^3-\lam^1)\\
	\label{lam-m3} \vr_1(\lam^2) &=& \Gamma_{21}^2(\lam^1-\lam^2)\\
         \nn                   \vr_3(\lam^2) &=& \Gamma_{23}^2(\lam^3-\lam^2)\\
         \nn                   \vr_1(\lam^3) &=& \Gamma_{31}^3(\lam^1-\lam^3)\\ 
	\nn                   \vr_2(\lam^3) &=& \Gamma_{32}^3(\lam^2-\lam^3).
\eea
Algebraic equations \eq{ic2} can be written as:
\beq\label{eq-alg-lam}
		A_\lambda\, \left[\begin{array}{c}\lam^1\\ \lam^2\\ \lam^3\end{array}\right]=0\,, \qquad\text{where}\qquad A_\lam=
		\left[\begin{array}{ccc}
		c_{23}^1 &  \Gamma_{32}^1 & -\Gamma_{23}^1 \\
	  	\Gamma_{31}^2&  c_{13}^2 &    - \Gamma_{13}^2\\
	    	\Gamma_{21}^3&-\Gamma_{12}^3   &c_{12}^3   
		\end{array}\right].
	\eeq
Condition \eq{eq-nec-sh3}  in Proposition~\ref{prop-nec} implies that \eq{ic3} is trivial.  We also  note that, since $\R$ is involutive and satisfies conditions in Proposition~\ref{prop-nec},  for all $i,j,k\in\{1,2,3\}$, the structure coefficients $c_{ij}^k$ and Chrisoffel symbols $\Gamma_{ij}^k$ are independent of the completion   of $\R$ to a frame, and therefore the system    \eq{lam-m3} -- \eq{eq-alg-lam} can be written out without specifying a completion to a full frame.  
Our goal is to  describe the solution set of \eq{lam-m3} -- \eq{eq-alg-lam}.

Looking more closely at matrix $A_\lam$ we make the following observations
\begin{itemize}
\item  From the symmetry of the connection it follows that the last column of $A_\lambda$ is the sum of the first two columns and therefore $\rank A_\lam\leq 2$.  
\item Non-richness of $\R$ implies that at least one of $c$'s appearing in  $A_\lambda$ is non zero and therefore $\rank A_\lam\geq 1$.  
\item Condition \eq{eq-nec-sh2}  in Proposition~\ref{prop-nec} implies that, for each row in $ A_\lam$, either all three entries are zero, or all three entries are non-zero. 
\end{itemize}
Following the same argument as in Section~3 of \cite{jk1}, one can show that if $\rank A_\lam=2$ at $\bar u$, then the three eigenfunctions must coincide in a neighborhood of $\bar u$, i.e.~ $\lam^1=\lam^2=\lam^3=\lam$ for some functions $\lam$, and, therefore, $\fs(\R)$ does not contain strictly hyperbolic fluxes. Moreover,
\eq{lam-m3}  imply that $\lam$ is constant along the integral manifolds of the involutive frame $\R$, and we can prescribe an arbitrary  value of $\lam$ along a manifold $\Xi$ transverse to $\R$.  Otherwise,  $\rank A_\lam=1$, and we may assume without loss of generality, that   $c_{23}^1\neq 0$. The  first equation in \eq{eq-alg-lam} can be solved for $\lam^1$ and this solution can be substituted in \eq{lam-m3}. After simplifications we get a system that  specifies the derivatives of the two unknown functions $\lam^2$ and $\lam^3$ on $\R^n$ along
a partial  involutive frame $\vr_1$, $\vr_2$ and $\vr_3$:
\begin{eqnarray}
	\nn   \vr _1(\lam^2) &=& \frac{\Gamma_{21}^2\Gamma_{23}^1}{c_{32}^1}(\lam^2-\lam^3)\,,\\
	      \nonumber        \vr_2(\lam^2) &= &\left[\frac{\Gamma_{23}^1}{\Gamma_{32}^1}(\Gamma_{32}^3-\Gamma_{12}^1)
	                                                             -\frac{c_{32}^1}{\Gamma_{32}^1}\,r_2\left(\frac{\Gamma_{32}^1}{c_{32}^1}\right)\right](\lam^2-\lam^3)\,,\\
	     \label{IIa-solved}         \vr_3(\lam^2)& =&- \Gamma_{23}^2(\lam^2-\lam^3)\,,\\
	      \nonumber        \vr_1(\lam^3) &=& \frac{\Gamma_{31}^3\Gamma_{32}^1}{c_{32}^1}(\lam^2-\lam^3)\,,\\
	      \nonumber       \vr_2(\lam^3)& =& \Gamma_{32}^3(\lam^2-\lam^3)\,,\\
	      \nonumber       \vr_3(\lam^3) &= & \left[\frac{\Gamma_{32}^1}{\Gamma_{23}^1}(\Gamma_{13}^1-\Gamma_{23}^2)
	+\frac{c_{32}^1}{\Gamma_{23}^1}\,r_3\left(\frac{\Gamma_{23}^1}{c_{32}^1}\right)\right](\lam^2-\lam^3)\,,\quad
	\end{eqnarray}
This system looks identical to the system  (3.22) in \cite{jk1}, however, in \cite{jk1}, we had $n=3$, while here $n\geq 3$ and, therefore, the classical  Frobenius theorem, used in \cite{jk1}, is not sufficient in this case, and, therefore, we appeal to a more general Theorem~\ref{th-pde-frob}.  To verify the integrability conditions   we rewrite  \eq{IIa-solved}  
\beq\label{simple}
	r_i(\lambda^s)=\phi_i^s(u)(\lambda^2-\lambda^3)\qquad \mbox{for $i=1,2,3$ and  $s=2, 3,$}
\eeq
where $\phi_i^s$ are known functions of $\Gamma$'s, given by the right-hand sides in \eq{IIa-solved}. 
Then the  integrability conditions   amount to:
\beq\label{comp2}
	\Big[r_i(\phi_j^s)-r_j(\phi_i^s)+\phi_j^s(\phi_i^2-\phi_i^3)-\phi_i^s(\phi_j^2-\phi_j^3)\Big](\lam^2-\lam^3)=
	\left[\sum_{k=1}^3c^k_{ij}\phi_k^s\right](\lam^2-\lam^3)\,, 
\eeq
where $1\leq i<j\leq 3$, $s=2, 3$ and $c_{ij}^k=\Gamma_{ij}^k-\Gamma_{ji}^k$.

These conditions are satisfied if $\lam^2=\lam^3$ in a neighborhood of $\bar u$, in which case, the  first equation in \eq{eq-alg-lam} implies  $\lam^1=\lam^2=\lam^3=\lam$, and, as above, the functions  $\lam$ must be constant along the integral manifolds of the involutive frame $\R$, and we can prescribe an arbitrary  value of $\lam$ along a manifold $\Xi$ transverse to $\R$.    
For a strictly  hyperbolic flux to exist the following six conditions must hold:
\begin{eqnarray}
\label{frob-comp1}	\vr_i(\phi_j^2)-\vr_j(\phi_i^2)
	&=&\phi_j^2\phi_i^3-\phi_i^2\phi_j^3+\sum_{k=1}^3c^k_{ij}\phi_k^2\qquad 1\leq i<j\leq 3,\label{s=2} \\
	\label{frob-comp2}\vr_i(\phi_j^3)-\vr_j(\phi_i^3)&=&
	\phi_j^2\phi_i^3-\phi_i^2\phi_j^3+\sum_{k=1}^3c^k_{ij}\phi_k^3\qquad 1\leq i<j\leq 3.\label{s=3} 
 \end{eqnarray}

Conditions~\eq{frob-comp1} -- \eq{frob-comp2}, in the case of full frames in $\RR^3$, were derived in \cite{jk1}, and  Examples~5.1~and~5.3 in \cite{jk1} show that these compatibility conditions may or may 
not be satisfied: they must be checked for each case individually. 
If these integrability conditions are met then, according to Theorem~\ref{th-pde-frob}, the general solution
to the $\lambda$-system depends on two functions of $n-3$ variables prescribing the values of $\lam^2$ and $\lam^3$ along any two $n-3$ dimensional manifold passing through $\bar u$ and transverse  to $\R$. 
Function $\lam^1$ is then determined by the first equation in \eq{eq-alg-lam}. Combining the above argument  with  Propositions~\ref{prop-inv} and  Propositions~\ref{prop-nec}  we arrive to the following theorem:
\begin{theorem}\label{thm-iia-gen} Assume $\R=\{\vr_1,\vr_2,\vr_3\}$ is a non-rich partial frame in involution, on a neighborhood $\Omega$  of $\bar u$, satisfying conditions \eq{eq-nec-sh2} and \eq{eq-nec-sh3}  in Proposition~\ref{prop-nec}. For $i,j,k \in \{1,2,3\}$, let $c_{ij}^k$ and $\Gamma_{ij}^k$ be defined by
$$ [\vr_i,\vr_j]=\sum_{k=1}^3c_{ij}^k\vr_k \qquad  \nabla_{\vr_i}\vr_j=\sum_{k=1}^3\Gamma_{ij}^k\vr_k.$$
Up to permutation of indices and by shrinking $\Omega$ we may assume $c_{23}^1$ is nowhere zero on $\Omega$.
\begin{itemize}
\item  If the matrix $A_\lam$ defined in \eqref{eq-alg-lam} has rank 1 and that \eq{frob-comp1} -  \eq{frob-comp2} are satisfied in a neighborhood of $\bar u$,
then the solution set $\fs(\R)$ of system \eqref{eq-fr} depends on $n+2$ arbitrary functions of $n - 3$ variables (2 of those determine the  values  $\lam^2$ and $\lam^3$, while $n$  of those determine the values  $\vf$ along an $(n-3)$-dimensional manifold passing through $\bar u$ and transverse  to $\R$). The set   $\fs(\R)$ contains strictly hyperbolic fluxes. 
\item  If the matrix $A_\lam$ defined in \eqref{eq-alg-lam} has rank 2 at $\bar u$  or  \eq{frob-comp1} -  \eq{frob-comp2} are not satisfied at $\bar u$, then  then the three eigenfunctions must coincide in a neighborhood of $\bar u$, i.e.~ $\lam^1=\lam^2=\lam^3=\lam$ for some functions $\lam$, such that  $\lam$ is constant along the integral manifolds of the involutive frame $\R$, and can take arbitrary  values along a manifold $\Xi$ transverse to $\R$.  The solution set $\fs(\R)$ of system \eqref{eq-fr} depends on $n+1$ arbitrary functions of $n - 3$ variables (1 of those determine the  values  $\lam$ and $n$  of those determine the values  $\vf$ along an $(n-3)$-dimensional manifold passing through $\bar u$ and transverse  to $\R$). The set   $\fs(\R)$ does not contain strictly hyperbolic fluxes. \end{itemize}
\end{theorem}
When  the partial frame $\R$ \emph{does not} satisfy the necessary conditions of Proposition~\ref{prop-nec}  for the existence of strictly hyperbolic fluxes, then the algebraic conditions \eq{ic2}  and \eq{ic3} force two or more of eigenfunctions to be equal to each other, and we can prove  the following result:
\begin{theorem}\label{thm-iib-gen} Assume $\R=\{\vr_1,\vr_2,\vr_3\}$ is a non-rich partial frame  in involution, on a neighborhood $\Omega$  of $\bar u$, such that $\R$ \emph{does not} satisfy condition \eq{eq-nec-sh2} or condition  \eq{eq-nec-sh3}  in Proposition~\ref{prop-nec}. Then there are exactly two possibilities:
\begin{description}
\item{either}  the  $\lam$-system \eq{ic1}  -- \eq{ic3}  implies that  $\lam^1=\lam^2=\lam^3=\lam$, in a neighborhood of $\bar u$,where a function $\lam$ is constant along the integral manifolds of the involutive frame $\R$ and  may take arbitrary values on an $(n-3)$-dimensional manifold  $\Xi_0$ passing through $\bar u$ and transverse  to $\R$.
\item{or}, up to permutation of indices,   the  $\lam$-system \eq{ic1}  -- \eq{ic3}    implies that  $\lam^1=\lam^2=\lam$, but allows the possibility that $\lam\neq \lam^3$ in a neighborhood of $\bar u$. In this case, the  function $\lam^3$ is uniquely determined by  its values on an $(n-2)$-dimensional manifold  $\Xi_1$ passing through $\bar u$ and transverse  to $\{\vr_1,\vr_2\}$ and the  function $\lam$ is uniquely determined by  its values on an $(n-3)$-dimensional manifold  $\Xi_2$ passing through $\bar u$ and transverse  to $\R$, and.
\end{description}
In both cases, the $\lam$-system   \eq{ic1}  -- \eq{ic3} has a locally unique solution with the data, described above, and for each such solution,  the $\fs(\R)$-system \eqref{eq-fr} has a locally unique solution determined by  the values of $\vf$ on an $(n-3)$-dimensional manifold  $\Xi$ passing through $\bar u$ and transverse  to $\R$. 
The set   $\fs(\R)$  contains \emph{no} strictly hyperbolic fluxes. 

\end{theorem} 

\begin{proof}\hfill
\begin{enumerate}
\item[1)]  If condition \eq{eq-nec-sh2} is not satisfied, then equations \eq{ic2}  imply that at least two functions among $\lam^1$, $\lam^2$  and $\lam^3$ are identically equal to each other in neighborhood of $\bar u$.
Similarly, if  condition \eq{eq-nec-sh3} is not satisfied than equations \eq{ic3} imply that at least two functions among $\lam^1$, $\lam^2$  and $\lam^3$  coincide  in neighborhood of $\bar u$. In either case the set   $\fs(\R)$  does not contain strictly hyperbolic fluxes.
\item[2)] If  \eq{ic2} and  \eq{ic3} imply that all three are equal, i.~e.~$\lam^1=\lam^2=\lam^3=\lam$, then, the differential part \eq{ic1} of the $\lam$-system, 
  implies that  the function $\lam$ is constant along the integral manifolds of the involutive frame $\R$.  In this case, the system  \eq{ic1} trivially satisfies the assumptions of Theorem~\ref{th-pde-frob}, which implies that for any assignment of $\lam$ along   an $(n-3)$-dimensional manifold  $\Xi_0$ passing through $\bar u$ and transverse  to $\R$, there is unique such function in a neighborhood of $\bar u$.
\item[3)]   If \eq{ic2}  and \eq{ic3}, imply that only two of $\lam$'s coincide, e.g. $\lam^1=\lam^2=\lam$, but they don't imply that they must be  equal to $\lam^3$, then one can argue that $c$'s and $\Gamma$ satisfy  the following conditions
\beq\label{cond-1=2}  c^3_{12}=0, \, \Gamma^2_{13}=0 \text{ and }  \Gamma^1_{23}=0,\eeq
and  the $\lam$-system \eq{ic1} --  \eq{ic3} becomes:
\bea
	\nn                  \vr_2(\lam) &=&0\\
	\nn                   \vr_3(\lam) &=& \Gamma_{13}^1(\lam^3-\lam)\\
	\label{lam-IIb} \vr_1(\lam) &=& 0\\
         \nn                   \vr_3(\lam) &=& \Gamma_{23}^2(\lam^3-\lam)\\
         \nn                   \vr_1(\lam^3) &=& \Gamma_{31}^3(\lam-\lam^3)\\ 
	\nn                   \vr_2(\lam^3) &=& \Gamma_{32}^3(\lam-\lam^3).
\eea
If   $\Gamma_{23}^2\neq \Gamma_{13}^1$, the the second and the fourth equations in the above system imply that $\lam=\lam^3$, and therefore again  $\lam^1=\lam^2=\lam^3=\lam$, and we arrive to the situation considered in part 2) of the proof.  If   
\beq\label{cond-2neq3} \Gamma_{23}^2= \Gamma_{13}^1\eeq we end up with the system
\bea
	\label{lam-IIb-1} \vr_1(\lam) &=& 0\\
	\label{lam-IIb-2}            \vr_2(\lam) &=&0\\
	\label{lam-IIb-3}  \vr_3(\lam) &=& \Gamma_{13}^1(\lam^3-\lam)\\
	\label{lam-IIb-4}          \vr_1(\lam^3) &=& \Gamma_{31}^3(\lam-\lam^3)\\ 
	\label{lam-IIb-5}               \vr_2(\lam^3) &=& \Gamma_{32}^3(\lam-\lam^3).
\eea
We subtract  equations \eq{lam-IIb-1} from \eq{lam-IIb-4}, and equation  \eq{lam-IIb-2} from \eq{lam-IIb-5}, and introduce a new unknown functions $\mu=\lam^3-\lam$.
We obtain:
\bea
	\label{lam-IIb-1} \vr_1(\lam) &=& 0\\
	\label{lam-IIb-2}            \vr_2(\lam) &=&0\\
	\label{lam-IIb-3}  \vr_3(\lam) &=& \Gamma_{13}^1\,\mu\\
	\label{lam-IIb-4}          \vr_1(\mu) &=&- \Gamma_{31}^3\,\mu \\ 
	\label{lam-IIb-5}               \vr_2(\mu) &=&- \Gamma_{32}^3\,\mu.
\eea
By assumption $\{\vr_1,\vr_2,\vr_3\}$ are in involution, the  first condition in \eq{cond-1=2} implies that the vector fields $\vr_1$ and $\vr_2$ are in involution. 
Thus we can first apply Theorem~\ref{th-pde-frob} to the sub-system \eq{lam-IIb-4} - \eq{lam-IIb-5}, whose integrability condition,
\beq \label{icIIb-3}  \vr_2(\Gamma^3_{31})-\vr_1(\Gamma^3_{32})
	= c_{21}^2\, \Gamma_{32}^3+ c_{21}^1\,\Gamma_{31}^3\eeq
is satisfied as shown   in Lemma~3.6 of \cite{jk1}, due to the flatness and symmetry property of the connection, 
combined with conditions \eq{cond-1=2} and \eq{cond-2neq3}. 
Thus there is unique solution $\mu$ for the subsystem  \eq{lam-IIb-4} - \eq{lam-IIb-5} with any data prescribed along an $(n-2)$-dimensional manifold $\Xi_1$  passing through $\bar u$ and  transversal to $\vr_1,\vr_2$. 
Any solution $\mu$ can be substituted into  \eq{lam-IIb-3}, and then we apply Theorem~\ref{th-pde-frob} to the sub-system \eq{lam-IIb-1} -- \eq{lam-IIb-3}, whose integrability condition
\begin{align}
	\label{icIIb-12} \vr_2(\Gamma^1_{13})&=\Gamma_{23}^3\, \Gamma_{13}^1\\
	\nn \vr_1(\Gamma^1_{13})&=\Gamma_{13}^3\, \Gamma_{13}^1
\end{align}
As it is shown in Lemma~3.6 of \cite{jk1}, conditions \eqref{icIIb-12} hold identically 
on $\Omega$ due to the flatness and symmetry property of the connection, 
combined with conditions \eq{cond-1=2} and \eq{cond-2neq3}. 
Then Theorem~\ref{th-pde-frob} guarantees that there exists a locally unique solution of system sub-system \eq{lam-IIb-1} -- \eq{lam-IIb-3}, with the  values of function $\lam$ prescribed along an $(n-3)$-dimensional manifold  $\Xi_2$ passing through $\bar u$ and transverse  to $\R$.  Recalling that  $\mu=\lam^3-\lam$, we conclude that  $\lam$ is uniquely determined by  its values on an $(n-3)$-dimensional manifold  $\Xi_1$ passing through $\bar u$ and transverse  to $\R$, and  function $\lam^3$ is uniquely determined by  its values on an $(n-2)$-dimensional manifold  $\Xi_2$ passing through $\bar u$ and transverse  to $\{\vr_1,\vr_2\}$.
\item[4)] It follows from Proposition~\ref{prop-inv} that for each solution of the $\lam$ system,   $\fs(\R)$-system \eqref{eq-fr} has a locally unique solution determined by  the values of $\vf$ on an $(n-3)$-dimensional manifold  $\Xi$ passing through $\bar u$ and transverse  to $\R$. 
\end{enumerate}
\end{proof}

\section{Non-involutive  partial frames  of two vector fields in $\RR^3$.}\label{sect-noninvol}
In  the non-involutive case, the differential consequences  \eqref{gc1} -- \eqref{gc3} of the $\fs(\R)$-system \eq{eq-fr} involve additional functions $\aaa$'s, so instead of the ``$\lam$-system", we get the ``$\lam$-$a$-system", and, moreover,   \eqref{gc1} -- \eqref{gc3} do not provide a complete set of the integrability conditions for the $\fs(\R)$-system. This makes the non-involutive case  to be  much harder to analyze  than  the involutive case,  and we are able to treat  only the lowest  dimension where such scenario can arise: $\R=\{\vr_1,\vr_2\}$ is a partial frame in $\RR^3$, such that at a fixed  point $\bar u\in \Omega$:
\beq\label{eq-non-inv} [\vr_1,\vr_2]_{\bar u}\notin\Span_{\RR}\{\vr_1|_{\bar u}, \vr_2|_{\bar u}\}.\eeq
The $\fs(\R)$-system then consists of two equations:
\beq\label{FR2}\tnabla_{\vr_1} \,\vf = \lambda^1 \, \vr_1\text{ and } \tnabla_{\vr_2} \,\vf = \lambda^2 \, \vr_2
\eeq
and the necessary conditions \eqref{eq-nec-shg} for  strict hyperbolicty become
\beq\label{eq-nec-sh-m2} \tnabla_{\vr_1}\vr_2|_{\bar u} \notin \Span_{\RR}\{\vr_1|_{\bar u}, \vr_2|_{\bar u}\}  \text{ and }\tnabla_{\vr_2}\vr_1|_{\bar u} \notin \Span_{\RR}\{\vr_1|_{\bar u}, \vr_2|_{\bar u}\}.
 \eeq
Below we state two theorems that describe the size and the structure of the flux space $\fs(\R)$ for partial frames $\R$ satisfying \eqref{eq-nec-sh-m2}. The proofs of the Theorems~\ref{thm-noninv-str} and \ref{thm-noninv-size}  rely on the  sequences of lemmas listed below.
We remind the reader that $\fstriv$ denotes the $4$-dimensional space of trivial fluxes.
\begin{theorem}\label{thm-noninv-str} Let $\R=\{\vr_1,\vr_2\} $  a non-involutive partial frame on an open neighborhood of $\bar u\in \RR^3$ satisfying conditions \eqref{eq-nec-sh-m2}. Then 
\begin{enumerate}
\item[1)] 
A non-zero flux $\vf\in \fs(\R)/\fstriv$ is either strictly hyperbolic or non-hyperbolic.
\item [2)]  If $\dim \fs(\R)/\fstriv>1$, then $\fs(\R)$ contains strictly hyperbolic fluxes.
\item  [3)]  If $\fs(\R)$ contains a non-hyperbolic flux, then for  any vector field $\vs$ completing  $\R$ to a local frame, the following identity holds on an open neighborhood of $\bar u$:
\beq\label{eq-nec-nonhyp} \Gamma_{12}^3\,\Gamma_{21}^3-2\,(c_{12}^3)^2=\Gamma_{11}^3\,\Gamma_{22}^3, \eeq
where $c$'s and $\Gamma$'s are structure components and Christoffel symbols for connection $\tnabla$ relative to the frame $\vr_1,\vr_2,\vs$.

\end{enumerate}
\end{theorem}
Although identity \eq{eq-nec-nonhyp} is a closed condition, and, therefore, is restrictive, Examples~\ref{ex:3} and~\ref{ex:4b} demonstrate that there are partial frames whose set of fluxes contains non-hyperbolic fluxes. On the other hand, Examples~\ref{ex:1}, \ref{ex:2}, \ref{ex:4a} and~\ref{ex-no-mist}  show  that there are partial frames  for which all non-trivial fluxes are strictly hyperbolic.

\begin{theorem}\label{thm-noninv-size}  Let $\R=\{\vr_1,\vr_2\} $ be a non-involutive partial frame on an open neighborhood of $\bar u\in \RR^3$ satisfying conditions \eqref{eq-nec-sh-m2}.   Let $\vs$ be any completion of $\vr_1$ and $\vr_2$ to  a local frame near $\bar u$ and let $\Gamma$'s be Christoffel symbols for connection $\tnabla$ relative to this frame. Assume further that the following condition is satisfied:
\beq\label{eq-mist}\Gamma^3_{22} (\bar u)\,\Gamma^3_{11}(\bar u) -9\,\Gamma^3_{12}(\bar u)\,\Gamma^3_{21}(\bar u)\neq0.\eeq
Then 
\begin{enumerate}
\item[1)]  $0\leq\dim \fs(\R)/\fstriv\leq 4$.
\item[2)]  For each $k=0, \dots, 4$ there exists $\R$, satisfying assumptions of the theorem, such that $\dim \fs(\R)/\fstriv=k$.

\end{enumerate}
\end {theorem}
Condition \eq{eq-mist} arises in the proof of Lemma~\ref{lem-lam}.  Example~\ref{ex-no-mist} illustrates that there are partial frames with non-trivial fluxes, for which   \eq{eq-mist} does not hold. However, from the  proof of Lemma~\ref{lem-lam}, one can see that analyzing the size $\fs(\R)$ in this case becomes rather technical and we left this non-generic case for the future work.

 \begin{lemma}\label{lem-indep}Conditions \eqref{eq-nec-nonhyp} and  \eqref{eq-mist} are independent of the choice of a vector-field $\vs$ that completes $\R$ to a frame. 
 
 \end{lemma}
 \begin{proof} Consider two completions of $\R$ to a local frame  in  a neighborhood $\Omega$ of $\bar u$. The first one is given by a vector field $\vs$, while the second one  is given by a vector field $\vs'$.  Then,  we can express $\vs'$ as linear combination of  $\{\vr_1,\vr_2, \vs\}$:
 $$ \vs'=\alpha\,\vr_1 + \beta\,\vr_2+\gamma\,\vs,$$
 for some smooth functions $\alpha$, $\beta$ and $\gamma$, such that $\gamma$ is nowhere zero on $\Omega$.
  Let $c$'s and $\Gamma$'s be the structure components and Christoffel symbols for connection $\tnabla$ relative to the frame $\vr_1,\vr_2,\vs$ and let $c'$'s and $\Gamma'$'s be the structure components and Christoffel symbols for connection $\tnabla$ relative to the frame $\vr_1,\vr_2,\vs'$. Then  for $i,j=1,2$:
\begin{align*}\tnabla_{\vr_i}\,{\vr_j}&= \Gamma^1_{ij}\,\vr_1+ \Gamma^2_{ij}\,\vr_2+\Gamma^3_{ij}\,\vs= \Gamma^1_{ij}\,\vr_1+ \Gamma^2_{ij}\,\vr_2+ \Gamma^3_{ij}\,\frac 1 \gamma (\vs'-\alpha \vr_1-\beta\,\vr_2)\\
&=\left(\Gamma^1_{ij}-\frac{\alpha}{\gamma}\Gamma^3_{ij}\right)\,\vr_1+\left(\Gamma^2_{ij}-\frac{\beta}{\gamma}\,\Gamma^3_{ij}\right)\,\vr_2+\frac 1 \gamma \Gamma^3_{ij}\, \vs'.
 \end{align*}
 Therefore,  $ \Gamma'^3_{ij}=\frac 1 \gamma \Gamma^3_{ij}$ and so $ c'^3_{ij}=\frac 1 \gamma c^3_{ij}$  for $i,j=1,2$. Then  \eq{eq-mist} and \eq{eq-nec-nonhyp} hold for $c'$'s and $\Gamma'$'s if and only if they hold for  $c$'s and $\Gamma$'s.
 \end{proof}
The fact that   conditions \eq{eq-nec-nonhyp} and \eq{eq-mist} are independent of the completion of $\R$ to a frame suggests that they  can be written as some relations among vector fields $\vr_1, \vr_2$, $[\vr_1, \vr_2]$, $\tnabla_{\vr_1} \vr_2$ and $\tnabla\vr_2, \vr_1$. However, we have not discovered such expressions.
\begin{lemma}\label{lem-a} Let $\R=\{\vr_1,\vr_2\}$ be a non-involutive partial frame satisfying \eqref{eq-nec-sh-m2}. Let $\vs=[\vr_1,\vr_2]$ and for $1\leq i,j,k\leq 3$, let  $c_{ij}^k$ denote the structure functions  and let  $\Gamma_{ij}^k$ denote the Cristoffel symbols for connection $\tnabla$ relative to the local frame $\{\vr_1,\vr_2,\vs\}$. For functions $\lam^1$ and $\lam^2$ defined on an open neighborhood $\Omega$ of $\bar u$ the following two conditions are equivalent:
\begin{enumerate}
\item[1)]  There is a solution $\vf$ of $\fs(\R)$ system \eq{FR2}, for the \emph{prescribed} functions $\lam^1$ and $\lam^2$.

\item[2)] Functions $\lam^1$ and  $\lam^2$, together with functions $\aaa^1$ and $\aaa^2$, defined by 
\begin{align}
\label{eq-a1}\aaa^1 &=-{\vr_2} (\lambda^1)- \Gamma^1_{12}\,(\lambda^1-\lambda^2) \\
\label{eq-a2}\aaa^2&=\hskip3mm {\vr_1} (\lambda^2)- \Gamma^2_{21}\,(\lambda^1-\lambda^2)
\end{align}
satisfy the following system of $6$ equations:

\begin{align}
\label{eq-r1lam1}\vr_1(\lambda^1)&=\frac 1 {\Gamma^3_{21}}\,\Big(\Upsilon_1\,(\lambda^1-\lambda^2)+\Gamma^3_{11}\,\aaa^1+2\,\Gamma^3_{12}\,\aaa^2\Big),\\
\nn\\
\label{eq-r2lam2} \vr_2(\lambda^2)&=\frac 1 {\Gamma^3_{12}}\,\Big(\Upsilon_2\,(\lambda^1-\lambda^2)-2\,\Gamma^3_{21}\, \aaa^1-\Gamma^3_{22}\,\aaa^2\Big),\\
\nn\\ 
\label{eq-r2a1}\vr_2(\aaa^1)&=(\Gamma^1_{23}\,\Gamma^3_{12}-\Gamma^1_{32})\,(\lam^1-\lam^2)+(c^3_{23}-\Gamma^1_{21})\,\aaa^1-\Gamma^1_{22}\,\aaa^2,\\
\nn\\ 
\label{eq-r1a2}{\vr_1} (\aaa^2)&=(\Gamma^2_{13}\,\Gamma^3_{21}+\Gamma^2_{31})\,(\lam^1-\lam^2)-\Gamma^2_{11}\,\aaa^1+(c^3_{13}-\Gamma^2_{12})\,a^2,\\
\nn\\ 
\label{eq-r1a1}\vr_1(\aaa^1)- {\vs} (\lambda^1)&= \Gamma^1_{13}\Gamma^3_{12}\,(\lam^1-\lam^2) -(\Gamma^1_{11}-c^3_{13})\,\aaa^1-\Gamma^1_{12}\,\aaa^2,\\
\nn\\ 
\label{eq-r2a2}{\vr_2} (\aaa^2)-{\vs} (\lambda^2)&=\Gamma^2_{23}\,\Gamma^3_{21}\,(\lam^1-\lam^2)-\Gamma^2_{21}\,\aaa^1+(c^3_{23}-\Gamma^2_{22})\aaa^2,
\end{align}
where
\beq \label{eq-upsilon} \Upsilon_1=\Gamma^3_{12}\,(\Gamma^2_{21}-\Gamma^3_{31})\,-\vr_1(\Gamma_{12}^3) \text { and } \Upsilon_2=\Gamma^3_{21}\,(\Gamma^3_{32}-\Gamma^1_{12} )+\vr_2(\Gamma^3_{12}).\eeq
\end{enumerate}
Moreover, for every $\lam^1$ and $\lam^2$ satisfying condition 2), there exists unique, up to adding a constant vector, flux $\vf$ satisfying \eq{FR2}.
\end{lemma}
\begin{proof}
 We note that due to the symmetry of $\tnabla$ and our definition of $\vs$ we have
\beq\label{eq-c-gamma-12} \Gamma_{12}^1-\Gamma_{21}^1=c_{12}^1=0,\quad  \Gamma_{12}^2-\Gamma_{21}^2=c_{12}^2=0, \quad  \Gamma_{12}^3-\Gamma_{21}^3=c_{12}^3=1.\eeq
\begin{description}
\item[($1\Longrightarrow 2$)]   Assume for $\lam^1$ and $\lam^2$, there exists $\vf$ such that \eqref{FR2} holds.
Then  flatness condition \eqref{R0} implies that
\beq\label{[r1r2]f}\fbox{$\tnabla_{[\vr_1,\vr_2]}\vf=\tnabla_{\vr_1}\tnabla_{\vr_2}\vf-\tnabla_{\vr_2}\tnabla_{\vr_1}\vf.$}\eeq
We recall that $\vs=[\vr_1,\vr_2]$, expend the right-hand side, substitute \eqref{FR2}, and use \eqref{eq-c-gamma-12}, to  derive that
\beq \label{vs(f)}\tnabla_{\vs}\,\vf= \aaa^1\, \vr_1+ \aaa^2 \,\vr_2  +\aaa^3 \,\vs,\eeq
whith  $\aaa^1$ and $\aaa^2$ given by \eqref{eq-a1} and \eqref{eq-a2}, and 
\beq\label{eq-a3} \aaa^3=\Gamma_{12}^3\,\lambda^2-\Gamma_{21}^3\,\lambda^1.\eeq

We record following simple consequences of  \eq{eq-a3} and the last equation in \eq{eq-c-gamma-12} that is repeatedly used below.
\beq\label{a3-lam} \lam^1-\aaa^3=\Gamma^3_{12}\,(\lam^1-\lam^2) \text{ and }  \lam^2-\aaa^3=\Gamma^3_{21}\,(\lam^1-\lam^2)
\eeq

By expanding the flatness identity
\beq \label{r1sf}\fbox{$\tnabla_{[\vr_1,\vs]}\vf=\tnabla_{\vr_1}\tnabla_{\vs}\vf-\tnabla_{\vs}\tnabla_{\vr_1}\vf$}\eeq
we obtain
\begin{align}
\label{[r1s]-r1}\vr_1(\aaa^1)&= {\vs} (\lambda^1) +\Gamma^1_{13}\,\lam^1 -(\Gamma^1_{11}-c^3_{13})\,\aaa^1-\Gamma^1_{12}\,\aaa^2-\Gamma^1_{13}\,\aaa^3&\text{(coefficient of $\vr_1$)}\\
\nn&={\vs} (\lambda^1) +\Gamma^1_{13}\Gamma^3_{12}\,(\lam^1-\lam^2) -(\Gamma^1_{11}-c^3_{13})\,\aaa^1-\Gamma^1_{12}\,\aaa^2,&\\
\nn\\ 
\label{[r1s]-r2}{\vr_1} (\aaa^2)&= \Gamma^2_{31}\,\lam^1+c^2_{13}\,\lam^2-\Gamma^2_{11}\,\aaa^1+(c^3_{13}-\Gamma^2_{12})\,a^2-\Gamma^2_{13}\,\aaa^3&\text{(coefficient of $\vr_2$)}\\
\nn &=(\Gamma^2_{13}\,\Gamma^3_{21}+\Gamma^2_{31})\,(\lam^1-\lam^2)-\Gamma^2_{11}\,\aaa^1+(c^3_{13}-\Gamma^2_{12})\,a^2,&\\
\nn\\ 
\label{[r1s]-s}{\vr_1} (\aaa^3)&=\Gamma^3_{31}\,\lam^1-\Gamma^3_{11}\,\aaa^1-\Gamma^3_{12}\,\aaa^2-\Gamma^3_{31}\,\aaa^3&\text{(coefficient of $\vs$)}\\
\nn &=\Gamma^3_{31}\,\Gamma^3_{12}(\lam^1-\lam^2)-\Gamma^3_{11}\,\aaa^1-\Gamma^3_{12}\,\aaa^2,
\end{align}
where $a^3$ was eliminated from the right-hand sides of the above equations using \eqref{a3-lam}.

Similarly,   identity 
\beq\label{r2sf}\fbox{$\tnabla_{[\vr_2,\vs]}\vf=\tnabla_{\vr_2}\tnabla_{\vs}\vf-\tnabla_{\vs}\tnabla_{\vr_2}\vf$}\eeq
leads to
\begin{align}
\label{[r2s]-r1}\vr_2(\aaa^1)&=c^1_{23}\,\lam^1+\Gamma_{32}^1\,\lam^2+(c^3_{23}-\Gamma^1_{21})\,\aaa^1-\Gamma^1_{22}\,\aaa^2-\Gamma^1_{23}\,\aaa^3&\text{(coefficient of $\vr_1$)}\\
\nn&= (\Gamma^1_{23}\,\Gamma^3_{12}-\Gamma^1_{32})\,(\lam^1-\lam^2)+(c^3_{23}-\Gamma^1_{21})\,\aaa^1-\Gamma^1_{22}\,\aaa^2,\\
\nn\\ 
\label{[r2s]-r2}{\vr_2} (\aaa^2)&= {\vs} (\lambda^2)+\Gamma^2_{23}\,\lam^2-\Gamma^2_{21}\,\aaa^1+(c^3_{23}-\Gamma^2_{22})\aaa^2-\Gamma^2_{23}\,\aaa^3&\text{(coefficient of $\vr_2$)}\\
\nn&={\vs} (\lambda^2)+\Gamma^2_{23}\,\Gamma^3_{21}\,(\lam^1-\lam^2)-\Gamma^2_{21}\,\aaa^1+(c^3_{23}-\Gamma^2_{22})\aaa^2,\\
\nn\\ 
\label{[r2s]-s}{\vr_2} (\aaa^3)&=\Gamma^3_{32}\,\lam^2-\Gamma^3_{21}\,\aaa^1-\Gamma^3_{22}\,\aaa^2-\Gamma^3_{32}\,\aaa^3&\text{(coefficient of $\vs$)}\\
\nn&=\Gamma^3_{32}\,\Gamma^3_{21}\,(\lam^1-\lam^2)-\Gamma^3_{21}\,\aaa^1-\Gamma^3_{22}\,\aaa^2.
\end{align}
We note that  \eq{[r2s]-r1},  \eq{[r1s]-r2}, \eq{[r1s]-r1}, \eq{[r2s]-r2} coincide with \eq{eq-r2a1}, \eq{eq-r1a2}, \eq{eq-r1a1}, and \eq{eq-r2a2}, respectively. 
To show the remaining two equations,  \eq{eq-r1lam1} and \eq{eq-r2lam2}, we note that equations \eq{[r1s]-s}  and \eq{[r2s]-s} express the derivatives of $a_3$ in the $\vr_1$ and $\vr_2$ directions, respectively. However, these derivatives can be also obtained by differentiating  \eq{eq-a3} and substituting  \eq{eq-a1} and  \eq{eq-a2}:
\begin{align}
\label{r1(a3)}\vr_1(\aaa_3)&=\Gamma_{12}^3\,\vr_1(\lambda^2)-\Gamma_{21}^3\,\vr_1(\lambda^1)+\vr_1(\Gamma_{12}^3)\,(\lam^2-\lam^1)\\
\nn&=\left(\Gamma_{12}^3\,\Gamma^2_{21}-\vr_1(\Gamma_{12}^3)\right)\,(\lambda^1-\lambda^2)+ \Gamma_{12}^3 \aaa^2-\Gamma_{21}^3\,\vr_1(\lambda^1),\\
\nn\\
\label{r2(a3)}\vr_2(\aaa_3)&=\Gamma_{12}^3\,\vr_2(\lambda^2)-\Gamma_{21}^3\,\vr_2(\lambda^1)+\vr_2(\Gamma_{12}^3)\,(\lam^2-\lam^1)\\
\nn&=\Gamma_{12}^3\,\vr_2(\lambda^2)+\left(\vr_2(\Gamma_{12}^3)-\Gamma_{21}^3\,\Gamma^1_{12}\right)\,(\lambda^2-\lambda^1) +\Gamma_{21}^3\, \aaa^1,
\end{align}
where we used that, due to the last equation in \eq{eq-c-gamma-12}, derivatives of $\Gamma_{12}^3$ and  $\Gamma_{21}^3$ are equal.
From \eq{[r1s]-s} and \eq{r1(a3)} we obtain:
\beq\label{eq-r1(lam1)}\Gamma_{21}^3\,\vr_1(\lambda^1)=\left(\Gamma^3_{12}\,(\Gamma^2_{21}-\Gamma^3_{31})\,-\vr_1(\Gamma_{12}^3)\right)\,(\lambda^1-\lambda^2)+\Gamma^3_{11}\,\aaa^1+2\,\Gamma^3_{12}\,\aaa^2\eeq
Similarly, from \eq{[r2s]-s} and \eq{r2(a3)} we obtain: 
\beq\label{eq-r2(lam2)}\Gamma_{12}^3\,\vr_2(\lambda^2)=\left(\vr_2(\Gamma_{12}^3)+\Gamma^3_{21}\,(\Gamma^3_{32}-\Gamma^1_{12}\right)\,(\lambda^1-\lambda^2)-2\,\Gamma^3_{21}\, \aaa^1-\Gamma^3_{22}\,\aaa^2.\eeq
Condition \eqref{eq-nec-shg} imply that  $\Gamma^3_{21}\neq 0$ and $\Gamma^3_{12}\neq 0$, and, therefore, 
we can solve   \eq{eq-r1(lam1)} and \eq{eq-r2(lam2)} for $\vr_1(\lam^1)$ and  $\vr_2(\lam^2)$, establishing  \eq{eq-r1lam1} and \eq{eq-r2lam2}.
\item[($2\Longrightarrow 1$ and uniqueness)] Given functions $\lam^1$ and $\lam^2$, let $\aaa^1$, $\aaa^2$ and $\aaa^3$ be defined by \eqref{eq-a1},  \eqref{eq-a2}, \eqref{eq-a3} respectively. Then  equations \eq{FR2},  \eq{vs(f)} constitute a Frobenius-type system on the three unknown functions -- the components of  the flux $\vf$. It is straightforward to check that the integrability conditions for this system coincide  with of the  flatness conditions  \eq{[r1r2]f}, \eq{r1sf} and \eq{r2sf}. Reversing  the proof of part 1), we see that they are satisfied provided  $\lam^1$ and $\lam^2$,   satisfy condition 2). Thus if $\lam^1$ and $\lam^2$ satisfy condition 2), then for any prescription of the initial value  $\vf(\bar u)$, there exists a unique $\vf$ satisfying  \eq{FR2} and \eq{vs(f)}. Moreover, since \eq{vs(f)} is a consequence  of \eq{FR2}, there is a unique $\vf$ satisfying  \eq{FR2}  for any prescription of the initial value  $\vf(\bar u)$. We, therefore, conclude that  the generic solution  \eq{FR2} depends on three arbitrary constants. We finally note that  if $\vf$ satisfies \eq{FR2}, then so does $\vf+\text{(a constant vector in $\RR^3$)}$, and, therefore, the three arbitrary constants in the generic solution correspond to the components of an arbitrary constant vector. Thus, for the given pair of functions $\lam^1$ and $\lam^2$, the solution  of the $\fs(\R)$-system \eq{FR2} is unique up to addition of a constant vector.
 \end{description}
\end{proof}
\begin{lemma} \label{lem-lam}Let $\R=\{\vr_1,\vr_2\}$ be a partial frame satisfying assumptions of the Theorem~\ref{thm-noninv-size}. 
Then the set of pairs of functions $\lam(\R)=\{(\lam^1,\lam^2)\}$ satisfying condition 2)~of Lemma~\ref{lem-a} is a real vector space of dimension at most 5.
\end{lemma}
\begin{proof} It is straightforward to check that $\lam(\R)$ is a vector space. To prove the bound on its dimension, we prolong the system of equations \eq{eq-a1} -- \eq{eq-r2a2}, listed in condition 2) of Lemma~\ref{lem-a}, to a system of the  Frobenius-type on 5 unknown functions $\lam^1$, $\lam^2$, $\aaa^1$, $\aaa^2$, and $\tau$, where we define 
\beq
\label{tau} \tau=\vs(\lam^2) \text { for } \vs=[\vr_1,\vr_2].
\eeq
This is done in the following steps.
\begin{enumerate}
\item[1)] By expanding the right-hand side of the commutator relationship 
$$\fbox{$\vs(\lam^1)=[\vr_1,\vr_2](\lam^1)$}$$
and substitution of the expressions for $\vr_1(\lam^1)$,   $\vr_2(\lam^1)$, $\vr_1(\lam^2)$,   $\vr_2(\lam^2)$, $\vr_1(\aaa^1)$,   $\vr_2(\aaa^1)$, $\vr_1(\aaa^2)$,   $\vr_2(\aaa^2)$ from \eq{eq-a1} -- \eq{eq-r2a2}, we obtain
 \beq\label{eq-s(lam)-1} 2 \,\Gamma^3_{21}\,\vs(\lam^1)+ 2\,\Gamma^3_{12} \,\vs(\lam^2)= \Gamma^3_{21}\,\left(A_1\,(\lam^1-\lam^2)+B_1\,\aaa^1+C_1\, \aaa^2\right),\eeq
where
 \begin{align}
  \label{eq-A1} A_1&=-\vr_2\,\left(\frac {\Upsilon_1} {\Gamma^3_{21}}\right)-\vr_1(\Gamma^1_{12})+\frac {\Upsilon_1\,\Upsilon_2} {\Gamma^3_{21}\,\Gamma^3_{12}}-\Gamma^3_{12}\,(\Gamma^1_{13}+2\,\Gamma^2_{23})-\frac{\Gamma^3_{11}}{\Gamma^3_{21}}\,(\Gamma^1_{23}\,\Gamma^3_{12}-\Gamma^1_{32})+\Gamma^1_{12}\,\Gamma^2_{21}\\
  B_1&= -\vr_2\left(\frac{\Gamma^3_{11}}{\Gamma^3_{21}}\right)-\frac{\Upsilon_1\,(\Gamma^3_{21}-1)} {\Gamma^3_{21}\,\Gamma^3_{12}}-\frac{\Gamma^3_{11}\,c^3_{23}}{\Gamma^3_{21}}+2\,\frac{\Gamma^3_{12}\,\Gamma^2_{21}}{\Gamma^3_{21}}+\Gamma^1_{11}-c^3_{13}
 \\
 C_1 &=2\,\frac{\vr_2(\Gamma^3_{12})}{(\Gamma^3_{21})^2}-\frac{\Upsilon_1\,\Gamma^3_{22}} {\Gamma^3_{21}\,\Gamma^3_{12}}- 2\,\frac{\Gamma^1_{12}}{\Gamma^3_{21}}+\frac{\Gamma^3_{11}\,\Gamma^1_{22}}{\Gamma^3_{21}}\,+2\,\frac{\Gamma^3_{12}}{\Gamma^3_{21}}\,(\Gamma^2_{22}-c^3_{23})
 \end{align}
  \item[2)] By expanding the right-hand side of the commutator relationship 
$$\fbox{$\vs(\lam^2)=[\vr_1,\vr_2](\lam^2)$}$$
and substitution of the expressions for $\vr_1(\lam^1)$,   $\vr_2(\lam^1)$, $\vr_1(\lam^2)$,   $\vr_2(\lam^2)$, $\vr_1(\aaa^1)$,   $\vr_2(\aaa^1)$, $\vr_1(\aaa^2)$,   $\vr_2(\aaa^2)$ from \eq{eq-a1} -- \eq{eq-r2a2}, we obtain
\beq\label{eq-s(lam)-2} 2 \,\Gamma^3_{21}\,\vs(\lam^1)+ 2\,\Gamma^3_{12} \,\vs(\lam^2)= \Gamma^3_{12} \left(A_2(\lam^1-\lam^2)+B_2\,\aaa^1+C_2\,\aaa^2\right).
  \eeq
  \begin{align}
\label{eq-A2}  A_2&= \vr_1\left(\frac {\Upsilon_2} {\Gamma^3_{12}}\right)-\vr_2\,\left(\Gamma^2_{21}\right)+\frac {\Upsilon_2\,\Upsilon_1} {\Gamma^3_{12}\,\Gamma^3_{21}}\,- {\Gamma^3_{21}}(\Gamma^2_{23}+2\, \Gamma^1_{13})-\frac{\Gamma^3_{22}}{\Gamma^3_{12}}\,(\Gamma^2_{13}\,\Gamma^3_{21}+\Gamma^2_{31})+\Gamma^2_{21}\,\Gamma^1_{12}\\
  B_2&=-2\,\frac{\vr_1\left(\Gamma^3_{21}\right)}{(\Gamma^3_{12})^2}+\frac {\Upsilon_2\,\Gamma^3_{11}} {\Gamma^3_{12}\,\Gamma^3_{21}}+2\,\frac{\Gamma^2_{21}}{\Gamma^3_{12}}+\frac{ \Gamma^3_{22}\,\Gamma^2_{11}}{\Gamma^3_{12}}+2\,\frac{\Gamma^3_{21}}{\Gamma^3_{12}}\,(\Gamma^1_{11}-c^3_{13})\\
  %
  %
  C_ 2&=-\vr_1\left(\frac{\Gamma^3_{22}}{\Gamma^3_{12}}\right)+\frac {\Upsilon_2\,(\Gamma^3_{12}+1)} {\Gamma^3_{12}\,\Gamma^3_{21}}-\frac{\Gamma^3_{22}\,c^3_{13}}{\Gamma^3_{12}}+\frac{2\,\Gamma^3_{21}\,\Gamma^1_{12}}{\Gamma^3_{12}}+\Gamma^2_{22}-c^3_{23}
  \end{align}
\item[3)]  Observe that \eq{eq-s(lam)-1} and \eq{eq-s(lam)-2} have identical left-hand sides, and so their right-hand side must be equal.
  It turns out that this indeed the case. In fact,
  \beq\label{ABCid}\Gamma^3_{21}A_1\equiv \Gamma^3_{12} A_2,\quad \Gamma^3_{21}B_1\equiv \Gamma^3_{12} B_2 \text{ and } \Gamma^3_{21} C_1\equiv \Gamma^3_{12} C_2\eeq
  due to flatness condition \eq{R0}.
To show  the A-identity in \eq{ABCid}, we first compute  $\Gamma^3_{21}A_1- \Gamma^3_{12} A_2$ by substituting $\Upsilon_1$ and $\Upsilon_2$  into \eq{eq-A1} and \eq{eq-A2}  and making various simplifications. We obtain
\begin{align}\nn
&\Gamma^3_{21}A_1- \Gamma^3_{12} A_2=
 -\vs(\Gamma^3_{12})+\Gamma^3_{12}\,\vr_2(\Gamma^3_{31})-\Gamma^3_{21}\,\vr_1(\Gamma^3_{32})-\Gamma^2_{21}\,\Gamma^3_{32}+\Gamma^3_{31}\,\Gamma^3_{32} -\Gamma^3_{31}\Gamma^1_{12} \\
\label{Acheck} &+ \Gamma^3_{21}\,\Gamma^3_{12}\,(\Gamma^1_{13}-\Gamma^2_{23})-{\Gamma^3_{11}}\,(\Gamma^1_{23}\Gamma^3_{12}-\Gamma^1_{32})+{\Gamma^3_{22}}\,(\Gamma^2_{13}\,\Gamma^3_{21}+\Gamma^2_{31}).
\end{align}
We then expand the identity
\beq\label{commA}\Gamma^3_{12}\,\Big(\nabla_{\vr_2}\,\nabla_{\vs}\vr_1-\nabla_{\vs}\,\nabla_{\vr_2}\vr_1-\nabla_{[\vr_2,\vs]}\,\vr_1\Big) -\Gamma^3_{21}\,\Big(\nabla_{\vr_1}\,\nabla_{\vs}\vr_2-\nabla_{\vs}\,\nabla_{\vr_1}\vr_2-\nabla_{[\vr_1,\vs]}\,\vr_2\Big)\equiv 0.\eeq
and observe that the coefficient of $\vs$ in \eq{commA} equals to the left hand side of \eq{Acheck}. Similarly, we use the $\vs$ coefficient
of the expanded identity  $\nabla_{\vr_1}\nabla_{\vr_2}\vr_1-\nabla_{\vr_2}\nabla_{\vr_1}\vr_1\equiv \nabla_{\vs}\vr_2$ to show the $B$-identity of
\eq{ABCid},  and the $\vs$ coefficient
of the expanded identity  $\nabla_{\vr_1}\nabla_{\vr_2}\vr_2-\nabla_{\vr_2}\nabla_{\vr_1}\vr_2\equiv\nabla_{\vs}\vr_2$ to show the $C$-identity of
\eq{ABCid}. 

  \item[4)] Introducing a new unknown function $\tau$,  defined by \eq{tau},
we solve   \eq{eq-s(lam)-1}  for $\vs(\lam^1)$:
\beq\label{s(lam1)} \vs(\lam^1)=-\,\frac{\Gamma^3_{12}}{\Gamma^3_{21}} \, \tau+ \frac 1 2 A_1\,(\lam^1-\lam^2)+\frac 12 \,B_1\,\aaa^1+\frac 12 C_1\, \aaa^2\eeq
and rewrite  \eq{eq-r1a1} and \eq{eq-r2a2} as
\begin{align}
\label{eq-r1(a1)}\vr_1(\aaa^1)&=\left( \frac 1 2 A_1 +\Gamma^1_{13}\Gamma^3_{12}\right)\,(\lam^1-\lam^2)+\left(\frac 12 \,B_1-\Gamma^1_{11}+c^3_{13}\right)\,\aaa^1+\left(\frac 12 C_1 -\Gamma^1_{12}\right)\, \aaa^2-\,\frac{\Gamma^3_{12}}{\Gamma^3_{21}} \, \tau,\\
\nn\\ 
\label{eq-r2(a2)}{\vr_2} (\aaa^2)&=\Gamma^2_{23}\,\Gamma^3_{21}\,(\lam^1-\lam^2)-\Gamma^2_{21}\,\aaa^1+(c^3_{23}-\Gamma^2_{22})\aaa^2+\tau,
\end{align}

 \item[5)] To complete the system \eq{eq-a1} - \eq{eq-r1a2}, \eq{tau}, \eq{s(lam1)}, \eq{eq-r1(a1)} and \eq{eq-r2(a2)}   to a system of the Frobenius type we need to express the remaining derivatives $\vs(a^1)$, $\vs(a^2)$,
$\vr_1(\tau)$, $\vr_2(\tau)$ and $\vs(\tau)$ as   functions of $\lam^1, \lam^2, \aaa^1,\aaa^2$ and $\tau$. For this purpose, we continue to consider various consequences  of commutator relationships.

From 
$$\fbox{$[\vr_1,\vs](\lam^2)=\vr_1(\vs(\lam^2))-\vs(\vr_1(\lam^2))$}$$
expending the left-hand side and  substituting \eq{tau} and \eq{eq-a2} into the right-hand side, we get:
\beq\label{[r1,s](l2)}  c^1_{13}\,\vr_1(\lam^2)+c^2_{13}\,\vr_2(\lam^2)+c^3_{13}\,\vs(\lam^2)=\vr_1(\tau)-\vs(\Gamma^2_{21}\,(\lambda^1-\lambda^2)+  \aaa^2).\eeq
By substituting the already known expressions of the directional derivatives, $\vr_1(\lam^2)$, $\vr_2(\lam^2)$, $\vs(\lam^2)$  and $\vs(\lam^1)$, given by \eq{eq-a2}, \eq{eq-r2lam2}, \eq{tau} and \eq{s(lam1)} into \eq{[r1,s](l2)}, respectively, we obtain:
\beq\label{eq-L1} \vr_1(\tau)-\vs (\aaa^2) =\cl_1(\lam^1-\lam^2, \aaa^1,\aaa^2, \tau),\eeq
where  $\cl_1$ is some known, linear in its arguments function with coefficients depending on $c$'s, $\Gamma$'s  and their derivatives.  The explicit expression of $\cl_1$ is  too long to be included.

From $$\fbox{$[\vr_2,\vs](\lam^2)=\vr_2(\vs(\lam^2))-\vs(\vr_2(\lam^2))$}$$
expending the left-hand side and  substituting \eq{tau} and \eq{eq-r2lam2} into the right-hand side, we get:
\beq\label{[r2,s](l2)}  c^1_{23}\,\vr_1(\lam^2)+c^2_{23}\,\vr_2(\lam^2)+c^3_{23}\,\vs(\lam^2)=\vr_2(\tau)-\vs\left(\frac 1 {\Gamma^3_{12}}\,\Big(\Upsilon_2\,(\lambda^1-\lambda^2)-2\,\Gamma^3_{21}\, \aaa^1-\Gamma^3_{22}\,\aaa^2\Big)\right). \eeq
By substituting the already known expressions of the directional derivatives, $\vr_1(\lam^2)$, $\vr_2(\lam^2)$, $\vs(\lam^2)$  and $\vs(\lam^1)$  into \eq{[r2,s](l2)}
we obtain \beq\label{eq-L2} {\Gamma^3_{12}} \,\vr_2(\tau)+2\, \Gamma^3_{21}\,\vs (\aaa^1)+\Gamma^3_{22}\,\vs (\aaa^2) =\cl_2(\lam^1-\lam^2, \aaa^1,\aaa^2, \tau),\eeq
where function $\cl_2$ is linear in its arguments with coefficients depending on $c$'s, $\Gamma$'s  and their derivatives.

Similarly from the commutator relationships 
$$\fbox{$[\vr_2,\vs](\lam^1)=\vr_2(\vs(\lam^1))-\vs(\vr_2(\lam^1))$}
\text{ and }
\fbox{$[\vr_1,\vs](\lam^1)=\vr_1(\vs(\lam^1))-\vs(\vr_1(\lam^1))$}$$
we obtain equations
\beq\label{eq-L3}- {\Gamma^3_{12}} \,\vr_2(\tau)+\Gamma^3_{21}\,\vs (\aaa^1) =\cl_3(\lam^1-\lam^2, \aaa^1,\aaa^2, \tau),\eeq
and 
\beq\label{eq-L4}\Gamma^3_{12} \,\vr_1(\tau)+\Gamma^3_{11}\,\vs (\aaa^1)+2\,\Gamma^3_{12}\,\vs (\aaa^2) =\cl_4(\lam^1-\lam^2, \aaa^1,\aaa^2, \tau),\eeq
where function $\cl_3$ and $\cl_3$  are linear in its arguments with coefficients depending on $c$'s, $\Gamma$'s  and their derivatives.

Equations \eq{eq-L1}, \eq{eq-L2},  \eq{eq-L3} and  \eq{eq-L4} can be viewed as a linear inhomegeneous system of four equations on the four unknowns
$\vs(\aaa^1)$, $\vs(\aaa^2)$, $\vr_1(\tau)$ and $\vr_2(\tau)$:
\beq\label{4matrix}
\left[
\begin{array}{cccc}
 1 & 0  &0   &-1\\
 0 &   \Gamma^3_{12}&   2\, \Gamma^3_{21} & \Gamma^3_{22}\\
  0&- \Gamma^3_{12}   &  \Gamma^3_{21}   &0\\
   \Gamma^3_{12}  & 0  &\Gamma^3_{11}     & 2\,\Gamma^3_{12} 
\end{array}
\right]\,
\left[\begin{array}{c}
\vr_1(\tau)  \\
\vr_2(\tau)\\
\vs(\aaa^1)\\
\vs(\aaa^2)
\end{array}
\right]\,
=\left[\begin{array}{c}
\cl_1\\
 \cl_2\\
 \cl_3\\
  \cl_4
\end{array}
\right].
\eeq
We find that 
$$\det(M)=\Gamma^3_{12}\,(9 \Gamma^3_{12} \,\Gamma^3_{21} -\Gamma^3_{11}\,\Gamma^3_{22}).$$   
We note that under the assumptions of Theorem~\ref{thm-noninv-size} $\det(M)\neq 0$, in a neighborhood of $\bar u$,
and, hence,  \eq{4matrix} can be solved to find  expressions $\vs(\aaa^1)$, $\vs(\aaa^2)$, $\vr_1(\tau)$ and $\vr_2(\tau)$ as linear functions of $\lam^1-\lam^2, \aaa^1,\aaa^2, \tau$, with coefficients depending on $c$'s, $\Gamma$'s  and their derivatives. Finally, 
\beq\label{s(tau)}\vs(\tau)=[\vr_1,\vr_2](\tau)=\vr_1(\vr_2(\tau))-\vr_2(\vr_1(\tau)),\eeq
after substitution of  already known expressions of the derivatives, $\vr_1(\tau)$, $\vr_2(\tau)$, $\vr_1(\lam^1)$, $\vr_2(\lam^1)$, $\vr_1(\lam^2)$, $\vr_2(\lam^2)$, $\vr_1(\aaa^1)$, $\vr_2(\aaa^1)$, $\vr_1(\aaa^2)$, and  $\vr_2(\aaa^2)$, also becomes a linear  function of $\lam^1-\lam^2, \aaa^1,\aaa^2, \tau$, with coefficients depending on $c$'s, $\Gamma$'s  and their derivatives. 
\item[6)] The fifteen equations \eq{eq-a1} - \eq{eq-r1a2}, \eq{tau}, \eq{s(lam1)}, \eq{eq-r1(a1)}, \eq{eq-r2(a2)}, \eq{4matrix} and \eq{s(tau)} can be used express all the directional derivatives of functions  $\lam^1$, $\lam^2$, $\aaa^1$, $\aaa^2$ and $\tau$ as linear combinations of   $\lam^1- \lam^2$, $\aaa^1$, $\aaa^2$ and $\tau$  with coefficients depending on $c$'s, $\Gamma$'s  and their derivatives. Therefore, we obtain a Frobenius-type system.   If the integrability conditions for this system are identically satisfied, the  generic solution depends on  5 constants -- the prescribed values of these functions at $\bar u$.   If the integrability conditions for this system are not identically satisfied, they will impose additional relationships on $\lam^1$, $\lam^2$, $\aaa^1$, $\aaa^2$ and $\tau$, thus reducing the size of the solution set.

\item[7)] The  Frobenius-type system  \eq{eq-a1} - \eq{eq-r1a2}, \eq{tau}, \eq{s(lam1)}, \eq{eq-r1(a1)}, \eq{eq-r2(a2)}, \eq{4matrix} and \eq{s(tau)}  was obtained as a consequence of condition 2)~of Lemma~\ref{lem-a}. Therefore,  the vector space of pairs
  functions $\lam(\R)=\{(\lam^1,\lam^2)\}$ satisfying this condition is at most 5.
\end{enumerate}
\end{proof}

\begin{lemma} \label{lem-non-hyp1}Let $\R=\{\vr_1,\vr_2\}$ be a non-involutive partial frame  satisfying conditions \eqref{eq-nec-sh-m2} at $\bar u\in \Omega$. Assume $\vf\in \fs(\R)$ is a non-hyperbolic flux. Then  the corresponding eigenfunctions $\lam^1$ and $\lam^2$, appearing in \eqref{FR2}, coincide and are non constant, i.e:
$$\lam^1=\lam^2=\lam, \text{ where } \lam \text{ is a non-constant function}. $$
  \end{lemma}
\begin{proof}
We recall that $\vf$  being non-hyperbolic means that the operator $\tnabla_{(\cdot)}\vf$, does not posses three real eigenfunctions. However, by \eqref{FR2}, it possesses two real  eigenfunctions $\lam^1$ and $\lam^2$. Since complex eigenfunctions come in conjugate pairs, in the $n=3$ case, a possibility of the third eigenfunction being complex is excluded, and, therefore, $\vf$ must posses a generalized eigenvector field, which we denote $\vs$.
Let $c$'s and $\Gamma$'s denote structure coefficients and Christoffel symbols for $\tnabla$, relative to the frame $\{\vr_1, \vr_2, \vs\}$.
\begin{enumerate}
\item ($\lam^1=\lam^2$) To prove by contradiction, we assume that $\lam^1\neq \lam^2$. Then either  
\beq\label{eq-genev} \tnabla_{\vs}\vf =\vr_1+\lam^1\,\vs\,\quad \text{ or }  \quad \tnabla_{\vs}\vf =\vr_2+\lam^2\,\vs.\eeq
Without loss of generality, we assume that the second equality holds (otherwise relabel $\vr_1$ and $\vr_2$).  
Then  \eq{[r1r2]f}
 together with \eq{FR2} and the second equation in \eqref{eq-genev} imply:
 \beq\label{eq-qw} 
 c^1_{12} \lam^1\vr_1+ c^2_{12} \lam^2\vr_2 + c^3_{12} (\vr_2+\lam^2\,\vs)=\tnabla_{\vr_1}(\lam^2 \,\vr_2)-\tnabla_{\vr_2}(\lam^1\,\vr_1). \eeq 
 Using \eq{FR2} again and collecting the coefficients with $\vs$, we obtain
 $$ c^3_{12}\, \lam^2=\lam^2\,\Gamma^3_{12}-\lam^1\,\Gamma^3_{21} \text{ or, equivalently, } \Gamma^3_{21}(\lam^2-\lam^1)=0. $$
 Conditions in \eqref{eq-nec-sh-m2} imply that $\Gamma^3_{21}\neq 0$ on an open neighborhood of $\bar u$, and, therefore,  $\lam^1=\lam^2$ on this neighborhood.
 \item($\lam\neq const$) Let $\lam^1=\lam^2=\lam$. Then, since $\vs$ is a generalized eigenvector field, we must have
  \beq\label{eq-genev2} \tnabla_{\vs}\vf =\alpha\vr_1+\beta \vr_2+\lam\,\vs,\eeq
 where $\alpha$ and $\beta$ are some functions, such that $\alpha(\bar u)$ or $\beta(\bar u)$ is non-zero.
 To prove by contradiction, we assume that $\lam$ is a constant function in a neighborhood of $\bar u$.
 Then \eq{[r1r2]f}, together with   \eq{FR2}, \eq{eq-genev2}  imply that
 \beq\label{eq-qw1} 
 c^1_{12} \lam\vr_1+ c^2_{12} \lam\vr_2 + c^3_{12} (\alpha\,\vr_1+\beta\,\vr_2+\lam\,\vs)=\lam\tnabla_{\vr_1}\,\vr_2-\lam\tnabla_{\vr_2}\,\vr_1. \eeq 
On the left hand-side of  \eq{eq-qw1}, we notice that  $c^1_{12} \lam\vr_1+ c^2_{12} \lam\vr_2 + c^3_{12} \lam\,\vs=\lam\, [\vr_1,\vr_2]$. At the same time, the right-hand side of  \eq{eq-qw1} equals to $\lam\, [\vr_1,\vr_2]$ due to the symmetry condition \eq{T0}.
 Then $\alpha\,\vr_1+\beta\,\vr_2=0$, which contradicts our assumption that vectors $\vr_1|_{\bar u}$ and $\vr_2|_{\bar u}$ are 
 independent and  $\alpha$ and $\beta$ are some functions such that   $\alpha(\bar u)$ or $\beta(\bar u)$ is non zero. Thus $\lam$ is a non-constant function.
 
 \end{enumerate}
 
\end{proof}
\begin{lemma} \label{lem-non-hyp2}Let $\R=\{\vr_1,\vr_2\}$ be a non-involutive partial frame  satisfying conditions \eqref{eq-nec-sh-m2}. Assume $\vf\in \fs(\R)$ is a non-hyperbolic flux. Then all other non-hyperbolic fluxes in   $\fs(\R)$  are of the form $c\,\vf+(\text{a trivial flux})$ where $c\neq 0\in\RR$.
\end{lemma}
\begin{proof} 
\begin{enumerate}
\item[1)] Let $\vf\in \fs(\R)$ be a non-hyperbolic flux.  From Lemma~\ref{lem-non-hyp1}, it follows that there exists
a non-constant function $\lam$ in a neighborhood of $\bar u$, such that $\vf$ and $\lam^1=\lam^2=\lam$ satisfy 
\eq{FR2}.
It is straightforward to check that  $c\,\vf+\bar\lam\bar\vf $, where $\bar\lam \in \RR$ and $\bar f\in \fsid$ (see \eq{def-idf} to recall the definition of identity fluxes) is a non-hyperbolic flux, which  together with
$\lam^1=\lam^2=c\,\lam+\bar\lam$ satisfy 
\eq{FR2}. Recalling \eq{def-triv}, we conclude that $c\,\vf+(\text{a trivial flux})$ belongs to $\fs(\R)$ and clearly those fluxes are non-hyperbolic. It remains to show that  any non-hyperbolic flux in $\fs(\R)$ is of this form.

\item[2)]  Lemma~\ref{lem-a} implies that
function $\lam$ together with functions $\aaa^1$ and $\aaa^2$, defined by 
\begin{align}
\label{eq-a1-nh}\aaa^1 &=-{\vr_2} (\lambda)\\
\label{eq-a2-nh}\aaa^2&=\hskip3mm {\vr_1} (\lambda)
\end{align}
satisfy the following equations (these are  \eq{eq-r1lam1}--\eq{eq-r2a2} in the case when $\lam^1=\lam^2$):
\begin{align}
\label{eq-r1lam}\vr_1(\lambda)&=\frac 1 {\Gamma^3_{21}}\Big(\Gamma^3_{11}\,\aaa^1+2\,\Gamma^3_{12}\,\aaa^2\Big),\\
\nn\\
\label{eq-r2lam} \vr_2(\lambda)&=-\frac 1 {\Gamma^3_{12}}\,\Big(2\,\Gamma^3_{21}\, \aaa^1+\Gamma^3_{22}\,\aaa^2\Big),\\
\nn\\ 
\label{eq-r2a1-nh}\vr_2(\aaa^1)&=(c^3_{23}-\Gamma^1_{21})\,\aaa^1-\Gamma^1_{22}\,\aaa^2,\\
\nn\\ 
\label{eq-r1a2-nh}{\vr_1} (\aaa^2)&=-\Gamma^2_{11}\,\aaa^1+(c^3_{13}-\Gamma^2_{12})\,a^2,\\
\nn\\ 
\label{eq-r1a1-nh}\vr_1(\aaa^1)- {\vs} (\lam)&= -(\Gamma^1_{11}-c^3_{13})\,\aaa^1-\Gamma^1_{12}\,\aaa^2,\\
\nn\\ 
\label{eq-r2a2-nh}{\vr_2} (\aaa^2)-{\vs} (\lam)&=-\Gamma^2_{21}\,\aaa^1+(c^3_{23}-\Gamma^2_{22})\aaa^2,
\end{align}
where,  $\Gamma$'s are Christoffel symbols for $\tnabla$ relative to the frame $\{\vr_1, \vr_2, \vs=[\vr_1,\vr_2]\}$.
Equations \eq{eq-a1-nh} and \eq{eq-a2-nh} immediately imply that:
\beq\label{eq-sl-nh}\vs(\lam)=\vr_1(\vr_2(\lam))-\vr_2(\vr_1(\lam))=-\vr_1(\aaa^1)-\vr_2(\aaa^2).\eeq
Then from \eq{eq-sl-nh}, together with  \eq{eq-r1a1-nh} and  \eq{eq-r2a2-nh}, we obtain:
\begin{align}
\label{eq-sl-nh1}\vs(\lam)&= \frac 1 3\,(\Gamma^1_{11}+\Gamma^2_{21}-c^3_{13})\,\aaa^1+\frac 1 3\,(\Gamma^2_{22}+\Gamma^1_{12}-c^3_{23})\,\aaa^2\,,\\
\label{eq-r1a1-nh1}\vr_1(\aaa^1)&= \frac 1 3\,(-2\,\Gamma^1_{11}+\Gamma^2_{21}+2\,c^3_{13})\,\aaa^1+\frac 1 3\,(\Gamma^2_{22}-2\,\Gamma^1_{12}-c^3_{23})\,\aaa^2,\\
\label{eq-r2a2-nh1}{\vr_2} (\aaa^2)&=\frac 1 3\,(\Gamma^1_{11}-2\,\Gamma^2_{21}-c^3_{13})\,\aaa^1+\frac 1 3\,(-2\,\Gamma^2_{22}+\Gamma^1_{12}+2\,c^3_{23})\,\aaa^2,
\end{align}
From Lemma~\ref{lem-non-hyp1}, we know that $\lam$ is a non-constant function, and, therefore, at least one of its  derivatives in the frame directions must be non-zero. Examining \eq{eq-a1-nh},  \eq{eq-a2-nh} and \eq{eq-sl-nh1}, we conclude that at least one of the functions $\aaa^1$ or $\aaa^2$ is non zero.  Without loss of generality, we assume that $\aaa^1\neq 0$ (otherwise, relabel $\vr_1$ and $\vr_2$).
\item[3)] 
Equations \eq{eq-a1-nh},  \eq{eq-a2-nh}, \eq{eq-r1lam}, \eq{eq-r2lam} imply:
\beq\label{eq-nh-lin}
\left[
\begin{array}{cc}
  \Gamma^3_{11} &   \Gamma^3_{12} +1\\
  \Gamma^3_{21} -1&   \Gamma^3_{22}  
\end{array}
\right]\,
\left[
\begin{array}{c}
  \aaa_1\\
  \aaa_2  
\end{array}
\right]=0
\eeq
Since  $[\aaa^1,\aaa^2]^T$ is a non-zero vector, 
 matrix $M=\left[
\begin{array}{cc}
  \Gamma^3_{11} &   \Gamma^3_{12} +1\\
  \Gamma^3_{21} -1&   \Gamma^3_{22}  
\end{array}
\right]$ must have rank less than 2  in order  for  \eq{eq-nh-lin} to have a solution, i.e.
\beq\label{nec-nh1}  \Gamma^3_{11}\, \Gamma^3_{22}- (\Gamma^3_{12} +1)(\Gamma^3_{21} -1)=0.\eeq
Substituting $c_{12}^3=1$ in \eq{nec-nh1} and simplifying, we get condition
\beq\label{nec-nh2} \Gamma_{12}^3\,\Gamma_{21}^3- (c_{12}^3)^2=\Gamma_{11}^3\,\Gamma_{22}^3.\eeq
\item[4)] At least one of the expressions  $\Gamma^3_{12} +1$ or $ \Gamma^3_{21} -1$ is non-zero (if both are zero, then $c^3_{21}=-2$, which contradicts our assumption that $c$'s are the structure constants for the frame $\{\vr_1, \vr_2, \vs=[\vr_1,\vr_2]\}$.
 Then \eq{eq-nh-lin} has a one parametric family of solutions. In part 2) of the proof, we argued that we  may assume that $\aaa^1\neq0$. Then, from \eq{eq-nh-lin}, we can express
 \beq\label{eq-a2a1} \aaa^2=\alpha\, \aaa^1,\eeq
 where $\alpha(u)$ is some \emph{known} function expressible in terms of $\Gamma$'s. (Explicitly, if $\Gamma^3_{12} \neq -1$, then $\alpha=\frac{ \Gamma^3_{11}}{\Gamma^3_{12} +1}$, otherwise, we can show that $\Gamma^3_{22}\neq 0$ and 
  $\alpha=\frac{ 3}{\Gamma^3_{22}}$.)
  
 Substitution of  \eq{eq-a2a1} into \eq{eq-a1-nh},  \eq{eq-a2-nh},  \eq{eq-r2a1-nh}, \eq{eq-sl-nh1},  and \eq{eq-r1a1-nh1}, leads to equations:
 \begin{align}
\label{eq-nh-r1lam}{\vr_1} (\lambda)&=\alpha\,\aaa^1,\\
\label{eq-nh-r2lam}\ {\vr_2} (\lambda)&=-\aaa^1,\\
\label{eq-nh-slam}{\vs} (\lambda)&=\alpha_1\,\aaa^1,\\
\label{eq-nh-r1a1}{\vr_1} (\aaa^1)&=\alpha_2\,\aaa^1,\\
\label{eq-nh-r2a1}{\vr_2} (\aaa^1)&=\alpha_3\,\aaa^1,
\end{align}
where $\alpha, \alpha_1,\alpha_2,\alpha_3$  are some known functions, expressible in terms of $\Gamma$'s and their directional derivatives. Substituting  \eq{eq-nh-r1a1} and \eq{eq-nh-r2a1}, in the commutator  relationship, we conclude that
\beq\label{eq-nh-sa1}{\vs} (\aaa^1)=\vr_1(\vr_2(\aaa^1))-\vr_2(\vr_1(\aaa^1))=\alpha_4\,\aaa^1,
 \eeq
  where $\alpha_4$  is another known function, expressible in terms of $\Gamma$'s and their directional derivatives. System \eq{eq-nh-r1lam} -- \eq{eq-nh-sa1} is a Frobenius-type system on two unknown functions, $\lam$ and $\aaa^1$, and so its solution  depends on at most two arbitrary constants. 
 \item[5)] From parts~1) and 2) of the proof, it follows that, there exists a non-constant functions $\lam$ and $\aaa^1=-\vr_2(\lam)$, satisfying  \eq{eq-nh-r1lam} -- \eq{eq-nh-sa1}, and then we immediately have a two parametric family of solution  $\lam_{c,\bar\lam}=c\,\lam +\bar\lam$, $\aaa^1_c=c\,\aaa^1$,  where $c,\bar\lam$ are arbitrary constants. From part 4), we conclude that there is no other solution. On the other hand, each $\lam_{c,\bar\lam}$, with $c\neq 0$, corresponds to a three-parametric family of  non-hyperbolic fluxes $c\,\vf+\bar\lam\bar \vf$, where $\bar\vf\in\fsid$.
 We conclude that   any non-hyperbolic flux in $\fs(\R)$ is of the form $c\,\vf+(\text{a trivial flux})$.
   \end{enumerate}
\end{proof}
\begin{remark}\label{rem-genev}From \eq{eq-a1-nh} and  \eq{eq-a2-nh} it follows that if $\vf$ is a nonhyperbolic flux for $\R=\{\vr_1,\vr_2\}$, then  $\vs=[\vr_1,\vr_2]$ is a generalized eigenvector field of $\vf$. Indeed,
\beq\label{eq-slam} \tnabla_{[\vr_1,\vr_2]}\vf=\tnabla_{\vr_1}\tnabla_{\vr_2}\vf-\tnabla_{\vr_2}\tnabla_{\vr_1}\vf=\tnabla_{\vr_1}(\lam\,\vr_2)-\tnabla_{\vr_2}(\lam\,{\vr_1})=a^1\,\vr_1+a^2\vr_2+\lam\,[\vr_1,\vr_2]. \eeq

\end{remark}
\paragraph{Proof of Theorem~\ref{thm-noninv-str}}
\begin{enumerate}
\item[1)] We want to show that  a non-zero flux $\vf\in \fs(\R)/\fstriv$ is either strictly hyperbolic or non-hyperbolic.
Assume that there exists a non-strictly hyperbolic flux $\vf\in \fs(\R)$. This means that $\vf$ has the third eigenvector field $\vr_3$  and at least two of the corresponding eigenvalue functions  $\lam^1$, $\lam^2$ and $\lam^3$ coincide in an neighborhood of a fixed point $\bar u\in \Omega$. Examining the $\vr_3$ component of the expnended flatness  condition  \eq{[r1r2]f}, we conclude that
\beq\label{eq-nsh}\Gamma^3_{12}\,\lam^2-\Gamma^3_{21}\,\lam^1=c^3_{12}\,\lam^3,\eeq
where here $c$'s and $\Gamma$'s denote structure coefficients and Christofel symbols for $\tnabla$, relative to the frame $\{\vr_1, \vr_2, \vr_3\}$. 
Equation \eq{eq-nsh} must hold as an identity in an neighborhood of $\bar u$, and it can be rewritten as
\beq\label{eq-nsh1}\Gamma^3_{12}\,(\lam^2-\lam^3) -\Gamma^3_{21}\,(\lam^1-\lam^3)\equiv 0,\eeq
From the assumption of the theorem it follows that $\Gamma^3_{12}\neq 0$, $\Gamma^3_{21}\neq 0$, and $\Gamma^3_{12}\neq \Gamma^3_{21}$. Then, from \eq{eq-nsh1}, we conclude that if any two of the functions  $\lam^1, \lam^2, \lam^3$ are equal then all there of them must be equal: $\lam^1(u)=\lam^2(u)= \lam^3(u)=\lam(u)$.
This implies that $\tnabla_{\vr}\vf=\lam\vr$ for any $r\in \mathcal X(\Omega)$. Therefore, from the flatness conditions
$$\tnabla_{[\vr_1,\vr_i]}\vf=\tnabla_{\vr_1}\tnabla_{\vr_i}\vf-\tnabla_{\vr_i}\tnabla_{\vr_1}\vf \text{ for } i=2,3$$
we can deduce that:
 $$\lam\, {[\vr_1,\vr_i]}={\vr_1}(\lam\,\vr_i)-{\vr_i}(\lam\,{\vr_1})  \text{ for } i=2,3.$$
Since the right-hand side of the above equality is $ \lam\, {[\vr_1,\vr_i]}+{\vr_1}(\lam)\,\vr_i-{\vr_i}(\lam)\,{\vr_1}$, and $\vr_1,\vr_2,\vr_3$ are independent  we conclude
$\vr_i(\lam)=0$ for $i=1,2,3$ and, therefore, $\lam\equiv \bar\lam\in\RR$ is a constant function. This implies that $\vf$ is a trivial flux, and  the statement is proven.%

\item[2)]  From Lemma~\ref{lem-non-hyp2},   if $\fs(R)$ contains strictly hyperbolic fluxes, then up to adding a trivial flux, it contains exactly one-parametric family of non-hyperbolic fluxes. Therefore,  if $\dim \fs(\R)/\fstriv>1$, then $\fs(\R)$ contains  hyperbolic fluxes, and, from the first statements of the theorem, we know that all non-trivial  hyperbolic fluxes in  $\fs(\R)$  are strictly hyperbolic.
\item[3)]  In the proof of Lemma~\ref{lem-non-hyp2}  (see \eq{nec-nh2}), we showed that if $\fs(\R)$ contains non-hyperbolic fluxes, then \eq{eq-nec-nonhyp}  holds with
 $c$'s and $\Gamma$'s being structure components and Christoffel symbols for the  connection $\tnabla$ relative to the frame $\{\vr_1,\vr_2,[\vr_1,\vr_2]\}$. Then Lemma~\ref{lem-indep} asserts  that \eq{eq-nec-nonhyp} holds with $c$'s and $\Gamma$'s corresponding to any completion   $\{\vr_1,\vr_2,\vs\}$ of $\R$ to a frame.

\end{enumerate}

\paragraph{Proof of Theorem~\ref{thm-noninv-size}}
\begin{enumerate}
\item  We want to show that $0\leq\dim \fs(\R)/\fstriv\leq 4$. Lemma~\ref{lem-lam} asserts that under  the assumptions of Theorem~\ref{thm-noninv-size},  the set of pairs of functions $\lam(\R)=\{(\lam^1,\lam^2)\}$ satisfying condition 2)~of Lemma~\ref{lem-a} is a real vector space of dimension at most 5. In addition, Lemma~\ref{lem-a} implies  for every $\lam^1$ and $\lam^2$ satisfying condition 2), there exists unique, up to adding a constant vector in $\RR^3$, flux $\vf$ satisfying \eq{FR2}. Thus $\dim\fs(\R)\leq 8$. On the other hand,  $\fs(\R)$ contains a 4-dimensional subspace of trivial fluxes and, therefore, the stated inequalities hold.
\item For $k=0, \dots, 4$, Examples~\ref{ex:0} -- \ref{ex:4a} exhibit  partial frames, satisfying the assumptions of the theorem, such that $\dim \fs(\R)/\fstriv=k$. 

\end{enumerate}
 
\section{Examples}\label{sect-ex}
The examples, provided in this section,  illustrate the main results of the paper and also provide a proof for the existence  statement  in Theorem~\ref{thm-noninv-size}. The computations were performed in the computer algebra system  {\sc Maple} by setting up systems of differential equations for $\vf$ and $\lam$'s and using a built in command {\tt pdsolove} to solve them.

\subsection{Rich partial frames}
For a rich partial frame satisfying conditions~\eqref{eq-rich-nec}, Theorem~\ref{thm-rich-sh}  describes the degree of freedom for prescribing $\lam$'s and $\vf$'s satisfying the $\fs(\R)$-system \eq{eq-fr}. The theorem also asserts that $\fs(\R)$ contains strictly hyperbolic fluxes. The following three examples demonstrate these results. They also underscore the following interesting phenomenon:  a hyperbolic flux corresponding to \emph{a rich partial frame} may have  \emph{a non-rich full frame}.   In fact, we found examples  with three different scenarios: in Example~\ref{ex:basis}, all strictly hyperbolic fluxes in $\fs(\R)$ are rich, in Example~\ref{ex1}, all hyperbolic (strictly and non-strictly)  fluxes in $\fs(\R)$ are non-rich, and finally in Example~\ref{ex2}, $\fs(\R)$ contains both rich and non-rich strictly hyperbolic fluxes.

In the following examples, $n=3$, $m=2$.    The standard affine coordinates in $\RR^3$ for the connection $\tnabla$ are denoted by $(u,v,w)$.
We start with a simple example, a partial frame given by the first two standard vectors in $\RR^3$:
\begin{example} \label{ex:basis}\rm 
  Let   $\vr_1=[1,0,0]^T$ and $\vr_2=[0,1,0]^T$ comprise a partial frame $\R$ on $\RR^3$. It is clear that $\R$ satisfies the assumptions of  Theorem~\ref{thm-rich-sh}, and as predicted by  this theorem $\lam^1$ and $\lam^2$, satisfying \eq{ic1} - \eq{ic3} are parametrized by two functions of two variables: 
  \beq\label{ex:basic-lam} \lam^1=\phi(u,w) \qquad \text{ and } \quad \lam^2=\psi(v,w).\eeq
  For each such pair of $\lam^1$ and $\lam^2$ we get a family of fluxes in $\fs(\R)$ parametrized by three arbitrary functions of one variable, $g, h $ and $k$:
   \begin{equation}\label{ex:basic-f} 
    \vf =
    \left[\begin{matrix}
  \displaystyle{\int_*^u}\phi(s,w)\,ds+g(w)\\
    \displaystyle{\int_*^v}\psi(s,w)\,ds+h(w) \\
k(w)
    \end{matrix}\right].
  \end{equation}
  On the other hand, we could start  by  parametrizing the set $\fs(\R)$ by  two arbitrary functions  $\Phi$ and  $\Psi$ of two variables and an arbitrary function $k$ of one variable: 
\beq\label{ex:basic-f-lam} \vf=[ \Phi(u,w), \Psi(v,w), k(w)] \quad \text{ with } \quad \lam^1=\pd{\Phi}{u},\quad \lam^2=\pd{\Psi}{v}.\eeq
Of course,  \eq{ex:basic-f-lam} is equivalent to \eq{ex:basic-lam} -- \eq{ex:basic-f}, but, in  \eq{ex:basic-f-lam}, arbitrary functions $g, h$ are absorbed into $\Phi$ and  $\Psi$.
Although  \eq{ex:basic-f-lam} is simpler,  \eq{ex:basic-lam} -- \eq{ex:basic-f}  more closely illustrate the argument in  the proof of Theorem~\ref{thm-rich-sh}. Obviously, for almost all choices of $\Phi$, $\Psi$ and $k$, the resulting flux is strictly hyperbolic.

We finally argue that all strictly hyperbolic fluxes in $\fs(\R)$ are rich. Let $\vr_3$ be the third eigenvector field of a  hyperbolic flux $\vf\in \fs(\R)$.
Since $\vr_3$ is linearly independent of $\vr_1$ and $\vr_2$,    it can be, up to rescaling, written as  
 $\vr_3=[a,b,1]^T$, where  $a$ and $b$ are some functions on $\RR^3$.
Then, since  $\tnabla_{\vr_3}\vr_1=\tnabla_{\vr_3}\vr_2=0$, we have, in particular, that 
\beq\label{eq:basic1}
\Gamma^2_{31}=\Gamma^1_{32}=0\quad \text{ and, therefore,} \quad c^2_{13}=\Gamma^2_{13} \text{ and }  c^1_{23}=\Gamma^2_{23}.\eeq
We also have  
  \beq\label{eq:basic2} \Gamma^3_{12}=\Gamma^3_{21}=c^3_{12}=0.\eeq
  Substituting \eq{eq:basic1} and \eq{eq:basic2} into \eq{ic2} produces two equations:
    \beq\label{eq:basic3}\Gamma^1_{23}\,(\lam^3-\lam^1)=0\quad\text{and}\quad \Gamma^2_{13}\,(\lam^3-\lam^2)=0.\eeq 
    If  $\Gamma^1_{23}\neq 0$ or $\Gamma^2_{13}\neq 0$, then \eq{eq:basic3} implies that  $\lam^3=\lam^1$ or
    $\lam^3=\lam^2$, and, therefore, $\vf$ is not strictly hyperbolic.  If $\Gamma^1_{23}= 0$ and  $\Gamma^2_{13}= 0$, then \eq{eq:basic2} implies that $c^1_{23}= 0$ and  $c^2_{13}= 0$, and therefore $\vf$ is rich.
    Thus $\fs(\R)$ does not contain non-rich strictly hyperbolic fluxes.

\end{example}


On the contrary, the following  example presents a rich pair of vector fields, satisfying \eqref{eq-rich-nec},  which admits  only  non-rich  hyperbolic
fluxes.
\begin{example}\rm \label{ex1}
Consider a partial frame $\R$, consisting  of the vector fields $\vr_1=[1,0,0]^T$ and $\vr_2=[w,1,0]^T$  on $\Omega\subset\RR^3$, such that $w\neq 0$. We have $[\vr_1,\vr_2]=0$,   $\tnabla\vr_1\,\vr_2=0$ and   $\tnabla\vr_2\,\vr_1=0$, and, therefore, we are
  in the case considered in Theorem~\ref{thm-rich-sh}. 
As predicted by Theorem~\ref{thm-rich-sh}, the degree of freedom for prescribing $\lam^1$ and $\lam^2$ consists of two arbitrary functions of two variables: 
\beq\label{ex1-lam} \lam^1=\phi\left(w, v-\frac{u} w\right) \qquad \text{ and } \quad \lam^2=\psi(v,w).\eeq
The corresponding family of fluxes is
 \begin{equation*}
    \vf =
    \left[\begin{matrix}
  w\,\displaystyle{\int_*^v}\psi(s,w)\,ds-w\,\displaystyle{\int_*^{v-\frac{u} w}}\phi(w,s)\,ds+g(w)\\
    \displaystyle{\int_*^v}\psi(s,w)\,ds+h(w) \\
k(w)
    \end{matrix}\right],
  \end{equation*}
  where $g, h $ and $k$ are arbitrary functions of one variable.
  
 Proposition~\ref{prop-suff} says that for any $\bar u\in \Omega$ and any choices  of $\phi$ and $\psi$, such that  the $\lam^1(\bar u)\neq \lam^2(\bar u)$, one can find functions  $h,g $ and $k$, so that the resulting flux is strictly hyperbolic. For a concrete example, let  $\phi\left(w, v-\frac{u} w\right)=-\frac 1 w$ and $\psi(v,w)=0$,  $g(w)=h(w)=0$ and $k(w)= -\frac 1 w - \log w$.  We observe that the flux 
   \begin{equation*}
    \vf =
    \left[\begin{matrix}
      v - \frac u w \\
      0 \\
      -\frac 1 w - \log w
    \end{matrix}\right]
  \end{equation*}
  is strictly hyperbolic with the eigenvalues 
  \begin{equation*}
    \lambda^1 = -\frac 1 w; \quad \lambda^2 = 0; \quad
    \lambda^3 = \frac{1-w}{w^2} .
  \end{equation*}
 and  with the third eigenvector given by $\vr_3=[u,0,1]^T$.

We now show that, although the partial frame $\R$ is rich, the corresponding set of fluxes $\fs(\R)$ does not contain any rich hyperbolic fluxes. Indeed, let $\vr_3$ be the third eigenvector of a strictly hyperbolic flux in $\fs(\R)$.  Up to a scaling,  any vector field, which is
  linearly independent from $\vr_1$ and $\vr_2$, is of the form
  $\vr_3=[a,b,1]^T$, where $a$ and $b$ are arbitrary functions on $\RR^3$. Since
  $[\vr_3,\vr_2]=[1,0,0]^T$, we have $c_{32}^1 = 1$, and, therefore, there is no rich hyperbolic 
fluxes in $\fs(\R)$.

\end{example}

Finally, we present an example of a rich partial frame $\R$, which admits both rich and non-rich strictly hyperbolic fluxes.
\begin{example} \label{ex2}\rm
  Consider a partial frame $\R$, consisting  of the vector fields $\vr_1=[1,-\sqrt{u},0]^T$ and $\vr_2=[1,-\sqrt{u},0]^T$  on $\Omega\subset\RR^3$, such that $u\neq 0$. One can directly check that   the assumption of Theorem~\ref{thm-rich-sh} are satisfied. 

Adjoining  the third vector field $\vr_3=[0,0,1]^T$,  we obtain a full rich  frame, which also satisfies hypothesis of Theorem~\ref{thm-rich-sh}, and therefore it admits strictly hyperbolic fluxes, all of which, by construction belong to $\fs(\R)$. We do not include the general explicit expression for these fluxes, which is rather long  and involves special functions.

  On the other hand, if we adjoin vector field $\tilde \vr_3=[1,0,-u]^T$, we obtain a non-rich full frame (with $c_{13}^2=-\frac{1}{4u}$), such that modulo $\fstriv$,  it has  a one parametric family of strictly hyperbolic  fluxes:
  \beq\label{eq-ex2-hyp}
   \vf =a\, [v, \frac{u^2}{2} + w, 0]^T,\text{ where } a\neq 0\in\RR,
  \eeq
  with the eigenvalues
  \begin{equation*}
    \lambda^1 = -\sqrt u; \quad \lambda^2 = \sqrt u; \quad
    \lambda^3 = 0 .
  \end{equation*}
 By construction, $\fs(\R)$ contains fluxes \eq{eq-ex2-hyp}, and, thus, it contains both rich and non-rich strictly hyperbolic fluxes.
\end{example}
\subsection{Non-involutive  partial frames  of two vectors fields in $\RR^3$.}
We now present examples of non-involutive partial frames  $\R=\{\vr_1,\vr_2\}$ on some open subsets of  $\RR^3$, which illustrate Theorems~\ref{thm-noninv-str} and~\ref{thm-noninv-size}. 
We continue with the examples, which satisfy all the hypothesis of   Theorem~\ref{thm-noninv-size}. These examples support the second claim of  this theorem, assuring that for each $k=0, \dots, 4$,  there exists $\R$, satisfying assumptions of the theorem, such that $\dim \fs(\R)/\fstriv=k$. 
\begin{example}[$\dim \fs(\R)/\fstriv=0$]\label{ex:0}\rm
For a partial frame $\R$ consisting of vector fields $\vr_1=[0,1,u]^T$ and $\vr_2=[w,0,1]^T$
 all fluxes are trivial.\end{example}
\begin{example}[$\dim \fs(\R)/\fstriv=1$]\label{ex:1}\rm
For a partial frame $\R$ consisting of vector fields $\vr_1=[v,u,w]^T$ and $\vr_2=[u,w,v]^T$,  on an open subset $\Omega\subset \RR^3$, where these vectors are independent,
the non-trivial fluxes form a one-parametric family:
  \[ \vf=
    \frac{c_1}{(u+v+w)^2}
   \left[ \begin{matrix}
      -\frac 1 2 u^2-uv \\
      -(u+v)(u+w) - \frac 1 2 v^2 \\
      vw+\frac 1 2 w^2
    \end{matrix}\right].
  \]
This frame does not satisfy condition \eq{eq-nec-nonhyp} and, therefore, in the  agreement with Theorem~\ref{thm-noninv-str} all non-trivial fluxes are strictly hyperbolic with eigenfunctions:
 $$
    \lambda^1 = c_1 \frac{u-v}{(u+v+w)^2}, \quad 
    \lambda^2 = c_1 \frac{v-w}{(u+v+w)^2}, \quad  \lambda^3 = 0. 
 $$
with the third  eigenvector equal to $\vr_3=[u,v,w]^T$.
\end{example}
\begin{example}[$\dim \fs(\R)/\fstriv=2$]\label{ex:2}\rm
  \newcommand{\Ei}{\operatorname{Ei}}
 For a partial frame $\R$ consisting of vector fields $\vr_1=[-1,0,v+1]^T$ and $\vr_2=[ \frac{w}{v^2-1},-1,u]^T$, defined on an appropriate open subset of $\RR^3$,
the set of  non-trivial fluxes forms a two-dimensional  vector space\footnote{technically, we should say ``the set of  non-trivial fluxes and the zero flux  form a two-dimensional  vector space.''}:
   
   \begin{align*}
  \vf &=
    c_1
   \left[ \begin{matrix}
      ((v-1)\,u+w)  \Ei(v-1) - e^{1-v} u \\
      \frac 1 2 \left[ (v-1)^2 \Ei(v-1)  - (3v+2)e^{1-v} \right] \\
      (v+1)((1-v)u-w) \Ei(v-1) + (2(v+1)u+ w)e^{1-v}
    \end{matrix}\right] \\
    &+
    c_2
    \left[\begin{matrix}
      uv+w \\
      \frac{v^2}{2} \\
      u\,(1-v^2)-vw
    \end{matrix}\right],
  \end{align*}
  where  $\Ei$ is the \emph{exponential integral}:
  \[ \Ei(x) = \int_1^\infty \frac{e^{-tx}}{t} \, dt . \]
This frame does not satisfy condition \eq{eq-nec-nonhyp}  and, therefore, in the  agreement with Theorem~\ref{thm-noninv-str} all non-trivial fluxes are strictly hyperbolic with eigenfunctions:
  \begin{align*}
    \lambda^1 &= -c_1 (2 \,\Ei(v-1) + e^{1-v}) - c_2 ; \\
    \lambda^2 &= c_1 ((v-1) \Ei(v-1)  + v\,e^{1-v}  ) + c_2 v; \\
    \lambda^3 &= c_1 e^{1-v} + c_2 .
  \end{align*}
  The third eigenvector of $[D\vf]$ is:
  $$\vr_3=
   \left[ \begin{matrix}
      c_1\,\Ei(v-1)+c_2 \\
      0 \\
       c_1\,\left(2\, e^{1-v}\,+(v-1)\Ei(v-1)\right)-c_2\,(v-1)
    \end{matrix}\right].
 $$     

\end{example}
\begin{example}[$\dim \fs(\R)/\fstriv=3$]\label{ex:3}\rm
  For a partial frame $\R$ consisting of vector fields $\vr_1=[1,\sqrt{w},0]^T$ and $\vr_2=[u,0,-w]^T$
the set of non-trivial fluxes forms a three-dimensional vector space:
  \[ \vf =
    c_1
    \left[\begin{matrix}
      3uv\sqrt{w} - v^2 - u^2w \\
      uvw \\
      vw^{3/2} - uw^2
    \end{matrix}\right]
    + c_2\,\left[
    \begin{matrix}
      v \\
      uw \\
      0
    \end{matrix}\right]
    + c_3
   \left[ \begin{matrix}
      u\sqrt{w} - v \\
      0 \\
      \frac{w^{3/2}}{3}
    \end{matrix}\right].
  \]
  In this case, when $c_1=0$ and $c_2 = \frac 1 2 c_3$, we obtain a one parametric family of non-hyperbolic fluxes:
  \[ \vf^{\rm nh} =
    c
   \left[ \begin{matrix}
      u\sqrt{w} - \frac 1 2\, v \\
     \frac 1 2\, u\,w \\
      \frac{w^{3/2}}{3}
    \end{matrix}\right].
  \]
  with the eigenfunctions:
  $$\lam^1=\lam^2=\frac 1 2 c\sqrt{w}.$$
 For an extra reassurance, we can confirm that $\R$ satisfies the  necessary condition  \eq{eq-nec-nonhyp}  for admitting  non-hyprebolic fluxes. To check this conditions we complete $\R$ to a frame, for instance,  by adjoining the vector  field $\vs=[\vr_1,\vr_2]=[1, \frac 1 2 \sqrt{w},0]^T$. One can also confirm an observation made in Remark~\ref{rem-genev} that $\vs$ is a generalized eigenvector:
 $$\tnabla_{\vs}\vf=\frac 1 2 c\sqrt{w} \,\vs+\frac 1 4 c\sqrt{w} \,\vr_1.$$

  In agreement with Theorem~\ref{thm-noninv-str},  all other fluxes are strictly hyperbolic with the eigenfunctions:
  \begin{align*}
    \lambda^1 &= c_1 (v \sqrt w + u w) + c_2 \sqrt w ; \\
    \lambda^2 &= c_1 \left(\frac 3 2 v \sqrt w - u w\right) +
      \frac{c_3}{2} \sqrt w ; \\
    \lambda^3 &= c_1 (2 v \sqrt{w} - 3 u w) - c_2 \sqrt w + c_3 \sqrt w .
  \end{align*}
The third eigenvector of $[D\vf]$ is:
  $$\vr_3= \left[ \begin{matrix}
     \,c_1^2\, (2\,u^2\,w-2\,v^2)+\,c_1\,c_2\,(3\,\sqrt{w}\,u-v)-\,c_1\,c_3\,(\sqrt{w}\,u+v)+c_2^2-c_2\,c_3\\
     \frac 1 2\,(c_1\,v+c_2)\,\left(c_1\,(\sqrt{w}\,v+u\,w)+c_2\,\sqrt{w}\right) \\
    2\, c_1\,w\,\left(c_1\,(\sqrt{w}\,v+u\,w)+c_2\,\sqrt{w}\right)
    \end{matrix}\right].
 $$

\end{example}
We observe that, in the last example, the relationship of $\vr_3$ on $c_1$, $c_2$ and $c_3$ is non-linear.
\begin{example} [$\dim \fs(\R)/\fstriv=4$]\label{ex:4a}\rm
  For a partial frame $\R$ consisting of vector fields $\vr_1=[1,0,v]^T$ and $\vr_2=[0,1,-u]^T$,
the set of non-trivial fluxes forms a four-dimensional vector space:
  \[ \vf =
    c_1
   \left[ \begin{matrix}
    2\,u\,(w+u\,v) \\
      2\,v\,(w-u\,v) \\
      w^2 +3u^2v^2
    \end{matrix}\right]
    + c_2
   \left[ \begin{matrix}
      2\,u^2 \\
      w-u\,v  \\
      2\,u^2\,v
    \end{matrix}\right]
    + c_3
   \left[ \begin{matrix}
    u\,v+w \\
      -2\,v^2 \\
      2\,uv^2
    \end{matrix}\right]
    + c_4
  \left[  \begin{matrix}
      0 \\
      2\,v \\
      w-uv
    \end{matrix}\right].
  \]
This frame does not satisfy condition \eq{eq-nec-nonhyp}  and, therefore, in the  agreement with Theorem~\ref{thm-noninv-str} all non-trivial fluxes are strictly hyperbolic with eigenfunctions:
   \begin{align*}
    \lambda^1 &= 2\,c_1 \left(  w+3  u\, v  \right) + 4\,c_2  u
      - 2\,{c_3} v  ; \\
    \lambda^2 &= 2\,c_1 \left( w-3 \, u\, v  \right) -2\, c_2 u
     -4\,c_3\, v +2\, c_4 ; \\
    \lambda^3 &= 2\,c_1 w + c_2\, u -c_3 v+c_4.
  \end{align*}
The third eigenvector of $[D\vf]$ is:
 $$\vr_3= \left[ \begin{matrix}
     -2\,c_1\,u-c_3 \\
    2\,c_1\,v+c_2 \\
  2\,c_1\,u\,v+2\,c_2\,u+2\,c_3\,v-c_4
    \end{matrix}\right].
 $$
\end{example}
For the contrast, we show another maximal-dimensional example, where $\fs(\R)$ contains non hyperbolic fluxes.
\begin{example}
 [$\dim \fs(\R)/\fstriv=4$]\label{ex:4b}\rm
  For a partial frame $\R$ consisting of vector fields $\vr_1=[1,0, 2\,v]^T$ and $\vr_2=[0,1,u]^T$,
the set of non-trivial fluxes forms a four-dimensional vector space:
  \[ \vf =
    c_1
\left[    \begin{matrix}
     u\,(uv-w)  \\
      -2\,v\,(2\,uv-w) \\
      -6\,u\,v(uv- w)-2\,w^2
    \end{matrix}\right]
    + c_2
    \left[\begin{matrix}
      {u^2} \\
       2\,(2uv-w) \\
      2\,u^2v
    \end{matrix}\right]
    + c_3
   \left[ \begin{matrix}
      w-uv \\
     2\,v^2 \\
     2\,uv^2
    \end{matrix}\right]
    + c_4
  \left[  \begin{matrix}
      0 \\
      v \\
      2uv-w
    \end{matrix}\right].
  \]
  In this case, when $c_1=c_3=c_4=0$ and $c_2 = 1$, we obtain a one parametric family of non-hyperbolic fluxes:
  \[ \vf^{\rm nh} =
    c
   \left[ \begin{matrix}
    {u^2} \\
      2\,(2uv-w) \\
      2\,u^2v
    \end{matrix}
    \right]
  \]
  with the eigenfucntions:
  $$\lam^1=\lam^2=2\,c\,u.$$
  In agreement with Theorem~\ref{thm-noninv-str},  all other fluxes are strictly hyperbolic with the eigenfunctions:
\begin{align*}
    \lambda^1 &= -c_1\,w +2\,c_2\,u+c_3\,v; \\
    \lambda^2 &= 2\,c_1\,(w-3\,u\,v)+2\,c_2\,u+4\,c_3\,v+c_4 ; \\
    \lambda^3 &= 2\,c_1\,(3\,u\,v-2\,w)+2\,c_2\,u-2\,c_3\,v-c_4 .
  \end{align*}
$$\vr_3= \left[ \begin{matrix}
     c_1\,u-c_3 \\
    c_1\,v-c_2 \\
  (5\,u\,v-3\,w)-\,c_2\,u-\,c_3\,v-c_4
    \end{matrix}\right].
 $$
\end{example}

\begin{remark} The following interesting property can be observed in Examples~\ref{ex:1} --\ref{ex:4b}, where $1\leq\dim \fs(\R)/\fstriv\leq 4$.
For basis fluxes  $\vf_1, \dots,\vf_k$  presented in these examples ($k=2,...,4$, depending on an example), the corresponding Jacobian matrices, $DF_1,\dots, DF_k$, have  the  additivity of eigenvalues property property, called the $L$-property in
  Motzkin's and Taussky's papers \cite{motzkin53,motzkin55}.  By construction, $\vr_1$ and $\vr_2$ are eigenvectors of   $DF_1,\dots, DF_k$, and, therefore, it is obvious, that if $\lam^1_1,\dots,  \lam^1_k$ are the eigenvalues for $\vr_1$ of  $DF_1,\dots, DF_k$, respectively, and $\lam^2_1,\dots,  \lam^2_k$ are the eigenvalues for $\vr_2$ of  $DF_1,\dots, DF_k$,  respectively, then for $\vf= c_1\vf_1+\dots+c_k\vf_k$, the Jacobian matrix $DF$ has the eigenvalue $\lam^1=c_1\lam^1_1+\dots+c_k\lam^1_k$  for  the eigenvector $\vr_1$ and the eigenvalue $\lam^2=c_1\lam^2_1+\dots+c_k\lam^2_k$  for   the eigenvectors $\vr_2$.
  However, it is surprising  that the third eigenvalues also ``add up''. Indeed,  is still true in all of the examples that $\lam^3=c_1\lam^3_1+\dots+c_k\lam^3_k$ is the third eigenvalue of $Df$, where 
  $\lam^3_1,\dots,  \lam^3_k$  are the third eigenvalues of  $DF_1,\dots, DF_k$, despite the fact that these matrices  have non-collinear third eigenvectors $\vr_{3,1},\dots,\vr_{3,k}$.
\end{remark}

We finish with an  example demonstrating that even when the first assumption of the Theorem~\ref{thm-noninv-size}, i.~e.~the necessary conditions   \eqref{eq-nec-sh-m2} for strict hyperbolicity, holds, the second assumption  given by the condition \eq{eq-mist} may not hold.   
\begin{example}\label{ex-no-mist}\rm  Consider a partial frame $\R$, defined on an open subset of $\RR^3$, where $w>0$, consisting of vector fields $\vr_1=[1,0,w]^T$ and $\vr_2=[0,1,-\frac 9 8\ln(w)+u ]^T$. This partial frame  satisfies the necessary condition for strict hyperbolicity. A vector field $\vs=[0,0,1]^T$ completes $\R$ to a frame, and one can easily verify that relative to this frame:
\beq \label{eq-no-mist} \Gamma^3_{22} ( u)\,\Gamma^3_{11}(u) -9\,\Gamma^3_{12}(u)\,\Gamma^3_{21}(u)\equiv 0.\eeq
In fact, this example was obtained by setting up a differential equation on the components of vector fields  $\vr_1$ and $\vr_2$,  induced by the identity \eq{eq-no-mist} and finding its particular solution.

For this partial frame the  vector space $\fs(\R)/\fstriv$ is one-dimensional:
 \[ \vf =
    c
\left[    \begin{matrix}
      \frac 1 8 e^{-u}  \\
      e^{-u}\,w \\
      e^{-u}\,w\,\left(u-\frac 9 8 \ln(w)+\frac 9 8\right)
    \end{matrix}\right].
    \]
This frame does not satisfy condition \eq{eq-nec-nonhyp} and, therefore, in the  agreement with Theorem~\ref{thm-noninv-str} all non-trivial fluxes are strictly hyperbolic with eigenfunctions:
    \begin{align*}
    \lambda^1 &=-\frac 1 8 c\, e^{-u}  ; \\
    \lambda^2 &=c\,e^{-u}\,\left(u-\frac 9 8 \ln(w)\right); \\
    \lambda^3 &= 0
  \end{align*}
and the third eigenvector field  is $\vr_3=[0,1,0]^T$.
    \end{example}

\bibliographystyle{plain}

\Addresses
\end{document}